\documentclass[10pt]{article}
\usepackage{amsmath}
\usepackage{amssymb}
\usepackage{ascmac}
\usepackage{amsthm}
\usepackage[final]{graphicx} 
\usepackage[dvips]{color}
\usepackage[dvipdfmx]{xcolor}
\usepackage{bm}
\usepackage{url}
\usepackage{tabularx}
\usepackage{fancyhdr}
\usepackage{tikz}
\usepackage{here}
\usetikzlibrary{calc}
\usetikzlibrary{positioning}
  \newcommand\C{\mathbb{C}}
  
  \newcommand\Z{\mathbb{Z}}

 \newtheorem{thm}{Theorem}[section]
 \newtheorem{prp}[thm]{Proposition}
 
 \newtheorem{lmm}[thm]{Lemma}
 \newtheorem{cor}[thm]{Corollary}
 
 \newtheorem{rmk}[thm]{\rm Remark}
 
 \newtheorem{dfn}[thm]{\rm Definition}

  \def\comment#1{ }
  \def\qed{\hfill $\square$ \bigbreak }
  \usepackage{lscape}

\newcommand*\pFqskip{8mu}
\catcode`,\active
\newcommand*\pFq{\begingroup
        \catcode`\,\active
        \def ,{\mskip\pFqskip\relax}%
        \dopFq
}
\catcode`\,12
\def\dopFq#1#2#3#4#5{%
        {}_{#1}F_{#2}\biggl(\genfrac..{0pt}{}{#3}{#4};#5\biggr)%
        \endgroup
}
\newcommand*\pfqskip{8mu}
\catcode`,\active
\newcommand*\pfq{\begingroup
        \catcode`\,\active
        \def ,{\mskip\pfqskip\relax}%
        \dopfq
}
\catcode`\,12
\def\dopfq#1#2#3#4#5{%
        {}_{#1}f_{#2}\biggl(\genfrac..{0pt}{}{#3}{#4};#5\biggr)%
        \endgroup
}

  \title{A study of a Fuchsian system of rank 8 in 3 variables and the ordinary differential equations as its restrictions}
  \author{Akihito Ebisu, Yoshishige Haraoka, Masanobu Kaneko, Hiroyuki Ochiai, \\Takeshi Sasaki and Masaaki Yoshida
}
  
\begin{document}
    \maketitle 
\begin{abstract}A Fuchsian system of rank 8 in 3 variables with 4 parameters is found. The singular locus consists of six planes and a cubic surface. The restriction of the system onto the intersection of two singular planes is an ordinary differential equation of order four with three singular points. A middle convolution of this equation turns out to be the tensor product of two Gauss hypergeometric equations, and another middle convolution sends this equation to the Dotsenko-Fateev equation. Local solutions of these ordinary differential equations are found. Their coefficients are {\it sums} of products of the Gamma functions.
These sums can be expressed as special values of the generalized hypergeometric series $_4F_3$ at 1.\end{abstract}
\par\noindent
    running title: Fuchsian differential equations

  \tableofcontents
  \newpage
\section*{Introduction}In Part I, we find a Fuchsian system $Z_3(A)$ of rank 8 in 3 variables $(t_1,t_2,t_3)$ with 4 parameters $A=(A_0,A_1,A_2,A_3)$. The singular locus consists of six planes and a cubic surface (the Fricke surface):
$$t_i=\pm1\quad (i=1,2,3), \qquad 1-t_1^2-t_2^2-t_3^2+2t_1t_2t_3=0.$$
In general, we can have ordinary differential equations as we like, we have only to give polynomials to coefficients, but to have a system in more than 2 variables of finite and non-zero rank, the coefficients must satisfy the integrability condition (a system of non-linear differential equations). Very few examples are known (most of them belong to the so-called hypergeometric family). Our system $Z_3(A)$, which is not hypergeometric, is an important example. 

The restriction of the system onto the plane, say $t_3=1$, is a system $Z_2(A)$ of rank 6 in 2 variables. Its singular locus consists of
$$t_i=\pm1,\quad (i=1,2)\qquad t_1=t_2,$$
which is equal to that of the well-known Appell's hypergeometric system $F_1$. 

In the future, we would like to study these systems; power series solutions, integral representations of the solutions, etc.

To understand these systems we begin with studying its several restrictions. The restriction of the system onto the diagonal $t_1=t_2=t_3$ is an ordinary differential equation $Z_{\Delta8}$ of order 8.

The restriction of the system $Z_2(A)$ onto the line, say,  $t_2=1$ is an ordinary differential equation $Z(A)$ of order 4 with three singular points $t_1=\pm1$ and $\infty$.  This equation $Z(A)$ has not been studied so far, to the authors knowledge.
\par\smallskip
While studying local solutions of $Z(A)$, which are fully presented in Part II, we find a power-series solution to $Z(A)$ at $t=1$, which is very similar to the product of two Gauss hypergeometric series.  This leads to the discovery that a middle convolution sends the equation $Z(A)$ to the tensor product of two Gauss hypergeometric equations, with special parameters. We also find that another middle convolution sends $Z(A)$ to the Dotsenko-Fateev equation. 
\par\medskip

In Part II, we study local solutions for the ordinary differential equation $Z(A)$ and for several related ones around their singular points.
We see the relation between $Z(A)$ and the tensor product of two specific Gauss hypergeometric equations.
At a singular point of $Z(A)$, say $t=1$, the coefficients of the holomorphic solution to $Z(A)$ satisfy a 3-term difference equation $Rc_{0}(A)$.
On the other hand, $_4F_{3}(*; 1)$, special values at the unit argument of the terminating generalized hypergeometric series $_4F_3$, satisfy a linear difference equation of order 2, if the parameters are carefully chosen.
Comparing the invariant of this difference equation with that of $Rc_{0}(A)$, we find solutions of $Rc_{0}(A)$  expressed in terms of $_4F_3(*;1)$.
From the observation that the special values $_4F_3(* ;1)$ appear as the coefficients of the product of two Gauss hypergeometric series, we notice that its product has relevance to the holomorphic solution to $Z(A)$ at $t=1$, which leads to the discovery stated above.

\par\smallskip
For most local solutions of the ordinary differential equations related to $Z(A)$, we can make use of middle convolutions connecting the equation and the tensor product of two Gauss equations to get explicit expressions for the solutions. But in these cases also, we  present a way to get them by using the difference equations for $_4F_3(*;1)$, because this method gives various expressions.
\par\smallskip
The coefficients of hypergeometric-type series are products of the Gamma functions. However for our  equation $Z(A)$ and the related ones including the Dotsenko-Fateev equation, the coefficients of local solutions are {\it sums} of products of the Gamma functions. These sums can be expressed as special values $_4F_3(*;1)$.%, which can also be expressed as the coefficients of the product of two specific Gauss hypergeometric series.
\par\smallskip
Solutions of the ordinary differential equations we studied in this paper admit Euler integral representations, which will be discussed elsewhere.

\newpage
\part{A Fuchsian system of rank 8 in 3 variables and its restrictions}
In \S1, we find a Fuchsian system $Z_3(A)$ of rank 8 in 3 variables $(t_1,t_2,t_3)$ with 4 parameters $A=(A_0,A_1,A_2,A_3)$. \par\noindent
In \S2, the restriction $Z_{\Delta8}(A)$ of the system $Z_3(A)$ onto the diagonal $t_1=t_2=t_3$ is studied.\par\noindent
In \S3, the restriction of the system $Z_3(A)$ onto the plane $t_3 = 1$ is studied and the system $Z_2(A)$ is found. The restriction of the system $Z_2(A)$ onto the diagonal $t_1 = t_2$ is also studied.\par\noindent
In \S4, the restriction of the system $Z_2(A)$ onto the line $t_2 = 1$ is studied and the system $Z(A)$ is found.
This  is an ordinary differential equation of order 4 with three singular points $t_1=\pm1$ and $\infty$.  \par\noindent
\S5 gives a relation between $Z(A)$ and the tensor product of two specific Gauss hypergeometric equations.\par\noindent
\S6 gives a relation between $Z(A)$ and the Dotsenko-Fateev equation. \par\noindent
\S7 illustrates the relation among the differential equations appeared above.\par\noindent
In \S8, a Pfaffian form of $Z_3(A)$ is presented.\par\noindent 
\S9 studies the tensor products of two Gauss hypergeometric equations.

\section{A Fuchsian system of rank 8 in 3 variables $Z_3(A)$}
We treat ideals of the ring of differential operators ${\C}[a_0,\dots, t_1,\dots, \partial/\partial t_1,\dots]$. We often call a set of generators of an ideal simply as a system, which sometimes also means the corresponding system of differential equations, after introducing an unknown, say $F, u,\dots. $
\par\medskip
In 2017, Don Zagier showed us a system generated by the differential operator
$$(1-t_1^2)\partial_{11}+2(t_3-t_1t_2)\partial_{12}+(1-t_2^2)\partial_{22}+a_{0}t_1\partial_1+a_{0}t_2\partial_2$$
and those obtained by a succession of the cyclic permutation $1\to2\to3\to1$ with a parameter $a_0$, where $\partial_1=\partial/\partial t_1, \partial_{12}=\partial^2/\partial t_1\partial t_2,\ {\rm etc}$. This system  in 3 variables  $(t_1,t_2,t_3)$ is Fuchsian of rank 8, and is highly reducible.  Hoping to have less reducible system of rank 8, we considered a bit general system with more parameters and got the following result.
\begin{thm} The system generated by the operator
 $$E_3=(1-t_1^2)\partial_{11}+2(t_3-t_1t_2)\partial_{12}+(1-t_2^2)\partial_{22}
  +a_{31}t_1\partial_1+a_{32}t_2\partial_2+a_{33}t_3\partial_3+a_{30}$$
and those obtained by a succession of the cyclic permutation $1\to2\to3\to1$ with constants $a_{ij}\ (i=1,2,3,\ j=0,1,2,3)$ is of rank $8$ if and only if 
 $$a_{11}=a_{22}=a_{33}=0,\quad a_{12}=a_{13}=a_{21}=a_{23}=a_{31}=a_{32}\ (=:a_0).$$
\end{thm} 
Set $a_1=a_{10}$, $a_2=a_{20}$,   $a_3=a_{30}$. Then the operators $E_1,E_2$ and $E_3$ are given as
$$\begin{array}{ll} &E_1=(1-t_2^2)\partial_{22}+2(t_1-t_2t_3)\partial_{23}+(1-t_3^2)\partial_{33}+a_{0}t_2\partial_2+a_{0}t_3\partial_3+a_{1},\\
 & E_2=(1-t_3^2)\partial_{33}+2(t_2-t_3t_1)\partial_{31}+(1-t_1^2)\partial_{11}+a_{0}t_3\partial_3+a_{0}t_1\partial_1+a_{2},
  \\ &E_3=(1-t_1^2)\partial_{11}+2(t_3-t_1t_2)\partial_{12}+(1-t_2^2)\partial_{22}+a_{0}t_1\partial_1+a_{0}t_2\partial_2+a_{3},\end{array}$$
with parameters $a=(a_0,a_1,a_2,a_3)$. We often use parameters $A=(A_0,A_1,A_2,A_3)$  related to $a$ by
$$a_0=2A_0-3,\ a_i=A_i^2-(A_0-1)^2\qquad i=1,2,3,$$
and name the system as $Z_3(A)$.

By using $b_i=(a_1+a_2+a_3)/2-a_i$ $(i=1,2,3)$ as parameters,   
 $F$ as unknown, and writing $F_1=\partial_1F, F_{12}=\partial_{12}F$, etc, this system can be also written as
  $$\begin{array}{ll}
 % E1&:=
&(t_1^2-1)F_{11}=(t_3-t_1t_2)F_{12}+(t_2-t_3t_1)F_{13}-(t_1-t_2t_3)F_{23}+a_0t_1F_1+b_1F,\\
 % E2&:=
&(t_2^2-1)F_{22}=(t_1-t_2t_3)F_{23}+(t_3-t_1t_2)F_{21}-(t_2-t_3t_1)F_{31}+a_0t_2F_2+b_2F,\\
 % E3&=:
&(t_3^2-1)F_{33}=(t_2-t_3t_1)F_{31}+(t_1-t_2t_3)F_{32}-(t_3-t_1t_2)F_{12}+a_0t_3F_3+b_3F.
  \end{array}$$

\begin{prp}
The system $Z_3(A)$ is Fuchsian, and the singular locus in the finite space consists of six planes and a cubic surface:
$$t_i=\pm1\quad (i=1,\ 2,\ 3), \qquad 1-t_1^2-t_2^2-t_3^2+2t_1t_2t_3=0.$$
The local exponents along the divisors are given as
\begin{eqnarray*}
  t_i=\pm 1: && 0,\  1,\  2,\  3,\  4,\  5,\  1/2\pm A_i,\\
  {\rm the\ cubic\ surface}: && 0,\  1,\  2,\  3,\  A_0,\  A_0+1,\  A_0+2,\  A_0+3,\\
  t_i=\infty: && 1-A_0\pm A_j,\  1 -A_0\pm A_k,\  2-A_0\pm A_j,\  2-A_0 \pm A_k\quad (\{i,j,k\}=\{1,2,3\}).
\end{eqnarray*}
The local exponents along a divisor are defined as those of the ordinary differential equation obtained by restricting the system onto a curve intersecting the divisor transversely at an ordinary point of the divisor.
\end{prp}
The singularities are known from the matrix 1-form $\omega$ in the next subsection.
If we restrict the system onto a generic line $t_2=$constant, $t_3=$constant, we get an ordinary differential equation of order 8 in $t:=t_1$ with polynomial coefficients:
$$(t+1)^2(t-1)^2(1-t^2-t_2^2-t_3^2+2t_2t_3t)^4P(t)\frac{d^8F}{dt^8}+\cdots =0,$$
where $P(t)$ is of degree 16, the number of apparent singular points, the local exponents at each points are $0,1,2,3,4,5,6$ and $8$. Though we omit the explicit expression of the coefficients of the ordinary equation above, we find the local exponents at the singular points as in the Proposition.
\begin{rmk}[Symmetry] The system $Z_3(A)$ is invariant under
$$(t_1,t_2,t_3)\to(\varepsilon_1t_1,\varepsilon_2t_2,\varepsilon_3t_3),\quad \varepsilon_i=\pm1,\ \varepsilon_1\varepsilon_2\varepsilon_3=1,$$
$$A_j\to -A_j\quad (j=1,2,3),$$ % and
$$(t_1,t_2,t_3,A_1,A_2,A_3)\to (t_{\sigma(1)},t_{\sigma(2)},t_{\sigma(3)},A_{\sigma(1)},A_{\sigma(2)},A_{\sigma(3)}),$$
where $\sigma$ is a permutation of $\{1,2,3\}$.
\end{rmk}

\def\diff{{\rm diff}}
\subsection{Outline of the poof of Theorem 1.1}Several integrable systems of partial differential equations with many variables are known; for example Appell-Lauricella's hypergeometric system $F_A$ in $n$ variables. The rank of $F_A$ is known to be $2^n$. The form of the equations tells immediately the rank does not exceed $2^n$. But it would be quite difficult to prove that the rank is exactly $2^n$  by manipulating the differential equations; this is proved by finding $2^n$ linearly independent hypergeometric series at a singular point.

In our case, no local solutions are known; so,
we are forced to check honestly the integrability condition. 
We transform the system $Z_3(A)$ into a Pfaffian form of size 8, and show the integrability.

Let $F$ be the unknown, $F_{ij..k}$
the partial derivative of $F$ by $t_i,t_j,\dots,t_k$, and set
\[ e= {}^{\rm tr}(F,F_{1},F_{2},F_{3},F_{12},F_{13},F_{23},DF_{123}),\qquad D:=-1 + t_1^2+ t_2^2 + t_3^2 -2t_1t_2t_3. \]
A computation shows that
the derivatives $F_{ij..k}$ can be written as
linear combinations of
  $F$, $F_1$, $F_2$, $F_3$, $F_{12}$, $F_{13}$, $F_{23}$ and $F_{123}$, and 
thus we get a Pfaffian system of the form 
\[de = \omega e,\]
where $\omega$ is an 8$\times$8-matrix 1-form given in \S \ref{88}. The integrability condition
of the system is written as
\[d\omega = \omega\wedge\omega,\]
and, by computation, we get Theorem 1.1.

 \section{Restriction of $Z_3(A)$ onto the diagonal $t_1=t_2=t_3$ }
  Let $F(t_1,t_2,t_3)$ be a solution of $Z_3(A)$.  The function $F(t,t,t)$ satisfies a Fuchsian ordinary differential equation. In this section, its singular points and the exponents  are described. Proofs are omitted.
\subsection{ $Z_{\Delta8}(A)$}
For generic parameters $A=(A_0,\dots,A_3)$ the function $F(t,t,t)$ satisfies an ordinary differential equation $Z_{\Delta8}(A)$ of order 8 with regular singular points at $-1,-1/2,1, \infty$ and apparent singular points at $-2$ and other 8 points. The local exponents are given as 
\[
\begin{array}{ll}%{|l||l|}
%\hline
t=-1 :& 0,\ 1,\ \tfrac12 \pm A_1,\  \tfrac12 \pm A_2,\  \tfrac12 \pm A_3 ,\\[2mm]
%\hline
t=-\tfrac12 :& 0,\ 1,\ 2,\ 3,\  A_0,\  A_0+1,\  A_0+2,\  A_0+3 ,\\[2mm]
%\hline
t=1 :& 0,\  2A_0,\  A_0-\frac12,\  A_0+\tfrac12,\  A_0+\tfrac32,\  A_0+\tfrac52,\  A_
0+\tfrac72,\ A_0+\tfrac92 ,\\[2mm]
%\hline
t=\infty :&  \tfrac12(3-3A_0 \pm A_1 \pm A_2 \pm A_3) ,\\[2mm]
%\hline 
t=-2 :& 0,\ 1,\ \ 3,\ 4,\ 5,\ 6,\ \ 8,\ 9 ,\\[2mm]
%\hline
t=\mbox{other 8 points} :& 0,\ 1,\ 2,\ 3,\ 4,\ 5,\ 6,\ \ 8. 
%\hline
\end{array}
\]
\normalsize
\subsection{ $Z_{\Delta6}(A)$}
If $A_3=A_2$ then $F(t,t,t)$ satisfies an ordinary differential equation $Z_{\Delta6}$ of order 6 with regular singular points at $-1,-1/2,1, \infty$ and apparent singular points at $-2$ and other 4 points. The local exponents are given as 
\[
\begin{array}{ll}%{|l||l|}
%\hline
t=-1 :& 0,\ 1,\ \tfrac12 \pm A_1,\  \tfrac12 \pm A_2 ,\\[2mm]
%\hline
t=-\tfrac12 :& 0,\ 1,\ 2,\  A_0,\  A_0+1,\  A_0+2 ,\\[2mm]
%\hline
t=1 :& 0,\  2A_0,\  A_0-\frac12,\  A_0+\tfrac12,\  A_0+\tfrac32,\  A_0+\tfrac52 ,\\[2mm]
%\hline
t=\infty :&  \tfrac12(3-3A_0 \pm A_1 \pm 2A_2),\   \tfrac12(3-3A_0 \pm A_1) 
,\\[2mm]
%\hline 
t=-2 :& 0,\ 1,\ \ 3,\ 4,\ 5,\ 6 ,\\[2mm]
%\hline
t=\mbox{other 4 points} :& 0,\ 1,\ 2,\ 3,\ 4,\ \ 6. 
%\hline
\end{array}
\]
\normalsize
\subsection{ $Z_{\Delta4}(A)$}
If $A_3=A_2=A_1$ then $F(t,t,t)$ satisfies an ordinary differential equation $Z_{\Delta4}$ of order 4 with regular singular points at $-1,-1/2,1, \infty$ and only one apparent singular point at $-2$. The local exponents are given as 
\[
\begin{array}{ll}%{|l||l|}
%\hline
t=-1 :& 0,\ 1,\ \tfrac12 \pm A_1 ,\\[2mm]
%\hline
t=-\tfrac12 :& 0,\ 1,\  A_0,\  A_0+1 ,\\[2mm]
%\hline
t=1 :& 0,\  2A_0,\  A_0-\frac12,\  A_0+\tfrac12 ,\\[2mm]
%\hline
t=\infty :&  \tfrac12(3-3A_0 \pm 3A_1),\   \tfrac12(3-3A_0 \pm A_1) ,\\[2mm]
%\hline 
t=-2 :& 0,\ 1,\ \ 3,\ 4 .
%\hline
\end{array}
\]
\normalsize
 \section{Restriction of $Z_3(A)$ onto the plane $t_3=1$ and  $Z_2(A)$}
 \subsection{Equation $Z_2(A)$}
 The restriction $Z_3(A)|_{t_3=1}$ of $Z_3(A)$ onto the plane $t_3=1$ is, by definition,  generated by the operators $P$, where
  $$P(t_1,t_2,\partial_1,\partial_2) + (t_3-1) Q,\quad \partial_i:=\partial/\partial t_i$$
 belongs to $Z_3(A)$ for some operator $Q=Q(t_1,t_2,t_3,\partial_1,\partial_2,\partial_3)$.
We find two such operators $P_1$ and $P_2$ as follows.
Since   $$E_3=(1-t_1^2)\partial_{11}+2(1-t_1t_2)\partial_{12}+(1-t_2^2)\partial_{22}+a_0(t_1\partial_1+t_2\partial_2)+a_3+2(t_3-1)\partial_{12},$$
we cut off the last term, and define $P_1$ as 
 $$P_1:=(1-t_1^2)\partial_{11}+2(1-t_1t_2)\partial_{12}+(1-t_2^2)\partial_{22}+a_0(t_1\partial_1+t_2\partial_2)+a_3.$$  
We next express $E_1$ and $E_2$ as
  $$\begin{array}{ll}
  &E_1=G_1+(t_3-1)R_1+2(t_1-t_2)\partial_{23}+a_0\partial_3,\\
  &E_2=G_2+(t_3-1)R_2+2(t_2-t_1)\partial_{13}+a_0\partial_3,\end{array}$$
  where
  $$G_1=(1-t_2^2)\partial_{22}+a_0t_2\partial_2+a_1,\quad R_1=-2t_2\partial_{23}-(1+t_3)\partial_{33}+a_0\partial_3;$$
  $G_2$ and $R_2$ are given by exchanging 1 and 2 in $G_1$ and $R_1$, respectively. Differentiate these:
  $$\begin{array}{ll}
  &E_{1,1}=G_{1,1}+(t_3-1)R_{1,1}+2\partial_{23}+2(t_1-t_2)\partial_{123}+a_0\partial_{13},\\
  &E_{2,2}=G_{2,2}+(t_3-1)R_{2,2}+2\partial_{13}+2(t_2-t_1)\partial_{123}+a_0\partial_{23},\\
  \end{array}$$
  where $E_{1,1}:=\partial_1E_1, G_{1,1}:=\partial_1G_1,$ etc, for example,
  $$G_{1,1}=(1-t_2^2)\partial_{122}+a_0t_2\partial_{12}+a_1\partial_1.$$
  We have
  $$\begin{array}{ll}
  E_1-E_2&\equiv\  G_1-G_2+2(t_1-t_2)(\partial_{23}+\partial_{13}),\\
  E_{1,1}+E_{2,2}&\equiv \ G_{1,1}+G_{2,2}+(2+a_0)(\partial_{23}+\partial_{13})\end{array}$$
  modulo $(t_3-1)$, and so
  $$2(t_1-t_2)(E_{1,1}+E_{2,2})-(2+a_0)(E_1-E_2)\equiv 2(t_1-t_2)(G_{1,1}+G_{2,2})-(2+a_0)(G_1-G_2).$$
Now we define the second operator  $P_2$ by the right hand-side of  this identity:
  $$P_2:=2(t_1-t_2)\{(1-t_2^2)\partial_{122}+a_0t_2\partial_{12}+a_1\partial_1+(1-t_1^2)\partial_{112}+a_0t_1\partial_{12}+a_2\partial_2\}$$
$$-(2+a_0)\{(1-t_2^2)\partial_{22}+a_0t_2\partial_2+a_1-(1-t_1^2)\partial_{11}-a_0t_1\partial_1-a_2\}.$$
Though we have no rigorous proof that $P_1$ and $P_2$ generate the ideal  $Z_3(a)|_{t_3=1}$, we study  
the system $Z_2(A)$ in $(t_1,t_2)$ generated by $P_1$ and $P_2$.
\begin{thm}The system $Z_2(A):=\langle P_1,P_2\rangle$ is of rank 6.
  The singular locus in $\mathbb{P}^1\times\mathbb{P}^1$ is given by
  $$ t_i=\pm1,\ \infty\ (i=1,\ 2), \quad t_1=t_2.$$
 % the finite plane is given by $$t_i=\pm1\ (i=1,\ 2),\quad t_1=t_2.$$
\end{thm}
\begin{prp}
The local exponents along the divisors above are given as
$$\begin{array}{ll}
t_1=\pm1:& 0,\ 1,\ 2,\ 3,\ \dfrac12\pm A_1,\ \\
t_2=\pm1:& 0,\ 1,\ 2,\ 3,\ \dfrac12\pm A_2,\ \\
t_1=t_2:& 0,\ 1,\ 2A_0,\ 2A_0+1,\ A_0\pm A_3,\ \\
t_i=\infty:&1-A_0\pm A_j,\ 2-A_0\pm A_j,\ 1-A_0\pm A_3\qquad (\{i,j\}=\{1,\ 2\}).\end{array}$$
\end{prp}
If we restrict the system $Z_2(A)$ further
onto a generic line $t_2=$constant, we get an ordinary differential equation of order 6 in $t:=t_1$ with polynomial coefficients:
$$(t+1)^2(t-1)^2(t-t_2)^4P(t)\frac{d^6F}{dt^6}+\cdots =0,$$
where $P(t)$ is of degree 6, the number of apparent singular points, whose local exponents are  $0,1,2,3,4$ and $6$. Though we omit the explicit expression of the coefficients of the ordinary equation above, we find the local exponents at the singular points as in the proposition.
\begin{rmk}Any set of six independent solutions defines a map from
$(t_1,t_2)$-space into the five dimensional projective space, whose image is regarded as a surface.
We remark that the operator $P_1$ implies that the second jet-space
of the surface is always degenerate; the system $Z_2(A)$ is not
general in this sense among those systems of rank 6.
\end{rmk}
\subsection{Outline of the proof of Theorem 3.1} 
Using unknown $F$, we rewrite the system in Pfaffian form relative to a frame
$$e_6={}^{\rm tr}(F,F_1,F_2,(t_1-t_2)F_{11},(t_1-t_2)F_{12},(t_1-t_2)^2F_{112}).$$
This time, 
by using $P_1=0$ and $P_2=0$, and their higher-order derivatives,
we can see that the derivatives $F_{ij..k}$, $i,j,k = 1, 2$,
can be written in terms of 
$F$, $F_{1}$, $F_{2}$, $F_{11}$, $F_{12}$ and $F_{112}$.
Thus, we get a Pfaffian form $\omega_6$ such that 
$de_6 = \omega_6 e_6$. It is a straightforward computation to
see that the integrability condition $d\omega_6 =\omega_6\wedge\omega_6$
holds.
The 6$\times$6-matrix 1-form $\omega_6$ is listed in \S \ref{66}.

\subsection{Restriction of $Z_2(A)$ onto the diagonal $t_1=t_2$}
Change the coordinates from $(t_1,t_2)$ to $(t,s)$ by $t_1=t, t_2=t+s.$ Then the operator $P_1$  becomes
$$\partial_{11}-t^2\partial_{11}+a_0t\partial_{1}+a_3
+s\left\{-s\partial_{22}-2t(\partial_{12}-\partial_{22})+a_0\partial_{2}-2t\partial_{22}\right\}.$$
Thus the restriction of $Z_2(A)$ to the diagonal $s=0$ is the ordinary differential equation
$$(1-t^2)F_{11}+a_0tF_1+a_3F=0,\qquad F_1=dF/dt.$$
The local exponents at $t=-1,1$ and $\infty$ are
$$0,\  A_0-\frac12; \quad 0,\  A_0-\frac12 \quad {\rm and} \quad 1-A_0\pm A_3,$$respectively.
%}
  \section{Restriction of $Z_2(A)$ onto the line $t_2=1$ and $Z(A)$}
  \subsection{Equation $Z(A)$}% of  $Z_2(A)$ onto the line $t_2=1$}
Express $P_1$ and $P_2$ as
  $$\begin{array}{ll}
  P_1&\equiv\ Q_1+2(1-t_1)\partial_{12}+a_0\partial_2,\\
  P_2&\equiv\  Q_2+2(t_1-1)\{a_0\partial_{12}+(1-t_1^2)\partial_{112}+a_0t_1\partial_{12}+a_2\partial_2\}-(2+a_0)a_0\partial_2\\
  &=\ Q_2+2a_0(t_1^2-1)\partial_{12}-2(t_1^2-1)(t_1-1)\partial_{112}+\{2a_2(t_1-1)-(2+a_0)a_0\}\partial_2
  \end{array}$$
  mod $(1-t_2)$, where
  $$\begin{array}{ll}
  &Q_1=(1-t_1^2)\partial_{11}+a_0t_1\partial_1+a_3,\\
  &Q_2=2(t_1-1)a_1\partial_1-(2+a_0)\{a_1-(1-t_1^2)\partial_{11}-a_0t_1\partial_1-a_2\}.\end{array}$$
  Differentiate $P_1$, and we have
  $$\begin{array}{ll}
  P_{1,1}:=\partial_1P_1&=\ Q_{1,1}-2\partial_{12}-2(t_1-1)\partial_{112}+a_0\partial_{12}\\
  &=\ Q_{1,1}+(a_0-2)\partial_{12}-2(t_1-1)\partial_{112},\end{array}$$
 where $Q_{1,1}:=\partial_1Q_1,\partial_{112}:=\partial_1\partial_{12}$.  Set
 $$P_3:=P_2-(t_1^2-1)P_{1,1}=Q_2-(t_1^2-1)Q_{1,1}+(t_1^2-1)(a_0+2)\partial_{12}+\{2a_2(t_1-1)-(2+a_0)a_0\}\partial_2,$$
  and differentiate:
  $$\begin{array}{ll}
  P_{3,1}&=\ Q_{2,1}-2t_1Q_{1,1}-(t_1^2-1)Q_{1,1,1}\\
  & +\ 2t_1(a_0+2)\partial_{12}+(t_1^2-1)(a_0+2)\partial_{112}+2a_2\partial_2+\{2a_2(t_1-1)-(2+a_0)a_0\}\partial_{12}.\end{array}$$
  By using $P_1,P_3$ and $P_{1,1}$, express $\partial_2,\partial_{12}$ and $\partial_{112}$ in terms of $Q_1,Q_{1,1}$. Substitute these expressions into $P_{3,1}$, and we get an ordinary differential operator  $Z(a)$ of order four. 
  \begin{prp} The equation $Z(a)$ is irreducible.
  \end{prp}
  Proof:\   If $Z(a)$ factors as $Z_1Z_2$, then the local exponents of $Z_2$ at $t=\pm1$ and $\infty$ are subsets of those of $Z(a)$.  Riemann relation says that the sum of the local exponents of $Z_2$ is an integer. The Riemann scheme of $Z(a)$ below shows that this can not happen if $Z_2$ is of order 1 or 3. Assume that the order of $Z_2$ is 2, and let $k$ be the number of apparent singular points. Then Riemann relation says that the sum of the local exponents of $Z_2$ is equal to $1+k$. On the other hand the Riemann scheme shows that the sum is greater than $3+k$.  \qed
  This assures that $Z(a)$ is the restriction of $Z_2(a)$ onto $t_2=1$.
  \begin{thm} The restriction $Z(a)$ of $Z_2(a)$ onto the line $t_2=1$
  is given by 
    $$Z(a):=p_0\partial^4+p_1\partial^3+p_2\partial^2+p_3\partial+p_4,$$
  where $\partial=d/dt, t=t_1$, and
  $$\begin{array}{lll}
  p_0&=&2(t+1)^2(t-1)^3,\\
  p_1&=&-4(t+1)(t-1)^2\{(2+a_0)+(a_0-2)t\},\\
  p_2&=&2(t-1)\{(a_0^2-2a_1+6a_0+2+a_2+a_3)+(3a_0^2+4a_0-4+2a_1)t\\&
  +&(a_0^2-4a_0+2-a_2-a_3)t^2\},\\
  p_3&=&(-4a_0^2-8a_0+4a_0a_1+4a_1-(2a_0+4)(a_2+a_3))\\&
  +&(-2a_0^3-6a_0^2-4a_0a_1-4a_1+4(a_2+a_3))t+2a_0(a_0+a_2+a_3)t^2,\\
  p_4&=&2a_2a_3t+(a_1-a_2-a_3)(a_0+2)^2-2a_2a_3.
    \end{array}$$
This equation has one accessory parameter; the local exponents do not change if we add a constant to $p_4$.
\end{thm}
We denote the operator $Z(a)$ with parameters $(A_0,A_1,A_2,A_3)$ by $Z(A)$, the explicit form of which will be given in \S 10.
The Riemann scheme of $Z(A)$ is given as
$$\left\{ \begin{array}{ccc}
t=1 & t=-1 & t=\infty \\
0&\tfrac12-A_1&1-A_0+A_2\\
A_0-\tfrac12&0&1-A_0-A_2\\
A_0+\tfrac12&1&1-A_0+A_3\\
2A_0 & \tfrac12+A_1 &1-A_0-A_3
\end{array} 
\right\} .$$
\begin{rmk}[Symmetry] $Z(A)$ is invariant under
$$A_j\to-A_j\quad (j=1,2,3)\quad {\rm and}\quad A_2\leftrightarrow A_3.$$
\end{rmk}
\begin{rmk}The Dotsenko-Fateev equation $(\S\ref{DoFa})$ appears as the restriction on a divisor of the Appell's equation $F_4$ in two variables $($see $\cite{Hara1})$.
\end{rmk}

\subsection{A small change $\tilde{Z}(A)$ of $Z(A)$}\label{cosmetic}
To increase symmetry, we introduce an operator  $\tilde{Z}(A)$ as
\footnote{Equivalent to changing the unknown $z$ of the equation $Z(A)$ to a new unknown $w$ by $z = (t-1)^{A_0-\tfrac12} w$}
\begin{align*}
\tilde{Z}(A) & := {\rm Ad}((t-1)^{-A_0+\frac12})Z(A)=(t-1)^{-A_0+\tfrac12} \circ Z(A) \circ 
(t-1)^{A_0-\tfrac12}.
\end{align*}
\def\pp{\tilde{p}}
(Ad stands for addition which will be recalled in \S 5.1.)
\footnote{Strictly speaking, $\tilde{Z}(A)=\frac12(t-1)^{-1}{\rm Ad}((t-1)^{-A_0+\frac12})Z(A)$.}
We further change the variable $t$, used for $Z(A)$ and $\tilde{Z}(A)$ etc,  into  the new variable
$$ x=\frac{1-t}2.$$
In $x$-coordinate, the differential operator $\tilde{Z}(A)$ changes into
\footnote{The equations $Z(A)$ and $\tilde{Z}(A)$ rewritten in the new variable $x$ will be denoted by the same notation.}
 \begin{equation}\label{tildeZ}
\tilde{Z}(A)=x^2(x-1)^2\partial^4+m_1(x)\partial^3+m_2(x)\partial^2+m_3(x)\partial+m_4(x),\quad \partial:=d/dx,
\end{equation}
where
$$
 \begin{aligned}
  m_1&=4 (x-1) x (2 x-1),\\
  m_2&=\frac{1}{4} \left(4 A_0^2 x-4 A_0^2-4 A_1^2 x-4 A_2^2
   x^2+4 A_2^2 x-4 A_3^2 x^2+4 A_3^2 x+58 x^2-58 x+9\right),\\
  m_3&=\frac{1}{2} \left(2 A_0^2-2 A_1^2-4 A_2^2 x+2 A_2^2-4
   A_3^2 x+2 A_3^2+10 x-5\right),\\
  m_4&=\left(A_2-\frac{1}{2}\right)\left(A_2+\frac{1}{2}\right)
  \left(A_3-\frac{1}{2}\right)\left(A_3+\frac{1}{2}\right).
 \end{aligned}
 $$
 The local exponents do not change if we add a constant to $m_3$; the constant term
%$$ A_0^2-A_1^2+A_2^2+A_3^2-\frac{5}{2} $$
 of $m_3$ is called the {\it accessory parameter}.
\begin{rmk}[Symmetry]\label{symmZA} $\tilde{Z}(A)$ is invariant under
$$A_j\to-A_j\quad (j=0,1,2,3)\quad {\rm and}\quad A_2\longleftrightarrow A_3$$
and 
$$(x,A_0,A_1)\longleftrightarrow (1-x,A_1,A_0).$$
\end{rmk}
 The Riemann scheme of $\tilde{Z}(A)$ is given as
$$\left\{ \begin{array}{ccc}
x=0 & x=1 & x=\infty \\
\tfrac12-A_0 & \tfrac12-A_1 & \tfrac12+A_2 \\
0 & 0 & \tfrac12-A_2 \\
1 & 1 & \tfrac12+A_3 \\
\tfrac12+A_0 & \tfrac12+A_1 & \tfrac12-A_3
\end{array} 
\right\}. 
$$

\subsection{Invariants of ordinary differential operators}
For a differential operator
$L=p_0\partial^4+p_1\partial^3+p_2\partial^2+p_3\partial+p_4$,
the operator
\[L^*=\partial^4\circ p_0-\partial^3\circ p_1+\partial^2\circ p_2-\partial\circ p_3+p_4,\qquad \partial=d/dx\]
is called the {\sl adjoint operator}.
\begin{prp}The equation $\tilde Z(A)$ is self-adjoint.
\end{prp}
To explain the meaning of this proposition, we recall some differential invariants of ordinary differential operators. 
An ordinary differential operator
 \[\partial^4 +Q_1\partial^3 +Q_2\partial^2 +Q_3\partial+Q_4 \]
 is transformed into the operator of the form
 \begin{equation}\label{orgeq}
 \partial^4 + q_2\partial^2 + q_3\partial  + q_4
\end{equation}
 which has no third-order term, by multiplying a non-zero function
 to the dependent variable. The coefficients $q_i$ are given as
\begin{eqnarray*}
&&  q_2=Q_2-\frac{3}{2}Q_1' -\frac{3}{8}Q_1^2, \\ 
&&  q_3=Q_3- \frac{1}{2}Q_1Q_2 +\frac{1}{8}Q_1^3-Q_1'', \\
  &&  q_4=Q_4-\frac{1}{4}Q_1Q_3 +\frac{1}{16}Q_1^2Q_2
  -\frac{3}{256}Q_1^4
  -\frac{1}{4}Q_2 Q_1' +\frac{3}{32}Q_1^2Q_1'
  + \frac{3}{16}(Q_1')^2- \frac{1}{4}Q_1'''. 
\end{eqnarray*}
It is known (\cite{Wilc}) that, for an appropriate choice of the dependent variable and the coordinate $y=y(x)$, the operator $(\ref{orgeq})$ can be transformed further into an operator 
\begin{equation} \label{normaleq}
  \partial^4+r_3 \partial + r_4, \qquad \partial=d/dy. 
\end{equation}
Though $r_3$ and $r_4$ are not unique, the forms
$$  \theta_3:=r_3dy^{\otimes3}=(q_3-q_2')dx^{\otimes3},\quad \theta_4:=\left(r_4-\frac12r_3'\right)dy^{\otimes4}=\left(q_4-\frac12q_3'-\frac9{100}q_2^2+\frac15q_2''\right)dt^{\otimes4}$$
are unique and are called the fundamental invariants of the operator $(\ref{orgeq})$.
   
By an easy calculation, we see that the adjoint operator of
$(\ref{orgeq})$ is
\[\partial^4 + q_2 \partial^2 + (2q_2'-q_3)\partial + q_4+q_2''-q_3'.\]
Hence, we have:
\begin{lmm} The operator $(\ref{orgeq})$ is self-adjoint if and only if
$\theta_3=0$. 
\end{lmm}  
\begin{rmk} The property that
$\theta_3\equiv 0$ is rephrased geometrically as follows:
Let $z_1$, ... , $z_4$ be linearly independent solutions of the
equation and let us consider $z=[z_1,\dots,z_4]$ as a curve
in the projective space ${\bf P}^3$. Then, we can see that,
when $\theta_3\equiv 0$, the curve formed by the tangent vectors
to this curve $z$, which lies in the $5$-dimensional projective space
of all lines in ${\bf P}^3$,
is degenerate in the sense that it lives in a $4$-dimensional
hyperplane. 
\end{rmk}  
 
\section{$\tilde{Z}(A)$ is related to the tensor product of two Gauss equations}
In \S\ref{SolveAc0}, we study local solutions of $\tilde{Z}(A)$ at $x=0$ and find that they are closely related to the product of two specific Gauss hypergeometric series. In this section we show that an addition and a middle convolution connect  $\tilde{Z}(A)$ with the tensor product of the two Gauss equations.
We begin with introducing two important operations for differential operators.
\par\medskip\noindent
Detailed study of the tensor product of two Gauss equations in general is made in the last section of Part 1.
\subsection{Definition and fundamental properties of addition and middle convolution}\label{midconv}
For a differential operator $P$ in $x$ and a function $f$ in $x$, the {\bf addition} by $f$ is defined as
$${\rm Ad}(f)P:=f\circ P\circ f^{-1},$$
which is already appeared in  \S4.2;  multiplying a non-zero function $f$
 to the dependent variable to get a new one.

For a differential operator $P$ in $x$ and a complex number $\mu$, the {\bf middle convolution} $mc_{\mu}P$ with parameter $\mu$ is defined symbolically (cf. Definition 2.3 in \cite{Hara2}) as
$$mc_{\mu}P:=\partial^{-\mu}\circ P\circ \partial^{\mu}, \quad \partial=\frac{d}{dx}.$$
Actual procedure is as follows:  Write the operator $P$ in the form
$$\sum p_{ij}x^i\partial^j.$$
Set $r=\max\{i-j,\ p_{ij}\not=0\}$, and multiply $\partial^r$ to $P$ from the left, then substitute $$x^i\partial^i=\theta(\theta-1)\cdots(\theta-i+1),\quad \theta=x\partial$$
to express $\partial^rP$ as a linear combination of $\{\theta^i\partial^j\}$. Then replace $\theta$ by $\theta-\mu$, and finally divide the operator by $\partial$ from the left as many times as possible to obtain $mc_\mu P$. 
%\par\smallskip\noindent
Fundamental properties:
%\begin{itemize}
$$mc_{\mu+\mu'}=mc_{\mu}\circ mc_{\mu'},\quad mc_{-\mu}=mc_{\mu}^{-1},\qquad
mc_\mu\theta=\theta-\mu,\quad mc_\mu\partial=\partial.$$
%\end{itemize}
\subsection{A middle convolution connects $\tilde{Z}(A)$ with the tensor product of two Gauss equations}
%Suggested by the formula in Remark \ref{discovermc} (\S\ref{SolveAc0}),
Recall the Gauss equation:
$$
E\left(\genfrac{}{}{0pt}{0}{a,b}{c};x\right)
=x(x-1)\partial^2+((a+b+1)x-c)\partial +ab.
$$
Let us  consider the tensor product
$$K=K(A)=K(A_{-+-+},A_{-++-},1-A_0,A_{----},A_{--++},1-A_0)$$
of the two Gauss equations
$$ 
E\left(\genfrac{}{}{0pt}{0}{A_{-+-+},A_{-++-}}{1-A_0};x\right)
\quad{\rm and}\quad 
E\left(\genfrac{}{}{0pt}{0}{A_{----},A_{--++}}{1-A_0};x\right)
,$$
where
$$A_{\varepsilon_0,\varepsilon_1,\varepsilon_2,\varepsilon_3}:=\frac{\varepsilon_0A_0+\varepsilon_1A_1+\varepsilon_2A_2+\varepsilon_3A_3+1}2\qquad \varepsilon_j=\pm.$$
  It is, by definition,  the differential equation satisfied by the {\it product of the solutions} of the two Gauss equations\footnote{Do not confuse this and the {\it product of two operators}}, and is given as follows (see \S 9):
$$
 K(A)=x^3(x-1)^2\partial^4+k_1(x)\partial^3+k_2(x)\partial^2+k_3(x)\partial+k_4(x)
,\qquad \partial=d/dx,$$
where
$$
 \begin{aligned}
  k_1&=(1-x) x^2 (4 A_0 x-4 A_0-10 x+5),\\
  k_2&=x (6 A_0^2 x^2-11 A_0^2 x+5 A_0^2-24 A_0 x^2+33
   A_0 x-9 A_0\\
   &\quad-A_1^2 x-A_2^2 x^2+A_2^2
   x-A_3^2 x^2+A_3^2 x+25 x^2-25 x+4),\\
  k_3&=\frac{1}{2} (-8 A_0^3 x^2+12 A_0^3 x-4 A_0^3+36
   A_0^2 x^2-39 A_0^2 x+6 A_0^2\\
   &\quad+4 A_0 A_1^2 x+4
   A_0 A_2^2 x^2-4 A_0 A_2^2 x+4 A_0 A_3^2
   x^2-4 A_0 A_3^2 x-56 A_0 x^2\\
   &\quad+42 A_0 x-2 A_0-3
   A_1^2 x-6 A_2^2 x^2+3 A_2^2 x-6 A_3^2 x^2+3
   A_3^2 x+30 x^2-15 x),\\
  k_4&=\frac{1}{2} (2 A_0^4 x-2 A_0^4-8 A_0^3 x+5 A_0^3-2
   A_0^2 A_1^2-2 A_0^2 A_2^2 x\\
   &\quad+2 A_0^2
  A_2^2-2 A_0^2 A_3^2 x+2 A_0^2 A_3^2+12
   A_0^2 x-4 A_0^2+A_0 A_1^2+4 A_0 A_2^2
   x\\
   &\quad-A_0 A_2^2+4 A_0 A_3^2 x-A_0 A_3^2-8
   A_0 x+A_0+2 A_2^2 A_3^2 x-2 A_2^2 x-2
   A_3^2 x+2 x).
 \end{aligned}
$$
To increase symmetry, we make a change (addition) as:
$$L(A)={\rm Ad}(x^{-A_0})K(A)=x^{-A_0}\circ K(A)\circ x^{A_0},$$
which can be expressed as
 $$
 L(A)=x^2(x-1)^2\partial^4+\ell_1(x)\partial^3+\ell_2(x)\partial^2+\ell_3(x)\partial+\ell_4(x),\qquad \partial=d/dx $$
where
$$
 \begin{aligned}
  \ell_1&=5 (x-1) x (2 x-1),\\
  \ell_2&=A_0^2 x-A_0^2-A_1^2 x-A_2^2 x^2+A_2^2
   x-A_3^2 x^2+A_3^2 x+25 x^2-25 x+4,\\
  \ell_3&=-\frac{3}{2} \left(-A_0^2+A_1^2+2 A_2^2 x-A_2^2+2
   A_3^2 x-A_3^2-10 x+5\right),\\
  \ell_4&=(A_2-1)(A_2+1)(A_3-1)(A_3+1).
 \end{aligned}
$$
\noindent
The constant term $\ell_{30}$ of $\ell_3(x)=\ell_{31}x+\ell_{30}$ is the accessory parameter.
The Riemann scheme of $L(A)$ is given as
\begin{align*}
&
x^{-A_0}
\left\{ \begin{array}{ccc}
x=0 & x=1 & x=\infty \\
0 & 0 & A_{-+-+} \\
A_0 & -A_1 & A_{-++-}
\end{array} 
\right\}
\times
\left\{ \begin{array}{ccc}
x=0 & x=1 & x=\infty \\
0 & 0 & A_{----} \\
A_0 & A_1 & A_{--++}
\end{array} 
\right\} \\
&=
x^{-A_0}
\left\{ \begin{array}{ccc}
x=0 & x=1 & x=\infty \\
0 & 0 & -A_0-A_2+1 \\
A_0 & -A_1 & -A_0+A_3+1 \\
A_0+1 & A_1 & -A_0-A_3+1 \\
2A_0 & 1 & -A_0+A_2+1
\end{array} 
\right\} %\\
=
\left\{ \begin{array}{ccc}
x=0 & x=1 & x=\infty \\
-A_0 & -A_1 & 1+A_2 \\
0 & 0 & 1-A_3 \\
1 & 1 & 1+A_3 \\
A_0 & A_1 & 1-A_2
\end{array} 
\right\}. 
\end{align*}
\begin{rmk}[Symmetry] $L$ is invariant under
$$A_j\to-A_j\quad (j=0,1,2,3)\quad {\rm and}\quad A_2\longleftrightarrow A_3$$
and 
$$(x,A_0,A_1)\longleftrightarrow (1-x,A_1,A_0).$$
\end{rmk}
This operator $L$ is connected with the operator $\tilde{Z}(A)$, introduced in \S 4.2, by the middle convolution as follows.
\begin{thm} $mc_{-\frac{1}{2}}\tilde{Z}(A)=L(A).$\end{thm}\label{midconvZto2G}
 We first express the operator $\tilde{Z}(A)$ defined in \S\ref{cosmetic}  (\ref{tildeZ}) as a polynomial in 
$\theta = x\partial_x$ and $\partial=\partial_x$:
\begin{eqnarray*}
\tilde{Z}(A) &=& (\theta+\frac{1}{2}-A_2)(\theta +\frac{1}{2}+A_2)
  (\theta + \frac{1}{2} - A_3)(\theta + \frac{1}{2} + A_3) \\
&& \quad -(\theta+1)(2\theta^2+4\theta+A_1^2-A_0^2-A_2^2-A_3^3+\frac{5}{2})
  \partial \\
&& \quad + (\theta +\frac{3}{2}+A_0)(\theta + \frac{3}{2}-A_0)\partial^2.
\end{eqnarray*}
Then, replacing $\theta$ by $\theta+1/2$, we have
\begin{eqnarray*}
  mc_{-1/2} \tilde{Z}(A)
  &=& (\theta+1-A_2)(\theta +1+A_2)(\theta + 1-A_3)(\theta + 1+ A_3) \\
&& \quad -(\theta+\frac{3}{2})(2\theta^2+6\theta+A_1^2-A_0^2-A_2^2-A_3^3+5)
  \partial \\
  && \quad + (\theta +2+A_0)(\theta + 2-A_0)\partial^2,
\end{eqnarray*}
which turns out to be $L(A)$. This expression leads to

\begin{cor} { Assume $A_2$ is a positive integer$:$ $A_2=m+1$ for $m\ge0$,
$($or a negative integer $A_2=-m-1$, $m\ge 0$$;$ namely, $A_2$ is a non-zero integer$)$.
Then, the equation $L(A)u=0$ has a polynomial solution of degree $m$.}
\end{cor}
\medskip
\noindent Proof. For an integer $k$, we see that
\[ L(A)x^k = p_kx^k+q_{k-1}x^{k-1}+r_{k-2}x^{k-2},\]
where
\begin{eqnarray*}
  p_k &=& \left((k+1)^2-A_2^2\right)\left((k+1)^2-A_3^2\right), \\
  q_{k-1} &=& -k(2k+1)(k^2+k-2+\alpha), \qquad \alpha = (A_1^2-A_0^2-A_2^2-A_3^2+5)/2, \\
  r_{k-2} &=& k(k-1)(k^2-A_0^2),
\end{eqnarray*}
and $q_1=r_2=r_1=0$. 
This implies that the  operator $L(A)$ sends the space of polynomials of
degree smaller than or equal to $m$, into itself.
Since $p_m=0$ by assumption, the image is a proper subspace; hence,
there exists a polynomial solution.

\bigskip

Let $u=\sum_{k=0}^m a_kx^k$ be such a solution. Then, we see
\begin{eqnarray*}
  &a_{m-1}p_{m-1}+a_mq_{m-1}=0, & \\
  &a_{m-2}p_{m-2}+a_{m-1}q_{m-2}+a_mr_{m-2}=0, & \\
  & \cdots &, \\
  &a_kp_k + a_{k+1}q_k + a_{k+2}r_k=0, & \\
  & \cdots &, \\
  & a_0p_0+a_1q_0+a_2r_0=0.&
\end{eqnarray*}
If $A_3\notin \Z$, for a given value $a_m$,
the other $a_k$ are uniquely determined as follows.
Let $M$ be a lower triple-triangular matrix defined as
\[
M =\left(
\begin{array}{cccccc}
  p_{m-1} &  \\
  q_{m-2} & p_{m-2} & \\
  r_{m-3} & q_{m-3} & p_{m-3} & \\
         & r_{m-4} & q_{m-4} & p_{m-4} &  \\
         &       &  \ddots & \ddots& \ddots & \\
         &        &         &         r_0 & q_0 & p_0 \\  
\end{array}
\right),
\]
and $a={}^{\rm tr}(a_{m-1}, a_{m-2},\dots, a_0)$ be a column vector.
Then, the linear equations above are written as
\[Ma = {}^{\rm tr}(-a_mq_{m-1}, -a_mr_{m-2},0, \dots, 0).\]
Hence, multiplying the inverse of $M$ yields the solution $a$;
note that, since $A_3\notin \Z$, the determinant of $M$ is non-vanishing. First few terms of $a$ are given as
\[ a_m:=1,\quad a_{m-1}=\frac{-m(m^2+m-2+\alpha)}{m^2-A_3^2},\]
  \[ a_{m-2}=\frac{m(m-1)\left\{(2m-1)(m^2+m-2\alpha)((m-1)^2+m-3+\alpha)
- (m^2-A_0^2)(m^2-A_3^2)\right\}}{2(2m+1)((m-1)^2-A_3^2)(m^2-A_3^2)}.\]

%\newpage
\section{Relation between $\tilde{Z}(A)$ and the Dotsenko-Fateev equation}
\subsection{The Dotsenko-Fateev equation}\label{DoFa}
The Dotsenko-Fateev operator (\cite{DF}) is an operator of order 3 defined as
$$S=S(a,b,c,g)=x^2(x-1)^2\partial^3+s_1\partial^2+s_2\partial+s_3,\quad \partial:=d/dx$$
where 
 $$\begin{array}{ll}
%   s_0& =(x-1)^2 x^2,\\[2mm]
   s_1& = -(-1+x)x(3ax+3bx+6cx+2gx-3a-3c-g),\\[2mm]
   s_2& = 2a^2x^2+4abx^2+12acx^2+3agx^2+2b^2x^2+12bcx^2+3bgx^2+12c^2x^2+8cgx^2\\[2mm]
   &+g^2x^2-4a^2x-4abx-16acx-4agx+ax^2-8bcx-2bgx+bx^2-12c^2x-8cgx+6cx^2\\[2mm]
  &-g^2x+gx^2+2a^2+4ac+ag-2ax+2c^2+cg-6cx-gx+a+c,\\[2mm]
   s_3& = -c(2a+2+2b+2c+g)(2ax+2bx+4cx+2gx-2a-2c-g+2x-1).
\end{array}$$   
\noindent
The constant term of $s_3$ is the accessory parameter.
%is $c(2a+2c+g+1)(2a+2b+2c+g+2)$ and
The Riemann scheme is 
\begin{align*}
\left\{ \begin{array}{ccc}
x=0 & x=1 & x=\infty \\
0 & 0& -2c \\
a+c+1 & b+c+1 &-a-b-2c-g-1 \\
2a+2c+g+2 & 2b+2c+g+2 & -2a-2b-2c-g-2
\end{array} 
\right\}. \end{align*}
\begin{rmk}[Symmetry] The adjoint operator of $S(a,b,c,g)$ is given by $S$ with the change:
$$(a,b,c,g)\to(-1-a,-1-b,-1-c,-g).$$
\end{rmk}

\subsection{A middle convolution and an addition
  send $\tilde{Z}(A)$ to the Dotsenko-Fateev equation}\label{Aabcg}
The equation $ \tilde{Z}(A)=x^2(x-1)^2\partial^4+\cdots$ 
has the Riemann scheme
$$
 \left\{
  \begin{matrix}
   x=0&x=1&x=\infty\\
   0&0&\frac{1}{2}-A_2\\
   1&1&\frac{1}{2}+A_2\\
   \frac{1}{2}-A_0&\frac{1}{2}-A_1&\frac{1}{2}-A_3\\
   \frac{1}{2}+A_0&\frac{1}{2}+A_1&\frac{1}{2}+A_3
  \end{matrix}
 \right\}.
$$

\begin{prp}\label{Q}A middle convolution with parameter $-\frac12-A_2$  sends $\tilde{Z}(A)$ to the equation defined by 
\footnote{It is known (cf. \cite{Hara2}) that  if we choose the parameter $\mu$ of a middle convolution as 
$$ \mu+1=\text{one\ of\ the\ local\ exponents\ of \ $\tilde{Z}$(A)\ at\ }\infty,$$
say, $ \mu=\frac{1}{2}-A_2-1=-\frac{1}{2}-A_2,$ 
the resulting equation is of order $3$; for generic parameter it is of order 4.} 
$$ Q(A)=mc_{-\frac{1}{2}-A_2}(\tilde{Z}(A))=x^2(x-1)^2\partial^3+q_1(x)\partial^2+q_2(x)\partial+q_3(x)
$$
of order 3, where
$$
 \begin{aligned}
%  q_0&=(x-1)^2 x^2,\\
  q_1&=(2A_2+3) (x-1) x (2 x-1),\\
  q_2&=A_0^2 x-A_0^2-A_1^2 x+5 A_2^2 x^2-5 A_2^2
   x+A_2^2+12 A_2x^2-12 A_2 x+2 A_2\\
   &\quad-A_3^2
   x^2+A_3^2 x+7 x^2-7 x+1,\\
   q_3&=\frac{1}{2} (2 A_2+1) \left(A_0^2-A_1^2+2 A_2^2
   x-A_2^2+4 A_2 x-2 A_2-2 A_3^2 x+A_3^2+2
   x-1\right).
 \end{aligned}
$$
The Riemann scheme of $Q(A)$ is
$$
 \left\{
  \begin{matrix}
   x=0&x=1&x=\infty\\
   0&0&1+2A_2\\
   -A_0-A_2&-A_1-A_2&1+A_2-A_3\\
   A_0-A_2&A_1-A_2&1+A_2+A_3
  \end{matrix}
 \right\}.
$$
\end{prp}
The coefficients are determined by the local exponents except $q_{30}$, where $q_3(x)=q_{31}x+q_{30}$, which is the accessory parameter.
\begin{rmk}[Symmetry]\label{symmetryQ} $Q(A)$ is invariant under
$$A_j\to-A_j\quad (j=0,1,3)\quad {\rm and}\quad (x,A_0,A_1)\longleftrightarrow (1-x,A_1,A_0).$$
Moreover, the change $(x,A_0,A_1,A_2,A_3)\to(\frac1x,A_3,A_1,A_2,A_0)$
takes $Q(A)$ into ${\rm Ad}(x^{1+2A_2})Q(A)$. 
\end{rmk}
\begin{rmk}%If $1+2A_2\not=0$,
  The symmetry of $Q(A)$ under $x\to1/x$ and its Riemann scheme  determine the accessory parameter as above, and so characterize the equation $Q(A)$.
\end{rmk}
\begin{prp}\label{QtoR}
The operator $R(A)$ defined by
$$ R(A):={\rm Ad}(x^{A_0+A_2}(x-1)^{A_1+A_2})Q(A)$$
has an expression:
$$
 R(A)=R(A_0,A_1,A_2,A_3)=x^2(x-1)^2\partial^3+r_1(x)\partial^2+r_2(x)\partial+r_3(x),
$$
where
$$
 \begin{aligned}
%  r_0&=(x-1)^2 x^2,\\
  r_1&=(1-x) x (3 A_0 x-3 A_0+3 A_1 x+2 A_2 x-A_2-6 x+3),\\
  r_2&=3 A_0^2 x^2-5 A_0^2 x+2 A_0^2+6 A_0 A_1 x^2-6
   A_0 A_1 x+4 A_0 A_2 x^2\\
   &\quad-6 A_0 A_2 x+2
   A_0 A_2-9 A_0 x^2+12 A_0 x-3 A_0+3 A_1^2
   x^2-A_1^2 x+4 A_1 A_2 x^2\\
   &\quad-2 A_1 A_2 x-9
   A_1 x^2+6 A_1 x+A_2^2 x^2-A_2^2 x-6 A_2 x^2+6
   A_2 x-A_2-A_3^2 x^2\\
   &\quad+A_3^2 x+7 x^2-7 x+1,\\
  r_3&=-\frac{1}{2} (2 A_0 x-2 A_0+2 A_1 x-2 x+1)
   (A_0+A_1+A_2-A_3-1)
   (A_0+A_1+A_2+A_3-1).
 \end{aligned}
$$
\end{prp}
The Riemann scheme of $R(A)$ is
$$
 \left\{
  \begin{matrix}
   x=0&x=1&x=\infty\\
   0&0&1-A_0-A_1\\
   2A_0&2A_1&1-A_0-A_1-A_2-A_3\\
   A_0+A_2&A_1+A_2&1-A_0-A_1-A_2+A_3
  \end{matrix}
 \right\},
$$
and the constant term of $r_3$ is the accessory parameter.

\begin{prp}\label{RtoDF}
Change the parameters $\{A_0,A_1,A_2,A_3\}$ to $\{a,b,c,g\}$ by
\begin{align*}
   A_{0}=\dfrac{2a+2c+g+2}{2},\quad
  A_{1}=\dfrac{2b+2c+g+2}{2},\quad
  A_{2}=\dfrac{-g}{2},\quad
  A_{3}=\dfrac{2a+2b+g+2}{2}.
\end{align*}
\label{lmm_relation_q_df}
Then $R(A_0,A_1,A_2,A_3)$ exactly coincides with the Dotsenko-Fateev equation $S(a,b,c,g)$.
\end{prp}
\begin{rmk}[Symmetry] $R(A)$ is invariant only under
$$A_j\to-A_j\quad (j=3)\quad {\rm and}\quad (x,A_0,A_1)\longleftrightarrow (1-x,A_1,A_0).$$
\end{rmk}

\section{Table of related differential equations}%\large
Though the  equations $Z(A), K(A)$ and $R(A)$ have origin in the Zagier system $Z_3(A)$, the hypergeometric equations and the Dotsenko-Fateev equation, respectively, the equations $\tilde{Z}(A), L(A)$ and $Q(A)$ are more accessible.  They are related as in the table below:
\medskip
$$\begin{array}{ccccll}
&Z_3(A)&&E_1,\ E_2&\\[4mm]
{\rm rest\ }(t_3=1)&\downarrow&&\downarrow&\otimes\\[4mm]
&Z_2(A)&&K(A)=E_1\otimes E_2&\\[4mm]
{\rm rest\ }(t_2=1)&\downarrow&&\downarrow&{\rm Ad}_1\\[4mm]
&Z(A)&&L(A)&\\[4mm]
{\rm Ad}_2&\downarrow& &\downarrow&mc_{\frac12}\\[4mm]
&\tilde Z(A)&\quad\sim&\tilde{Z}(A)&{\rm in}\quad x=\frac{1-t}2\\[4mm]
 &&&\downarrow& mc_{-\frac12-A_2}\\[4mm]
&&\quad&Q(A)&&\\[4mm]
& &&\downarrow&{\rm Ad}_3\\[4mm]
&&&R(A)&=\ S(a,b,c,g):DF&\\[4mm]
\end{array}$$
\par\medskip\noindent 
Here $E_1$ and $E_2$ are Gauss hypergeometric equations:
$$
 E_1
 =
 E\left(\genfrac{}{}{0pt}{0}{A_{-+-+},A_{-++-}}{1-A_0};x\right)
 ,\quad
 E_2
 =
 E\left(\genfrac{}{}{0pt}{0}{A_{----},A_{--++}}{1-A_0};x\right)
 .$$
 The additions ${\rm Ad}_j$ are given as
$${\rm Ad}_1={\rm Ad}(x^{-A_0}),\quad {\rm Ad}_2={\rm Ad}((t-1)^{\frac12-A_0}),\quad {\rm Ad}_3={\rm Ad}(x^{A_0+A_2}(x-1)^{A_1+A_2})$$
are used just for cosmetic changes.
Since a middle convolution is additive and invertible, from
$$
 \begin{aligned}
  &mc_{-\frac{1}{2}-A_2}(\tilde{Z})=Q,\quad
  {\rm Ad}_3(Q)=R,\quad
  {\rm Ad}_1(K)=L,\quad 
  mc_{\frac{1}{2}}(L)=\tilde{Z},
 \end{aligned}
$$
we have 
$$
 Q=mc_{-\frac{1}{2}-A_2}(\tilde{Z})=\left(mc_{-\frac{1}{2}-A_2}\circ mc_{\frac{1}{2}}\right)(L)
 =mc_{-A_2}(L),
$$
and
$$
  R={\rm Ad}_3(Q)
  =\left({\rm Ad}_3\circ mc_{-A_2}\right)(L)
  =\left({\rm Ad}_3\circ mc_{-A_2}\circ{\rm Ad}_1\right)(K),
$$
and conversely, 
$$
 K=\left({\rm Ad}_1^{-1}\circ mc_{A_2}\circ{\rm Ad}_3^{-1}\right)(R).
$$
Relation of the system of parameters: $a_*, \ A_*, \ A_{\pm\pm\pm\pm}$ and $(a,b,c,g)$:
\begin{align*}
 & a_0=2A_0-3,\quad a_i=A_i^2-(A_0-1)^2,\qquad i=1,2,3,\\[2mm]
 & A_{\varepsilon_0,\varepsilon_1,\varepsilon_2,\varepsilon_3}:=(\varepsilon_0A_0+\varepsilon_1A_1+\varepsilon_2A_2+\varepsilon_3A_3+1)/2,\qquad \varepsilon_j=\pm,\\[2mm]
&
  a=\frac{A_{+-++}-1}2, \quad
  b=\frac{A_{-+++}-1}2,\quad
  c=\frac{A_{+++-}-1}2,\quad
  g=-2A_{2},
  \\[2mm]
  &
  A_{0}=\dfrac{2a+2c+g+2}{2},\quad
  A_{1}=\dfrac{2b+2c+g+2}{2},\quad
  A_{2}=\dfrac{-g}{2},\quad
  A_{3}=\dfrac{2a+2b+g+2}{2}.
\end{align*}

\newpage
\section{Explicit expressions of matrix 1-forms}
\subsection{$8\times8$-matrix form $\omega=M_1dt_1+M_2dt_2+M_3dt_3$}\label{88}
In \S1.1 the system $Z_3(a)$ is transformed into the Pfaffian form $de=\omega e$ with the frame 
$$e={}^{\rm tr}(F,F_{1},F_{2},F_{3},F_{12},F_{13},F_{23},DF_{123}),\qquad D=-1 + t_1^2+ t_2^2 + t_3^2 -2t_1t_2t_3.$$ We express in this subsection the $8\times8$-matrix 1-form $\omega=M_1dt_1+M_2dt_2+M_3dt_3$. We use parameters

 $$b_1=( -a_1+a_2+a_3)/2,\  b_2= (a_1-a_2+a_3)/2,  \ b_3=(a_1+a_2-a_3)/2.$$
%\bigskip

\bigskip\noindent $M_1=$

\[\hskip-60pt
\left(
\begin{array}{cccccccc}
% 1, 1-8
0& 1&  0&   0 &   0&    0&    0&    0   \\\noalign{\smallskip}
% 2,1-4
\frac{b_1}{(t_1^2-1)}&   \frac{a_0t_1}{(t_1^2-1)}&  0&   0 &
% 2, 5-8
 \frac{-(t_1t_2-t_3)}{(t_1^2-1)}&   \frac{-(t_1t_3-t_2)}{(t_1^2-1)}&
   \frac{t_2t_3-t_1}{(t_1^2-1)}&   0   \\\noalign{\smallskip}
% 3,1-4 
% 3, 5-8
 0&   0&   0&   0 &  1&  0&  0&   0   \\
% 4,1-4 % 4,5-8
 0&   0&   0&   0 &  0&   1&   0&   0   \\
  \noalign{\smallskip}

% 5,1-4
 0&   
% {\small \begin{array}{c}  \end{array}
%   \over  } &
 \frac{\mbox{\scriptsize $p152$} }{(t_1^2-1)D} &
 \frac{\mbox{\scriptsize $p153$}  }{(t_1^2-1)D} &
 \frac{\mbox{\scriptsize $p154$}  }{(t_1^2-1)D} &
% 5, 5-8
 \frac{\mbox{\scriptsize $p155$} }{(t_1^2-1)D} &
 \frac{1}{(t_1^2-1)}&   
 \frac{p157}{(t_1^2-1)D}&   
 \frac{(t_2-t_3t_1)}{(t_1^2-1)D}\\
\noalign{\smallskip}\noalign{\smallskip} 
%

% 6, 1-4
 0&   
 \frac{\mbox{\scriptsize $p162$} }{(t_1^2-1)D} &
 \frac{\mbox{\scriptsize $p163$} }{(t_1^2-1)D} &
 \frac{\mbox{\scriptsize $p164$} }{(t_1^2-1)D} &
% 6, 5-8
 \frac{\mbox{\scriptsize $1$} }{(t_1^2-1)} &
 \frac{p166}{(t_1^2-1)D}&   
 \frac{\mbox{\scriptsize $p167$} }{(t_1^2-1)D} &
 \frac{-(t_1t_2-t_3)}{(t_1^2-1)D} \\
\noalign{\smallskip} 

% 
% 7, 1-4
 0&   0&   0&   0 &
% 7, 5-8
   0&   0&    0& \frac{1}{D} \\ \noalign{\smallskip}\noalign{\smallskip}
   % 8, 1-4
   
 \frac{m181}{(t_1^2-1)D}&   
 \frac{m182}{D}&   
 \frac{m183}{(t_1^2-1)D}&   
 \frac{m184}{(t_1^2-1)D}&
% 8, 5-8
 \frac{m185}{(t_1^2-1)}&   
 \frac{m186}{(t_1^2-1)}&   
 \frac{m187}{(t_1^2-1)D}&   
 \frac{m188}{(t_1^2-1)D}
\end{array}
\right),
\]

\begin{eqnarray*}
 p152 &=& -b_2(t_1t_2-t_3) + b_3t_1(t_1t_3-t_2),  \\
 p153 &=& b_1(t_2^2-1) - (b_1+b_3)t_1(t_2t_3-t_1), \\
 p154  &=&  b_2(t_2t_3-t_1)   - b_1t_3(t_1t_3-t_2),  \\
 p155  &=&  -a_0(t_1^2-1)(t_2t_3-t_1) + (t_1t_2-t_3)(t_1t_3-t_2),  \\
 p157  &=&  -(1+a_0)(t_2t_3-t_1)(t_1t_3-t_2),  \\   %
 p162  &=& \sigma_{23}\circ p152 = -b_3(t_1t_3-t_2) + b_2t_1(t_1t_2-t_3), \\
 p163  &=& \sigma_{23}\circ p154 =  b_3(t_2t_3-t_1) -b_1t_2(t_1t_2-t_3), \\
 p164  &=& \sigma_{23}\circ p153 =   b_1(t_3^2-1) -(b_1+b_2)t_1(t_2t_3-t_1), \\
 p166  &=& \sigma_{23}\circ p155 = p155=  -a_0(t_1^2-1)(t_2t_3-t_1) +(t_1t_2-t_3)(t_1t_3-t_2), \\ %
 p167  &=& \sigma_{23}\circ p157 = -(1+a_0)(t_2t_3-t_1)(t_1t_2-t_3),   
\end{eqnarray*}
where we use the permutation $\sigma_{23}$:
$(\sigma_{23}\circ P)(t_1,t_2,t_3,b_1,b_2,b_3) = P(t_1,t_3,t_2,b_1,b_3,b_2)$.

\def\rem#1{}

\begin{eqnarray*}
m181 &=& -b_2b_3(t_1^2-1)(t_2t_3-t_1)+b_1(b_2+b_3)(t_1t_2-t_3)(t_1t_3-t_2), \\
m182 &=& b_3(1-t_3^2)+b_2(1-t_2^2), \\              
m183 &=& a_0b_1t_2(t_1t_2-t_3)(t_1t_3-t_2)-a_0b_3(t_1t_3-t_2)(t_2t_3-t_1)
+b_1t_3(t_1^2-1)(t_2^2-1), \\
m184 &=& a_0b_1t_3(t_1t_3-t_2)(t_1t_2-t_3)-a_0b_2(t_1t_2-t_3)(t_2t_3-t_1)
+b_1t_2(t_1^2-1)(t_3^2-1), \\
m185 &=& a_0(t_2-t_3t_1)+b_3(t_2-t_3t_1)+b_2(t_2t_1^2-t_3t_1), \\              
m186 &=& a_0(t_3-t_1t_2)+b_3(t_1^2t_3-t_1t_2)+b_2(t_3-t_1t_2), \\              
m187 &=& a_0^2(t_2t_3-t_1)(t_1t_3-t_2)(t_1t_2-t_3)
-(a_0+b_2+b_3)t_1(t_2t_3-t_1)D + b_1(t_1^2-1)D\\
&&\quad +(t_1^2-1)(t_2^2-1)(t_3^2-1), \\
m188 &=& a_0(-t_1t_2^2+2t_2t_3-t_3^2t_1+t_1^3-t_1)+2(t_1^2-1)(t_2t_3-t_1).
\end{eqnarray*}
Note that $m181$, $m182$, $m187$  and $m188$ are $\sigma_{23}$-invariant,
and $\sigma_{23}\circ m183 = m184$, $\sigma_{23}\circ m185=m186$.
\newpage%\bigskip

\noindent $M_2=$

\[\hskip-60pt
 \left( \begin{array}{cccccccc}
 0&   0&   1&   0 &  0&   0&   0&   0 \\
% 2, 1-4
 0&   0&   0&   0  & 1&   0&   0&   0 \\\noalign{\smallskip}
% 3,1-4
\frac{b_2}{(t_2^2-1)}&    0&   
        \frac{a_0t_2}{(t_2^2-1)}&    0 & 
% 3, 5-8
   \frac{-(t_1t_2-t_3) }{(t_2^2-1)}&   \frac{t_1t_3-t_2}{(t_2^2-1)}&   
   \frac{-(t_2t_3-t_1)}{(t_2^2-1)}&   0 \\\noalign{\smallskip}
%4, 1-4, 5-8
   0&   0&    0&   0 &
    0&   0&   1&   0 \\
\noalign{\smallskip}
% 5, 1-4
 0&   
 \frac{\mbox{\scriptsize $p252  $}}{(t_2^2-1)D} &
\frac{\mbox{\scriptsize $p253 $}}{(t_2^2-1)D} &
\frac{\mbox{\scriptsize $p254 $}}{(t_2^2-1)D} & 
% 5, 5-8
\frac{\mbox{\scriptsize $p255$}}{(t_2^2-1)D} &
\frac{p256 }{(t_2^2-1)D}&   
\frac{1}{(t_2^2-1)}&
\frac{-(t_2t_3-t_1) }{(t_2^2-1)D}\\\noalign{\smallskip}
\noalign{\smallskip}
% 6, 1-4
 0&   0&   0&    0 & 
% 6, 5-8
 0&   0&   0& \frac{1}{D} \\\noalign{\smallskip}
% 7, 1-4
 0&   
\frac{\mbox{\scriptsize $p272$}}{(t_2^2-1)D} &
\frac{\mbox{\scriptsize $p273$}}{(t_2^2-1)D} &
\frac{\mbox{\scriptsize $p274$}}{(t_2^2-1)D} &
% 7, 5-8
\frac{1}{(t_2^2-1)} &
\frac{p276  }{ (t_2^2-1)D} &
\frac{\mbox{\scriptsize $p277$}}{(t_2^2-1)D} &
 \frac{-(t_1t_2-t_3)}{(t_2^2-1)D} \\
\noalign{\smallskip}\noalign{\smallskip}
% 8, 1-4
\frac{m281}{(t_2^2-1)D}&   
\frac{m282}{ (t_2^2-1)D}&   
\frac{m283}{ D}&   
\frac{m284}{(t_2^2-1)D} &
% 8, 5-8
\frac{m285}{(t_2^2-1)}&   
\frac{m286}{ (t_2^2-1)D}&   
\frac{m287}{(t_2^2-1)}&   
\frac{m288}{(t_2^2-1)D}
\end{array}
\right),
\]
\medskip%%%
\[
\begin{array}{lclclcl}
p252 &=& \sigma_{12}\circ p153, & & p272 &=& \sigma_{12}\circ p163, \\
p253 &=& \sigma_{12}\circ p152, & & p273 &=& \sigma_{12}\circ p162, \\
p254 &=& \sigma_{12}\circ p154, & & p274 &=& \sigma_{12}\circ p164, \\
p255 &=& \sigma_{12}\circ p155, & & p276 &=& \sigma_{12}\circ p167, \\
p256 &=& \sigma_{12}\circ p157 = p157, & & p277 &=& \sigma_{12}\circ p166 = p255,
\\
m281 &=& \sigma_{12}\circ m181, & & m285 &=& \sigma_{12}\circ m185, \\
m282 &=& \sigma_{12}\circ m183, & & m286 &=& \sigma_{12}\circ m187, \\
m283 &=& \sigma_{12}\circ m182, & & m287 &=& \sigma_{12}\circ m186, \\
m284 &=& \sigma_{12}\circ m184, & & m288 &=& \sigma_{12}\circ m188,
\end{array}
\]
where
$(\sigma_{12}\circ P)(t_1,t_2,t_3,b_1,b_2,b_3) = P(t_2,t_1,t_3,b_2,b_1,b_3)$.

\bigskip
\noindent $M_3=$

\[ \hskip-60pt
\left( \begin{array}{cccccccc}
 0&   0&   0&   1  &   0&   0&   0&   0 \\
% 2
 0&   0&   0&   0  &   0&   1&   0&   0 \\
% 3
 0&    0&   0&   0  &   0&   0&   1&   0 \\\noalign{\smallskip}
% 4, 1-4
\frac{b_3 }{(t_3^2-1)}&   0&    0&   
\frac{a_0t_3}{(t_3^2-1)} &
% 4, 5-8
\frac{t_1t_2-t_3 }{(t_3^2-1)}&   \frac{-(t_1t_3-t_2)}{(t_3^2-1)}&   
  \frac{-(t_2t_3-t_1)}{(t_3^2-1)}&    0 \\\noalign{\smallskip}
% 5,
 0&   0&   0&   0 & 
% 5, 5-8
 0&   0&    0&  \frac{1}{D} \\
\noalign{\smallskip}
% 6, 1-4
 0&   
\frac{\mbox{\scriptsize $p362$}}{(t_3^2-1)D} &
\frac{\mbox{\scriptsize $p363$}}{(t_3^2-1)D} &
\frac{\mbox{\scriptsize $p364$}}{(t_3^2-1)D} &
% 6, 5-8
\frac{ p365 }{(t_3^2-1)D}&   
\frac{\mbox{\scriptsize $p366$}}{(t_3^2-1)D} &
\frac{1}{(t_3^2-1)}&
\frac{-(t_2t_3-t_1)}{(t_3^2-1)D}\\
\noalign{\smallskip}\noalign{\smallskip}
% 7, 1-4
 0&   
 \frac{\mbox{\scriptsize $p372$}}{(t_3^2-1)D } &
\frac{\mbox{\scriptsize $p373 $}}{(t_3^2-1)D } &
\frac{\mbox{\scriptsize $p374 $}}{(t_3^2-1)D } &
% 7, 5-8
\frac{ p375 }{ (t_3^2-1)D}&   
      \frac{ 1}{(t_3^2-1)}&   
\frac{\mbox{\scriptsize $p377$}}{(t_3^2-1)D} &
\frac{-(t_1t_3-t_2) }{(t_3^2-1)D} \\
\noalign{\smallskip}\noalign{\smallskip}
% 8, 1-4
\frac{m381}{ (t_3^2-1)D}&  
\frac{m382}{(t_3^2-1)D}&   
\frac{m383}{ (t_3^2-1)D}&   
\frac{m384}{  D} &
% 8,
\frac{m385}{(t_3^2-1)D}&   
\frac{m386}{ (t_3^2-1)}&   
\frac{m387}{ (t_3^2-1)}&   
\frac{m388}{ (t_3^2-1)D}\\
\end{array} \right),
\]

\medskip

\[
\begin{array}{lclclcl}
p362 &=& \sigma_{13}\circ p164, && p372 &=& \sigma_{13}\circ p154, \\
p363 &=&  \sigma_{13}\circ p163, && p373 &=& \sigma_{13}\circ p153, \\
p364 &=&  \sigma_{13}\circ p162, && p374 &=& \sigma_{13}\circ p152, \\
p365 &=& \sigma_{13}\circ p167 = p167, && p375 &=& \sigma_{13}\circ p157, \\
p366 &=& \sigma_{13}\circ p166, && p377 &=& \sigma_{13}\circ p155 = p366, \\
m381 &=& \sigma_{13}\circ m181, && m385 &=& \sigma_{13}\circ m187, \\
m382 &=& \sigma_{13}\circ m184, && m386 &=& \sigma_{13}\circ m186, \\
m383 &=& \sigma_{13}\circ m183, && m387 &=& \sigma_{13}\circ m185, \\
m384 &=& \sigma_{13}\circ m182, && m388 &=& \sigma_{13}\circ m188.
\end{array}
\]
where $(\sigma_{13}\circ)P(t_1,t_2,t_3,b_1,b_2,b_3) = P(t_3,t_2,t_1,b_3,b_2,b_1)$.
\begin{rmk}
Though all the poles of the  entries of $M_i$ are simple, $d\omega\not=0$ .
\end{rmk}

\newpage
\subsection{$6\times6$-matrix form $\omega_6=N_1dt_1+N_2dt_2$}\label{66}
In \S3.1 the system $Z_2(A)$ is transformed into the Pfaffian form $de_6=\omega_6 e_6$ with the frame 
$$e_6={}^{\rm tr}(F,F_1,F_2,(t_1-t_2)F_{11},(t_1-t_2)F_{12},(t_1-t_2)^2F_{112}).$$ We express in this subsection the $6\times6$-matrix 1-form $\omega_6=N_1dt_1+N_2dt_2$.\par\bigskip

\noindent $N_1=$
\[ \hskip-30pt\left(
\begin{array}{cccccc}
0&   1&   0&   0&   0&   0 \\\noalign{\smallskip}
0&   0&   0&  \frac{1}{t_1-t_2}&   0&   0 \\\noalign{\smallskip}
0&   0&   0&   0&  \frac{1}{t_1-t_2} &   0 \\
\noalign{\smallskip}
\frac{-b_1(2+a_0)}{t_1^2-1} &    
\frac{n142 }{t_1^2-1} &
\frac{-(b_1+b_3)(t_1-t_2)}{t_1^2-1}& 
\frac{n144 }{(t_1^2-1)(t_1-t_2)} &
\frac{n145 }{(t_1^2-1)(t_1-t_2)} &
\frac{1-2t_1t_2+t_1^2 }{(t_1^2-1)(t_1-t_2)} \\
\noalign{\smallskip}\noalign{\smallskip}
0&   0&   0&   0& \frac{1}{t_1-t_2} & \frac{1}{t_1-t_2} \\
\noalign{\smallskip}\noalign{\smallskip}
\frac{n161 }{t_1^2-1}&
\frac{n162}{(t_1^2-1)(t_1-t_2)} &
\frac{n163 }{(t_1^2-1)(t_1-t_2)} &  
\frac{n164 }{(t_1^2-1)(t_1-t_2)} &  
\frac{n165 }{(t_1^2-1)(t_1-t_2)} & 
\frac{n166 }{(t_1^2-1)(t_1-t_2)}
\end{array}
\right)
\]

\begin{eqnarray*}
n142 &=& -a_0(t_1+t_2)-a_0^2t_1+(b_1-b_3)(t_1-t_2), \\    
n144 &=& t_1^2+2t_1t_2-3 +a_0(2t_1^2-t_1t_2-1), \\    
n145 &=& 2(t_2^2-1)-a_0(t_1^2-2t_1t_2+1), \\    
n161&=& (2+a_0)b_1-(b_1+b_2)(b_1+b_3), \\
n162 &=&-a_0(b_1+b_3)t_1(t_1-t_2)+a_0(2+a_0)t_1(t_1-t_2)
-a_0b_2(t_1^2-2t_1t_2+1)+2b_2(t_2^2-1)+2b_3(t_1-t_2)^2, \\
n163 &=& -a_0b_3t_2(t_1-t_2)+2b_3(t_1-t_2)^2+2b_1(1-2t_1t_2+t_1^2)
+a_0b_1(1-t_1t_2), \\
n164 &=& 2-2t_1^2+a_0(1-t_1^2)+2b_3t_1(t_1-t_2)+b_2(1-2t_1t_2+t_1^2)-b_1(1-t_1^2), \\
n165&=& 2-2t_2^2+a_0(3t_1^2-4t_1t_2-t_2^2+2)+a_0^2(1-t_1t_2)+b_3(t_1^2-t_2^2)+b_1(t_1^2+t_2^2-2), \\
n166 &=&  -2t_1^2+4t_1t_2-2 +a_0(t_1^2+t_1t_2-2), \\    
\end{eqnarray*}

\noindent $N_2=$

\[ \hskip-30pt\left(
\begin{array}{cccccc}
0&   0&   1&   0&   0&   0 \\\noalign{\smallskip}
0&   0&   0&   0& \frac{1}{t_1-t_2} &   0 \\ \noalign{\smallskip}
\frac{b_1+b_2}{t_2^2-1}&
\frac{a_0t_1 }{t_2^2-1}&
\frac{a_0t_2 }{t_2^2-1}&
\frac{ -(t_1^2-1)}{(t_2^2-1)(t_1-t_2)} &
\frac{ 2(1-t_1t_2)}{(t_2^2-1)(t_1-t_2)} &
0 \\
\noalign{\smallskip}\noalign{\smallskip}
0&   0&   0& \frac{ -1}{t_1-t_2} &   0&\frac{ -1}{t_1-t_2} \\
\noalign{\smallskip}\noalign{\smallskip}
\frac{b_1(2+a_0)}{t_2^2-1}&
\frac{n252}{t_2^2-1 } &
\frac{(b_1+b_3)(t_1-t_2) }{t_2^2-1} & 
\frac{ -(2+a_0)(t_1^2-1) }{(t_2^2-1)(t_1-t_2)} &
\frac{n255 }{(t_2^2-1)(t_1-t_2)}&
\frac{-(t_1^2-1)}{(t_2^2-1)(t_1-t_2)}\\
\noalign{\smallskip}\noalign{\smallskip}
\frac{n261 }{t_2^2-1 }&
\frac{n262 }{(t_2^2-1)(t_1-t_2)}&
\frac{n263 }{(t_2^2-1)(t_1-t_2)}& 
\frac{n264 }{(t_2^2-1)(t_1-t_2)} &
\frac{n265 }{(t_2^2-1)(t_1-t_2)}&
\frac{-2(t_2^2-1)+a_0(2-t_1^2-t_2^2)}{(t_2^2-1)(t_1-t_2)}
\end{array}
\right)
\]

\begin{eqnarray*}
n252 &=& (2+a_0)a_0t_1+(b_2+b_3)(t_1-t_2), \\
n255 &=& -t_2^2-2t_1t_2+3 +a_0(t_1^2-t_2^2+1-t_1t_2), \\    
n261 &=& (2+a_0)a_0b_1+(b_1+b_3)(b_1+b_2), \\
n262 &=& (2+a_0)a_0^2t_1(t_1-t_2)+a_0b_1t_1(t_1-t_2)+a_0b_2(1-2t_1t_2+t_1^2)
+a_0b_3(2t_1^2-3t_1t_2+t_2^2)+2b_2(1-t_2^2), \\
n263&=& a_0b_3t_1(t_1-t_2)-2b_1(1-t_2^2)+a_0b_1(t_2^2+t_1^2-1-t_1t_2), \\
n264&=& (2a_0+a_0^2+b_1+b_3)(1-t_1^2)-(b_2+b_3)(1-t_2^2), \\
n265 &=& a_0(t_2^2-2t_1t_2+1)+a_0^2t_1(t_1-t_2)-2b_3t_2(t_1-t_2)+2b_1(1-t_1t_2).
\end{eqnarray*}
\begin{rmk}
Though all the poles of the  entries of $N_i$ are simple, $d\omega_6\not=0$ .
\end{rmk}
Interested readers may refer to our list of data in Maple format
\footnote{ http://www.math.kobe-u.ac.jp/OpenXM/Math/Z12data/z12.html}
for the  matrix $1$-forms in this section as well as related ordinary
differential equations in Sections 1-3. 
\newpage
\section{Tensor product of two Gauss equations}
Consider two differential equations
\[ z_1'' = S_1z_1\qquad {\rm and}\qquad z_2''=S_2z_2\quad (z':=dz/dx),\]
with dependent variables $z_1$ and $z_2$. If $S_1\not=S_2$, the product $w=z_1z_2$ satisfies the fourth-order differential equation $K_{S_1,S_2}w=0$, where
\[ K_{S_1,S_2}:= \partial^4+ f_3 \partial^3 + f_2 \partial^2 + f_1 \partial + f_0,\quad (\partial:=d/dx)\]
and
\begin{eqnarray*}
  f_3&:=& -\frac{S_1'-S_2' }{S_1-S_2}, \quad f_2:= -2(S_1+S_2), \\
  f_1&:=&-\frac{S_1S_1'-S_2S_2'+5(S_1S_2'-S_1'S_2)}{S_1-S_2}, \\
  f_0&:=& -S_1''-S_2''+(S_1-S_2)^2+\frac{(S_1')^2-(S_2')^2}{S_1-S_2}.
\end{eqnarray*}
If  $S_1=S_2=S$, $w=z_1^2$ satisfies the third-order equation $K_Sw=0$, where
\begin{equation}\label{K_S}K_S=\partial^3-4S\partial-2S'.\end{equation}
For two Gauss equations $$y_j''+p_jy_j'+q_jy_j=0,\quad p_j=\frac{c_j-(a_j+b_j+1)x}{x(1-x)},\quad q_j=-\frac{a_j b_j}{x(1-x)},\quad (j=1,2)$$
we let
\[ y_j=\lambda_j z_j, \qquad \lambda_j=x^{-c_j/2} (x-1)^{(c_j-a_j-b_j-1)/2}.\]
Then, $z_j$ satisfies the equation
\[ z_j'' =S_j z_j,\qquad S_j=-q_j+\frac{1}{4}p_j^2 + \frac{1}{2}p_j',\]
with the Riemann scheme
$$
\left\{\begin{array}{ccc}x=0&x=1&x=\infty\\
c_j/2&(a_j+b_j-c_j+1)/2&(a_j-b_j-1)/2\\
1-c_j/2&(c_j-a_j-b_j+1)/2&(b_j-a_j-1)/2\end{array}\right\}.$$
From the equation $K_{S_1,S_2}$ satisfied by $z_1z_2$, we get the differential equation
$$K=K(a_1,b_1,c_1,a_2,b_2,c_2)={\rm Ad}(\lambda_1\lambda_2)K_{S_1,S_2}$$
satisfied by $y_1y_2 = (\lambda_1\lambda_2) z_1z_2$. Though we omit the explicit form of $K$, if $(a_1,b_1,c_1)\not=(a_2,b_2,c_2)$,
it is of order four and has {\it generically} two apparent singular points say $\{x_1,x_2\}$ other than $\{0,1,\infty\}$, and the Riemann scheme is given as
$$
\left\{\begin{array}{ccc}x=0&x=1&x=\infty\\
0&0&a_1\\
1-c_1&c_1-a_1-b_1&b_1\end{array}\right\}\times
\left\{\begin{array}{ccc}x=0&x=1&x=\infty\\
0&0&a_2\\
1-c_2&c_2-a_2-b_2&b_2\end{array}\right\}$$$$=
\left\{\begin{array}{ccccc}x=0&x=1&x=\infty&x=x_1&x=x_2\\
0&0&a_1+a_2&0&0\\
1-c_1&c_1-a_1-b_1&a_1+b_2&1&1\\
1-c_2&c_2-a_2-b_2&b_1+a_2&2&2\\
2-c_1-c_2&c_1+c_2-a_1-b_1-a_2-b_2&b_1+b_2&4&4\end{array}
\right\}.$$
Note that though the Gauss equations have no accessory parameters, $K$ has one.
The two apparent singular points $\{x_1,x_2\}$ are the roots of the following quadratic form:
$$
 \begin{aligned}
 &App(x):=(a_1+a_2-b_1-b_2)(a_1-a_2-b_1+b_2)x^2\\
 &\qquad+2\bigl(2a_1b_1-2a_2b_2+(1-a_1-b_1)c_1-(1-a_2-b_2)c_2\bigr)x+(c_1-c_2)(c_1+c_2-2).
 \end{aligned}
 $$
 If for example, $c_1=c_2$, then $App$ is divisible by $x$, and so we set $x_2=0$, and if the  other parameters remain generic, the Riemann scheme becomes
$$\left\{\begin{array}{cccc}x=0&x=1&x=\infty&x=x_1\\
0&0&a_1+a_2&0\\
1-c_1&c_1-a_1-b_1&a_1+b_2&1\\
2-c_1&c_1-a_2-b_2&b_1+a_2&2\\
2-2c_1&2c_1-a_1-b_1-a_2-b_2&b_1+b_2&4\end{array}
\right\}.$$ 
 \subsection{Tensor product without apparent singularities} 
 There are several choices of parameters that the tensor product has no apparent singularities, that is the cases $App$ reduces to constant times
 $$x^2,\quad (x-1)^2,\quad 1,\quad{\rm and}\quad x-1,\quad x,\quad x(x-1),$$
 corresponding to $$\{x_1,x_2\}\to0,1,\infty, \quad{\rm and}\quad \{x_1,x_2\}\to\{0,\infty\},\{1,\infty\},\{0,1\},$$ respectively.
Thanks to the symmetry of the Gauss equations on the three singular points, we consider only two cases: $x^2$ and $x(x-1)$. The first case occurs only when 
$$\begin{array}{ll}
(1.1)\ & \{ c_1 = c_2,\ a_1 = (2a_2b_2-a_2c_2+b_1c_2-b_2c_2)/(2b_1-c_2)\},\quad{\rm or}\\ [2mm]
  (1.2)\ &\{ c_1 = 2-c_2,\ a_1 = (2a_2b_2-a_2c_2-b_1c_2-b_2c_2+2b_1+2c_2-2)/(2b_1+c_2-2)\};\end{array}$$
  and the second case, 
$$\begin{array}{ll}
(2.1)\ &\{ c_1 = c_2,\ a_1 = -b_1+2c_2-a_2-b_2\},\quad{\rm or\quad(this\ is\ used\ below)}\qquad\qquad\qquad\quad\\[2mm]
(2.2)\ &\{ c_1 = c_2,\ a_1 = a_2-b_1+b_2\},\quad{\rm or}\\[2mm]
(2.3)\ &\{ c_1 = 2-c_2,\ a_1 = -b_1-2c_2+2+a_2+b_2\},\quad{\rm or}\\[2mm]
(2.4)\ &\{ c_1 = 2-c_2,\ a_1 = -a_2-b_1-b_2+2\}.
\end{array}$$
  For the first case, two of the local exponents at $x=0$ differ by 2, and for the second case,  two of the local exponents at $x=0$ and at $x=1$ differ by 1.
  Thanks to the adjoint symmetry (cf. \S 4.3)
  \[  G(\alpha,\beta,\gamma) \longleftrightarrow G(1-\alpha, 1-\beta, 2-\gamma)\]
  of the Gauss equation\footnote{Adjoint equation of  $G(\alpha,\beta,\gamma)$ is $G(1-\alpha, 1-\beta, 2-\gamma)$.}, we study only two cases (1.1) and (2.1); in these cases we have
  $$c_1=c_2=:c.$$
\noindent
For each case, in the following, we consider the renormalized equation (addition by $x^{c-1}$) 
\[L:= {\rm Ad}(x^{c-1})(K) =x^{c-1}\circ K\circ x^{1-c}.\]
After cancelling the common factor (denoted also by $L$), it is of the form
$$L=\left\{\begin{array}{ll}
&x^3(x-1)^3\partial^4+\cdots \qquad{\rm in\ case\ }(1.1),\\[2mm]
   &x^2(x-1)^2\partial^4+\cdots \qquad{\rm in\ case\ }(2.1).\end{array}\right.$$
\subsection{Tensor product without apparent singularities Case 1}
In this subsection we assume
\[a_1 = (2a_2b_2-a_2c_2+b_1c_2-b_2c_2)/(2b_1-c_2),\quad c_1 = c_2,\]
and study the middle convolution of $L$.
The local exponents of $L$ are given as follows:
\begin{eqnarray*}
&x=0:& [0, 2, c_2-1, -c_2+1], \\
&x=1:& [0, c_1-a_1-b_1, c_2-a_2-b_2, c_1-a_1-b_1+c_2-a_2-b_2],  \\
&x=\infty:& [a_1+a_2-c_2+1, a_1+b_2-c_2+1, b_1+b_2-c_2+1, b_1+a_2-c_2+1],
\end{eqnarray*}
where 
$a_1 = (2a_2b_2-a_2c_2+b_1c_2-b_2c_2)/(2b_1-c_2)$ and $c_1 = c_2$
should be assumed.
We follow the recipe of making the middle convolution: we consider $L6:=\partial^2\circ L$, and express it in terms of 
$(x\partial_x)^i\circ (\partial_x)^j$ with constant coefficients
and replace $x\partial_x$ by $x\partial_x - m$, where $m$ is a constant (parameter of middle convolution).
The resulting operator $M6$ is of order 6, with parameter $m$, written as
\[ M6= cm_6\partial^6 + cm_5\partial^5 + cm_4\partial^4
+ cm_3\partial^3+ cm_2\partial^2+ cm_1\partial^1+ cm_0,\]
where 
\begin{eqnarray*}
cm_6 &=& (2b_1-c_2)^2x^3(x-1)^3, \\
cm_5 &=& x^2(x-1)^2(2b_1 -c_2)(4xb_1^2-8xb_1c_2+4xa_2b_1+44xb_1\\
&&\quad -12xb_1m+4xb_2b_1+4c_2^2x+4a_2b_2x-4c_2b_2x-4c_2a_2x-22c_2x\\
&&\quad +6xmc_2-20b_1+6b_1m+10c_2-3mc_2), \\
cm_4 &=& x(x-1)P_2(x),\quad cm_3=P_3(x),\quad cm_2=P_2(x),\quad cm_1= (m-1)P_1(x), \\
cm_0 &=& (m-1)(m-2)(b_1+b_2-c_2+1-m)(b_1+a_2-c_2+1-m)\\
&&\quad  \times (-2b_1m+mc_2+2b_1-c_2-b_1c_2+c_2^2+2a_2b_1-2c_2a_2+2a_2b_2-c_2b_2)\\
&&\quad \times(-2b_1m+mc_2+2b_1-c_2-b_1c_2+c_2^2-c_2a_2+2b_2b_1-2c_2b_2+2a_2b_2),
\end{eqnarray*}
where $P_k(x)$ denotes symbolically a polynomial of degree $k$ in $x$.
The local exponents of $\partial^2\circ L$ are 
\begin{eqnarray*}
&x=0:&[0,\ 1,\ 2,\ 2,\ c_2-1,\ -c_2+1], \\
&x=1:&[0, 1, 2, c_2-a_2-b_2, c_1-a_1-b_1, c_1-a_1-b_1+c_2-a_2-b_2], \\
&x=\infty:& [1, 2, 1+b_1+b_2-c_2, 1+b_1+a_2-c_2, 1+a_1+a_2-c_2, 1+a_1+b_2-c_2],
\end{eqnarray*}
and those of $M6$ are 
\begin{eqnarray*}
&x=0:&[0,\ 1,\ 2,\ m+2,\ c_2-1+m,\ 1-c_2+m], \\
&x=1:& [0,\ 1,\ 2,\ m+c_2-a_2-b_2, m+c_1-a_1-b_1, m+c_1-a_1-b_1+c_2-a_2-b_2], \\
  &x=\infty:& [-m+1,\ -m+2,\ -m+1+b_1+b_2-c_2,\ -m+1+b_1+a_2-c_2,\\
    &&\quad -m+1+a_1+a_2-c_2,\ -m+1+a_1+b_2-c_2].
\end{eqnarray*}
The difference is
\[ [0, 0, 0, m, m, m], \quad  [0, 0, 0, m, m, m],\quad 
 [-m, -m, -m, -m, -m, -m].\]
\par\medskip
 Though $L$ would be generically irreducible, since the local exponents at $x=0$ are $[0,2,c_2-1,-c_2+1]$, the resulting $M6$ might be reducible (cf. \cite{Hara2}). In fact $M6$ breaks as
\[ M6= Y\circ M5,\quad M5=u_5\partial^5+u_4\partial^4+u_3\partial^3+u_2\partial^2+u_1\partial+u_0, \]
where $Y$ is a first order operator, and $u_0,\dots, u_5$ are polynomials in $x$ and the parameters, and $u_5$ is given as
\[
 u_5 = (2b_1-c_2)^2x^2(x-1)^3(b_1-c_2+b_2)(b_1-c_2+a_2)(x+(m-1)\lambda),
\]
where
$$\lambda =-\frac{2b_1-c_2}{(b_1-c_2+b_2)(b_1-c_2+a_2)}.$$
This shows that $M5$ has one apparent singularity at
$p:=-(m-1)\lambda.$
The local exponents of $M5$ are
$$\begin{array}{ll}
x=0:&  [0,\ 1,\ 2,\ c_2-1+m,\ -c_2+1+m], \\
x=1:&  [0, 1, -b_2+c_2+m-a_2, c_1-a_1-b_1+m, c_1-a_1-b_1+c_2-a_2-b_2+m], \\
x=\infty:&  [1-m, b_1+b_2-c_2+1-m, b_1+a_2-c_2+1-m, a_1+b_2-c_2+1-m, a_1+a_2-c_2+1-m], \\
x=p:&  [0,\ 1,\ 2,\ 3,\ 5].\end{array}$$
\subsubsection{Why $M6$ is divisible from the left by a first-order operator}
Note first that the operator $L=x^3(x-1)^3\partial^4+\cdots$ can be written as
$$L=x^2Q_1(\theta,\partial)+\{\lambda(\theta-1)+x\}Q_2(\theta,\partial),\quad \theta=x\partial.$$
Since
$$\partial^2x^2=(\theta+1)(\theta+2),\quad \partial^2\{\lambda(\theta-1)+x\}=\{(\theta+1)(\lambda \partial+1)+1\}\partial,$$
$L6=\partial^2L$ is written in terms of $\theta$ and $\partial$:
$$L6=(\theta+1)(\theta+2)Q_1(\theta,\partial)+\{(\theta+1)(\lambda \partial+1)+1\}\partial Q_2(\theta,\partial).$$
So $M6$ can be obtained from $L6$ by replacing $\theta$ by $\theta-m$:
$$M6=(\theta+1-m)(\theta+2-m)Q_1(\theta-m,\partial)+\{(\theta+1-m)(\lambda \partial+1)+1\}\partial Q_2(\theta-m,\partial).$$On the other hand we have the following formulae:
$$\begin{array}{ll}
  (\theta+1-m)(\theta+2-m)&=\displaystyle\left(\theta+1-m+\frac{x}{x+u}\right)\left(\theta+2-m-\frac{x}{x+u}\right)\quad{\rm \ for\ any\ constant\ }u,\\
  (\theta+1-m)(\lambda\partial+1)+1&=\displaystyle\left(\theta+1-m+\frac{x}{x+(m-1)\lambda}\right)\left(\lambda\partial+1-\frac{\lambda}{x+(m-1)\lambda}\right).\end{array}$$
Applying these by putting $u=(m-1)\lambda$, we see that $M6$ can be divisible from the left by 
$$Y=\theta+1-m+\frac{x}{x+(m-1)\lambda}.$$
\subsubsection{When $M5$ has no apparent singular point}
Further to forget the singularity $p$, we assume, for example,
$b_1=c_2-b_2$. Then we also have $a_1=c_2-a_2$,
namely, the two Gauss are equal up to a Euler
transformation: $F(a,b,c;x)=(1-x)^{c-a-b}F(c-a,c-b,c;x)$.
But $M5$ remains to be an equation of order 5 (see \S \ref{KS} for an analogous phenomenon), with the local exponents
$$\begin{array}{ll}%\begin{verbatim}
x=0:&  [0,\ 1,\ 2,\ c_2-1+m, -c_2+1+m], \\
x=1:&  [0,\ 1,\ m,\ b_2+a_2-c_2+m, c_2-b_2+m-a_2], \\
x=\infty:&    [2-m,\ 1-m,\ 1-m,\ 1-b_2+a_2-m, 1-a_2+b_2-m],
\end{array}$$%\end{verbatim}
and is reducible of type [14]%.
\footnote{An operator $P$ decomposes (is reducible) for example of type [14] means $P$ can be written as $P_1\circ P_2$, where $P_1$ is of order 1 and $P_2$ is of order 4. [4] part means $P_2$. The decomposition is not necessarily an irreducible decomposition, which is not unique.}. We have $[1]=-x(x-1)^2( x(x-1)\partial - (m-6)x-2 )$, and 
\[x^2(x-1)^2{\rm [4]} =x^2(x-1)^2\partial^4+ m_{3}\partial^3+ m_{2}\partial^2+ m_{1}\partial^1+ m_{0},\]
where
\begin{eqnarray*}
m_{3} &=& -x(2x-1)(x-1)(2m-5), \\
m_{2} &=& 2c_2b_2x+2a_2b_2x^2+2c_2a_2x+3+25x^2-x^2b_2^2-24mx^2-x^2a_2^2+6m^2x^2\\
&& \quad -2c_2x-4a_2b_2x-c_2^2+2c_2-4m+m^2-24x-6m^2x+24mx, \\
m_{1} &=& (-3+2m)(-2m^2x+m^2-3m+6mx-c_2b_2+xa_2^2-c_2a_2+2-2a_2b_2x+2a_2b_2+c_2+xb_2^2-5x), \\
m_{0} &=& -(-1+m)^2(a_2-1+m-b_2)(a_2+1-m-b_2).
\end{eqnarray*}
The local exponents are
$$\begin{array}{ll}
x=0:&  [0,\  1,\  c_2-1+m,\  -c_2+1+m],\\
x=1:&  [0,\  1,\  -a_2+c_2-b_2+m,\  -c_2+b_2+a_2+m],\\
x=\infty:&  [1-m,\  1-m,\  -a_2+1-m+b_2,\  a_2+1-m-b_2].
\end{array}$$

\subsection{Tensor product without apparent singularities Case 2 }
In this section we assume
$$a_1 = -b_1+c_1+c_2-a_2-b_2,\quad c_1 = c_2=:c,$$
and study the middle convolution of $L$. 
This case happens to connect Gauss equations and the equation $\tilde{Z}(A)$, because the above assumption exactly fits the parameter change (see \S 5.2):
$$a_1=A_{-+-+},\quad b_1=A_{-++-},\quad a_2=A_{----},\quad b_2=A_{--++},\quad c=1-A_0.$$
The middle convolution $ML4(m)$ of $L$ with parameter $m$ is now computed without multiplying $\partial$ from the left, and we get
\[ ML4(m):=x^2(x-1)^2\partial^4 + m\ell_3 \partial^3+ m\ell_2 \partial^2+ m\ell_1 \partial+ m\ell_0,\] 
where
\begin{eqnarray*}
m\ell_3 &=& -x(2x-1)(x-1)(2m-5), \\
m\ell_2 &=& 3+2b_1xb_2+2b_1xa_2+25x^2-24x+2c-c^2-4m+2xb_1^2+2ca_2x^2 \\
&&\quad +2cb_2x^2-2a_2b_2x-2cx+2c^2x-2x^2b_1^2-x^2a_2^2-x^2b_2^2
-2c^2x^2 \\
&& \quad -4cb_1x+6m^2x^2-24mx^2-6m^2x+24mx+m^2-2x^2b_1b_2
 -2x^2b_1a_2+4x^2b_1c, \\
m\ell_1 &=& 
 (-3+2m)(m^2-2m^2x-3m+6mx+2+2c^2x+2b_1xa_2+2b_1xb_2+c-c^2+a_2b_2\\
&& \quad -5x-b_1^2+2xb_1^2-2ca_2x+2b_1c-2cb_2x-b_1a_2-b_1b_2
 -4cb_1x+xb_2^2+xa_2^2), \\
m\ell_0 &=& (b_1-1+m+b_2-c)(b_1+1-m+b_2-c)(b_1+1-c-m+a_2)(b_1-1-c+m+a_2).
\end{eqnarray*}
The local exponents of $ML4(m)$ are 
\begin{eqnarray*}
  &x=0:& [0, 1, c-1+m, 1-c+m],  \\
  &x=1:& [0, 1, b_2-c+m+a_2, -b_2-a_2+c+m], \\
  &x=\infty:& [-b_1+1+c-b_2-m, b_1+1+b_2-c-m, -b_1+1+c-m-a_2, 1+b_1-c-m+a_2].
\end{eqnarray*}
Recall that  $L(A)={\rm Ad}(x^{-A_0})(K(A))$. Since $c-1=-A_0$, our $L$ just agrees with $L(A)$.
The middle convolution $ML4(m)$ of $L=L(A)$ is given by
\begin{eqnarray*}
m\ell_3 &=& -2x(2x-1)(x-1)(-2+(m-1/2)), \\
m\ell_2 &=& 29/2x^2-29/2x-3(m-1/2)+9/4-6x(m-1/2)^2+18x(m-1/2)+6x^2(m-1/2)^2\\
 && \quad -18x^2(m-1/2)+(m-1/2)^2-A_0^2+xA_0^2+xA_3^2-xA_1^2+xA_2^2-x^2A_2^2-x^2A_3^2, \\
m\ell_1 &=& -1/2(-1+(m-1/2))(10x+8(m-1/2)-5+8x(m-1/2)^2-16x(m-1/2)-4(m-1/2)^2\\
 && \quad +2A_0^2+2A_3^2+2A_2^2 -2A_1^2-4xA_3^2-4xA_2^2),\\
m\ell_0 &=& 1/16(2A_2-1+2(m-1/2))(-2A_2-1+2(m-1/2))(-2A_3-1+2(m-1/2))\\
&&\quad \times(2A_3-1+2(m-1/2)).
\end{eqnarray*}
Thus we rediscovered Theorem 5.2: $ ML4(1/2)=\tilde{Z}(A).$

\subsubsection{$K_{S_1,S_2}$ when $A_2=0$}\label{KS}
The quadratic form $App$ giving the two apparent singularities of $K$ and $K_{S_1,S_2}$ reduces to $ -4A_2A_3x(x-1);$ note that
 $$S_1-S_2=-\frac{A_2A_3}{x(x-1)}.$$
The coefficients of $K_{S_1,S_2}=\partial^4+f_3dx^3+\cdots$ are given as
$$\begin{array}{ll}
f_3& := (2x-1)/(x(-1+x)),\\[2mm]
f_2& := (-A_2^2x^2-A_3^2x^2+A_0^2x-A_1^2x+A_2^2x+A_3^2x-A_0^2+x^2-x+1)/(x^2(-1+x)^2),\\[2mm]
f_1& := -(-2A_2^2x^3-2A_3^2x^3+5A_0^2x^2-5A_1^2x^2+3A_2^2x^2+3A_3^2x^2-9A_0^2x+A_1^2x\\[2mm]&-A_2^2x-A_3^2x+2x^3+4A_0^2-3x^2+9x-4)/(2x^3(-1+x)^3),\\[2mm]
f_0& := (2A_2^2A_3^2x^4-4A_2^2A_3^2x^3+2A_2^2A_3^2x^2-2A_2^2x^4-2A_3^2x^4+6A_0^2x^3\\[2mm]&-6A_1^2x^3+4A_2^2x^3+4A_3^2x^3-15A_0^2x^2+3A_1^2x^2-3A_2^2x^2-3A_3^2x^2+2x^4+13A_0^2x\\[2mm]&-A_1^2x+A_2^2x+A_3^2x-4x^3-4A_0^2+15x^2-13x+4)/(2x^4(-1+x)^4).
\end{array}$$
When $A_2=0$, $K_{S_1,S_2}$ decomposes of type [13]:
$$\begin{array}{ll}
  &[1]=\partial+(2x-1)/(x(-1+x)),\\[3mm]
     &[3]= \partial^3+(-A_3^2x^2+A_0^2x-A_1^2x+A_3^2x-A_0^2+x^2-x+1)\partial/(x^2(-1+x)^2)\\[2mm]
     &\qquad-(-2A_3^2x^3+3A_0^2x^2-3A_1^2x^2+3A_3^2x^2-5A_0^2x+A_1^2x-A_3^2x+2x^3\\[2mm]
     &\qquad+2A_0^2-3x^2+5x-2)/(2x^3(-1+x)^3).
\end{array}$$
When $A_1=A_2=0$ $((a_1,b_1,c_1)=(a_2,b_2,c_2))$, the 3rd order operator [3] above is nothing but the equation (\ref{K_S}): $K_S=\partial^3-4S\partial-2S'$.

\newpage
\normalsize
\part{Local solutions of ordinary differential equations related to the Dotsenko-Fateev equation}
In Part I, we found a Fuchsian system of rank 8 in 3 variables with 4 parameters, and an ordinary differential equation $Z(A)$ of order 4 with three singular points by restricting the system on a projective line. In Part II, we study the ordinary differential equation $Z(A)$ and several related ones, around their singular points.\par\noindent
In \S10, we study a linear difference equation $Rc_0(A)$ of order 2, which is the recurrence relation satisfied by the coefficients of a power series solution to $Z(A)$ at $x=0$. \par\noindent
\S10.1 introduces the invariant of such a difference equation. \par\noindent
\S10.3 introduces  special values $_4F_3(*;1)$ of terminating generalized hypergeometric series $_4F_3$ at 1 satisfying a linear difference equation $Rc^{(0)}$ of order 2.\par\noindent
In \S10.4,  by studying the invariant of this difference equation and that of $Rc_0(A)$, we find solutions of $Rc_0(A)$ expressed in terms of $_4F_3(*;1)$. This expression is very near to the product of two Gauss hypergeometric series. This observation leads to the discovery: A middle convolution sends $Z(A)$ to the product of two Gauss hypergeometric equation.\par\noindent
In \S10.5 and 6, invariants of the difference equations are  used to get local solutions of $Z(A)$ at $x=0,1$ and $\infty$. \par\noindent
In \S10.7, Riemann-Liouville transformation is recalled.\par\noindent
In \S10.8 and 9, we get local solutions of $Z(A)$ at $x=0,1$ and $\infty$ by using middle convolution. \par\noindent
In \S11, we study a linear difference equation $Rc_1(A)$ of order 2, which is the recurrence relation satisfied by the coefficients of difference equation a solution to $Z(A)$ at $x=0$ with exponents $A_0\pm1/2$. \par\noindent
\S11.2 introduces  special values $_4F_3(*;1)$ of non-terminating generalized hypergeometric series $_4F_3$ at 1 satisfying a linear difference equations of order 2.\par\noindent
\S11.3 studies the difference equation $Rc^{(1)}$ satisfied by this special values.\par\noindent
In \S11.4,  by studying the invariant of this difference equation and that of $Rc_1(A)$, we find solutions of $Rc_1(A)$ expressed in terms of non-terminating  $_4F_3(*;1)$. \par\noindent
\S12and 13 give a set of local solutions of $Q(A)$ and the Dotsenko-Fateev equation.\par\noindent

\section{Local solutions of $Z(A)$ at $x=0$ with exponent 0 and $2A_0$, and those at infinity}
The explicit form of the differential equation in question is given as
$$Z(A)=p_0\partial^4+p_1\partial^3+p_2\partial^2+p_3\partial+p_4,\quad \partial=d/dx,$$
where
$$\begin{array}{ll}
p_0&=\ x^3(x-1)^2,\\
p_1&=\ -2x^2(x-1)(2A_0x-2A_0-5x+3),\\
p_2&=\ x(p_{22}x^2+p_{21}x+p_{20}),\\
&p_{22}=\ 6A_0^2-A_2^2-A_3^2-24A_0+25,\\
&p_{21}=\ -11A_0^2-A_1^2+A_2^2+A_3^2+36A_0-59/2,\\
&p_{20}=\ (10A_0-9)(2A_0-3)/4,\\
p_3&=\ p_{32}x^2+p_{31}x+p_{30},\\
&p_{32}=\ -(2A_0-3)(2A_0^2-A_2^2-A_3^2-6A_0+5),\\
&p_{31}=\ (A_0-1)(6A_0^2+2A_1^2-2A_2^2-2A_3^2-16A_0+11),\\
&p_{30}=\ -((2A_0-3)(2A_0-1)^2)/4,\\
p_4&=\ p_{41}x+p_{40},\\
&p_{41}=\ (A_0-1+A_3)(A_0-1-A_3)(A_0-1+A_2)(A_0-1-A_2),\\
&p_{40}=\ -(2A_0-1)^2(A_0^2+A_1^2-A_2^2-A_3^2-2A_0+1)/4,
\end{array}$$
its Riemann scheme is given as
$$\left\{ \begin{array}{ccc}
x=0 & x=1 & x=\infty \\
0&\tfrac12-A_1&1-A_0+A_2\\
A_0-\tfrac12&0&1-A_0-A_2\\
A_0+\tfrac12&1&1-A_0+A_3\\
2A_0 & \tfrac12+A_1 &1-A_0-A_3
\end{array} 
\right\} .$$

In this section we give the two local solutions around the point $x = 0$,
whose exponents are 0 and $2A_0$. The other local solutions around $x = 0$, whose exponents are $A_0\pm1/2$, will be given in the next section.
The coefficients of the power series satisfy a 3-term recurrence relation (a homogeneous linear difference equation of order 2), say $Rc_0(A)$. 
\par\noindent
On the other hand, under some condition, special values of the terminating generalized hypergeometric series $_4F_3$ at 1 satisfy a linear difference equation of order 2. By making use of this fact, we solve the equation $Rc_0(A)$ in terms of the special values of $_4F_3$ at 1. The result suggests a relation between $Z(A)$ and the tensor product of two Gauss hypergeometric differential equations: they are connected by a middle convolution. This method will be applied to other equations at other singular points.
\par\noindent
The local exponents of the equation $\tilde{Z}(A)$ at $x=0$ are
$$0,\quad 1,\quad 1/2+A_0,\quad 1/2-A_0;$$
the holomorphic solution of $Z(A)$ at $x=0$ corresponds to the solution of $\tilde{Z}(A)$ with local exponent $1/2-A_0$, so the solution of $Z(A)$ at $x=0$ with exponent $2A_0$ corresponds to the solution of $\tilde{Z}(A)$ with exponent $1/2+A_0$, which can be obtained from that with exponent $1/2-A_0$ by changing the sign of $A_0$. 
\par\smallskip
But for the other exponents $A_0\pm1/2$ of $Z(A)$ at $x=0$ ($\{0,1\}$ for $\tilde{Z}(A)$), we need to make use of non-terminating series $_4F_3$. These are studied in the next section.
\subsection{Invariants of linear difference equations}\label{ebisinv} 
In this subsection, we introduce the invariants of linear difference equations following \cite{E}.
Let us consider a homogeneous linear difference equation of order 2 
$$P: C_n=p_1(n)C_{n-1}+p_2(n)C_{n-2},$$
where $p_1$ and $p_2$ are rational functions in $n$.
The quantity
$$H:=\frac{p_1(n)p_1(n+1)}{p_2(n+1)}$$
is called the {\bf invariant} of the difference equation $P$.
Consider another such equation
$$Q: D_n=q_1(n)D_{n-1}+q_2(n)D_{n-2}.$$
The two equations are said to be {\bf essentially the same} if there is a homogeneous linear
difference equation of order 1:
$$\lambda(n-1)=\mu(n)\lambda(n),\quad \mu(n):{\text{\rm a rational function in }}n$$
and a solution $\lambda(n)$ so that 
$$\{C_n\}=\lambda(n)\{D_n\},$$
where $\{C_n\}$ is the set of solutions of $P$, and $\{D_n\}$ that of solutions of $Q$.
\begin{prp}
The two equations $P$ and $Q$ are essentially the same if and only if the two invariants agree.
\end{prp}
In fact, substituting 
$$C_n=\lambda(n)D_n,\quad  C_{n-1}=\lambda(n-1)D_{n-1},\quad  C_{n-2}=\lambda(n-2)D_{n-2},$$
into $P$ and equalizing with $Q$, we have
$$\frac{q_1(n)}{p_1(n)}=\mu(n), \quad \frac{q_2(n)}{p_2(n)}=\mu(n)\mu(n-1),$$
and, by eliminating $\mu$,
$$\frac{p_1(n)p_1(n-1)}{p_2(n)}=\frac{q_1(n)q_1(n-1)}{q_2(n)}. $$
On the other hand if we assume the last identity, we can trace back the argument up to $C_n=\lambda(n)D_n.$ 
\begin{cor}\label{solspace}
For two essentially the same equations $P$ and $Q$, the spaces of solutions $\{C_n\}$ and $\{D_n\}$ are related as
\begin{align*}
  \{C_{n}\} 
  =
  \dfrac{1}{w^{n+1}} \cdot \dfrac{\prod _{j} \Gamma (n+v _{j}+1)}
{\prod _{i} \Gamma (n+u _{i}+1)}\cdot \{D_{n}\},
\end{align*}
where the rational function $\frac{q_1(n)}{p_1(n)}$ is factorized as 
$$
\frac{q_1(n)}{p_1(n)}
=
w \cdot\frac{\prod_i(n+u_i)}{\prod_j(n+v_j)},\qquad w:{\rm\ independent \ of\ }n.
$$
\label{ratio_coeff_difference_eq}
\end{cor}
\subsection{Recurrence relation $Rc_0(A)$ for the coefficients of a holomorphic solution of $Z(A)$ at $x=0$}
Let 
$$f^{(0,0)}(A; x):=\sum_{n=0}^\infty  C_nx^n,\quad  C_0=1$$
be the (normalized) power series solution to $Z(A)$ at $x=0$. Substituting this expression into $Z(A)$, we see that the coefficients $C_{n}$ satisfy the following recurrence relation $Rc_0(A)$ ($ C_{-1}=0, C_0=1$):
\begin{align*}
  Rc_{0}(A):
\begin{array}{ll}
  &C_{n}
  =
  \dfrac{
  \left\{2(n-A_0)-1\right\}^2 \left\{2n^2-4A_0n+A_{0} ^2+A_{1} ^2-A_{2} ^2-A_{3}^2+1-2(n-A_{0})\right\} 
  }{n(n-2A_0)(2n-2A_0-1)(2n-2A_0+1)}  C_{n-1}  \\[3mm]
  &\phantom{C_{n}}-
  \dfrac{4(n-A_0-A_2-1)(n-A_0+A_2-1)(n-A_0-A_3-1)(n-A_0+A_3-1)
  }{n(n-2A_0)(2n-2A_0-1)(2n-2A_0+1)}
  C_{n-2}.
\end{array}
\end{align*}
Thus we see that the { invariant of $Rc_0(A)$} is given as
\begin{align*}
H_0(n;A)
=
\dfrac{
 -\left\{
   4(n-A_{0})^{2}-1
 \right\}
 \left\{
   \left( 2n^{2}-4A_{0}n+A_{0}^{2}+A_{1}^{2}-A_{2}^{2}-A_{3}^{2}+1 \right)^{2}
   -
   4(n-A_{0})^{2}
 \right\}
 }{
   4n(n-2A_{0})
   (n-A_{0}-A_{2})(n-A_{0}+A_{2})
   (n-A_{0}-A_{3})(n-A_{0}+A_{3})
 }
 .
\end{align*}
\begin{rmk}[Symmetry]\label{symAc0} $Rc_0(A)$ is invariant under
  $$A_j\ \to\  -A_j\ (j=1,2,3)\quad{\rm and}\quad A_2\leftrightarrow A_3.$$
  \end{rmk}   \label{symmetryRc0}

\subsection{3-term relation for the special values (at 1) of balanced terminating hypergeometric series $_4F_3$}\label{3termrelfor4F3}
We consider 
the special value (at 1) of the generalized hypergeometric series: 
$$
{}_4F_3({\bm{\alpha}};1)
=
\sum_{k=0}^\infty
\frac{(\alpha_{0})_{k}(\alpha_{1})_{k}(\alpha_{2})_{k}(\alpha_{3})_{k}}
{(\beta_{1})_{k}(\beta_{2})_{k}(\beta_{3})_{k}}
\frac{1}{k!},
\qquad {\bm{\alpha}}
=
\begin{pmatrix}
    \alpha_{0},\alpha_{1},\alpha_{2},\alpha_{3}\\
    \beta_{1},\beta_{2},\beta_{3}
\end{pmatrix}.
$$
Here we assume $\alpha_0$ is a non-positive integer (the series  terminates), and also assume it is {\it balanced}:
\begin{equation*}%\label{balanced}
  \beta_{1}+\beta_{2}+\beta_{3}-\alpha_{0}-\alpha_{1}-\alpha_{2}-\alpha_{3}=1.\end{equation*}
Set \begin{align*}
    \bm{e}_{1} &:=
  \begin{pmatrix}
    1,0,0,0\\
    1,0,0
  \end{pmatrix},
  \quad
  \bm{e}_{2}
  :=
  \begin{pmatrix}
    0,1,0,0\\
    0,1,0
  \end{pmatrix},
  \quad
  \bm{e}_{12}
  :=
  \begin{pmatrix}
    1,1,0,0\\
    1,1,0
  \end{pmatrix}.
\end{align*}
We have
\begin{prp}[\cite{AAR} (\S3.7), see also \cite{Wils}]\label{AARW}
\begin{align}
  _{4}F_{3}({\bm{\alpha}} ;1)
  &=
  U^{(0)}_{1}({\bm{\alpha}})
  {}
  _{4}F_{3}({\bm{\alpha}}+{\bm{e}_{1}} ;1)
  +
  V^{(0)}_{1}({\bm{\alpha}})
  {}
  _{4}F_{3}({\bm{\alpha}}+{\bm{e}}_{12} ;1),
  \label{seed_gn1}
  \\
  _{4}F_{3}({\bm{\alpha}} ;1)
  &=
  U^{(0)}_{2}({\bm{\alpha}})
  {}
  _{4}F_{3}({\bm{\alpha}}+{\bm{e}_{2}} ;1)
  +
  V^{(0)}_{2}({\bm{\alpha}})
  {}
  _{4}F_{3}({\bm{\alpha}}+{\bm{e}}_{12} ;1),
  \label{seed_gn2}
\end{align}
where
\begin{align}
  &U_{1}^{(0)}({\bm{\alpha}})
  :=
  \dfrac{-\left( \beta_{1}-\alpha_{1} \right)\left( \beta_{1}+\beta_{2}-\alpha_{2}-\alpha_{3} \right)}
  {\beta_{1}\left( \beta_{3}-\alpha_{0}-1 \right)}
  ,\quad
  &
  &V_{1}^{(0)}({\bm{\alpha}})
  :=
  \dfrac{-\alpha_{1} \left( \beta_{2}-\alpha_{2} \right)\left( \beta_{2}-\alpha_{3} \right)}
  {\beta_{1}\beta_{2}\left( \beta_{3}-\alpha_{0}-1 \right)},
  &
  \label{coef_seed_gn1}
  \\
  &U_{2}^{(0)}({\bm{\alpha}})
  :=
  \dfrac{-\left( \beta_{2} -\alpha_{0} \right)\left( \beta_{1}+\beta_{2}-\alpha_{2}-\alpha_{3} \right)}
  {\beta_{2}\left( \beta_{3}-\alpha_{1}-1 \right)},
  \qquad
  &
  &V_{2}^{(0)}({\bm{\alpha}})
  :=
  \dfrac{-\alpha_{0} \left( \beta_{1}-\alpha_{2} \right)\left( \beta_{1}-\alpha_{3} \right)}
  {\beta_{1}\beta_{2}\left( \beta_{3}-\alpha_{1}-1 \right)}.
  &
  \label{coef_seed_gn2}
\end{align}
\end{prp}
%This fact can be derived from , since the series  terminates and is  balanced.
Perform a change 
$
\displaystyle
{\bm{\alpha}}
\mapsto
{\bm{\alpha}}+\bm{e}_{12}
$ 
in (\ref{seed_gn1}), and we have
\begin{align}
  _{4}F_{3}({\bm{\alpha}}+{\bm{e}_{12}} ;1)
  &=
  U^{(0)}_{1}({\bm{\alpha}}+{\bm{e}_{12}})
  {}
  _{4}F_{3}({\bm{\alpha}}+2{\bm{e}_{1}}+{\bm{e}_{2}} ;1)
  +
  V^{(0)}_{1}({\bm{\alpha}}+{\bm{e}_{12}})
  {}
  _{4}F_{3}({\bm{\alpha}}+2{\bm{e}_{12}} ;1),
  \label{twist_gn1}
\end{align}
and a change 
$
\displaystyle
{\bm{\alpha}}
\mapsto
{\bm{\alpha}}+{\bm{e_{1}}}
$ 
in (\ref{seed_gn2}) to get
\begin{align}
  _{4}F_{3}({\bm{\alpha}}+{\bm{e}_{1}} ;1)
  &=
  U^{(0)}_{2}({\bm{\alpha}}+{\bm{e}_{1}})
  {}
  _{4}F_{3}({\bm{\alpha}}+{\bm{e}_{12}} ;1)
  +
  V^{(0)}_{2}({\bm{\alpha}}+{\bm{e}_{1}})
  {}
  _{4}F_{3}({\bm{\alpha}}+2{\bm{e}_{1}}+{\bm{e}_{2}} ;1).
  \label{twist_gn2}
\end{align}
Eliminating 
$_{4}F_{3}({{\bm{\alpha}}}+2{\bm{e}_{1}}+{\bm{e}_{2}} ;1)$ 
and 
$_{4}F_{3}({\bm{\alpha}}+{\bm{e}_{1}} ;1)$ 
from (\ref{seed_gn1}), (\ref{twist_gn1}) and (\ref{twist_gn2}), we eventually get
\begin{lmm}\label{q0}$$
{}_{4}F_{3}({\bm{\alpha}} ;1)
=q^{(0)}_1({\bm{\alpha}})   {}_{4}F_{3}({\bm{\alpha}}+{\bm{e}_{12}} ;1)+
q^{(0)}_2({\bm{\alpha}})  {}_{4}F_{3}({\bm{\alpha}}+2{\bm{e}_{12}} ;1),
$$
where
\begin{align*} 
q^{(0)}_1({\bm{\alpha}})&=
  U^{(0)}_{1}({\bm{\alpha}})
  U^{(0)}_{2}({\bm{\alpha}}+{\bm{e}_{1}})
  +
  V^{(0)}_{1}({\bm{\alpha}})
  +
  \dfrac
  {U^{(0)}_{1}({\bm{\alpha}})V^{(0)}_{2}({\bm{\alpha}}+{\bm{e}_{1}})}
  {U^{(0)}_{1}({\bm{\alpha}}+{\bm{e}_{12}})}, 
     \\
q^{(0)}_2({\bm{\alpha}})&=
  \dfrac{ 
  -
    U^{(0)}_{1}({\bm{\alpha}})
  V^{(0)}_{1}({\bm{\alpha}}+{\bm{e}_{12}})
  V^{(0)}_{2}({\bm{\alpha}}+{\bm{e}_{1}})
}
  {U^{(0)}_{1}({\bm{\alpha}}+\bm{e}_{12})}. 
    %\label{seed_contig_rel}
\end{align*}
\end{lmm}
If we change $\bm{\alpha}$ into
  $$   (n;\hat{\bm{\alpha}}):=
  \begin{pmatrix}
    -n, \alpha_{1}-n,\alpha_{2},\alpha_{3}\\
    \beta_{1}-n,\beta_{2}-n,\beta_{3}
    \end{pmatrix}
  =
  -n\bm{e}_{12}
  +
  \hat{\bm{\alpha}},
  \quad
  \text{where \ } 
  \hat{\bm{\alpha}}
  =
  \begin{pmatrix}
    0, \alpha_{1},\alpha_{2},\alpha_{3}\\
    \beta_{1},\beta_{2},\beta_{3}
    \end{pmatrix}
%  \end{align*}
  $$
in the lemma above, we get
\begin{prp}\label{Rc0}\begin{gather*}
    D(n;\hat{\bm{\alpha}}):=\pFq{4}{3}{
    -n,\alpha_{1}-n,\alpha_{2},\alpha_{3}
  }{\beta_{1}-n,\beta_{2}-n,\beta_{3}}{1}
  \quad
  \text{\rm with \ $\beta_{1}+\beta_{2}+\beta_{3}-\alpha_{1}-\alpha_{2}-\alpha_{3}=1$}
\end{gather*}
satisfies the homogeneous linear difference equation
\begin{gather}
  Rc ^{(0)}({\hat{\bm{\alpha}}}):\
  D_{n}
  =
  q^{(0)}_{1}(n;{\hat{\bm{\alpha}}})D_{n-1}
  +
  q^{(0)}_{2}(n;{\hat{\bm{\alpha}}})D_{n-2}.
  \label{diff_eq_4f3_2200220}
\end{gather}  
\end{prp}\label{diffeqforterm4F3}

\normalsize\noindent
The { invariant of the difference equation} 
$Rc ^{(0)}({\hat{\bm{\alpha}}})$ 
is
\begin{gather*}
  H^{(0)}(n;\hat{{\bm{\alpha}}})
  :=
  \dfrac{q^{(0)}_{1}(n;{\hat{\bm{\alpha}}})\, q^{(0)}_{1}(n+1;{\hat{\bm{\alpha}}})}
  {q^{(0)}_{2}(n+1;{\hat{\bm{\alpha}}})}.
  %\label{H_diff_eq_4f3_2200220}
\end{gather*}
\begin{rmk}[Symmetry] $Rc ^{(0)}({\hat{\bm{\alpha}}})$ is invariant under
\begin{align}
  \alpha_{2}\leftrightarrow \alpha_{3}
  \quad
  \text{\rm and}
  \quad
  \beta_{1}\leftrightarrow \beta_{2}.
  \label{symmetry_trivial_terminate_4f3_2200220}
\end{align}
\end{rmk}
\subsection{Local solutions at zero I: solving $Rc_{0}(A)$}\label{SolveAc0}
\begin{prp}\label{hatal}The two invariants $H_{0}(n;A)$ and $H^{(0)}(n;\hat{{\bm{\alpha}}})$ of the difference equations  $Rc_{0}(A)$ and $Rc ^{(0)}({\hat{\bm{\alpha}}})$ agree, as functions in $n$, if and only if 
     \begin{align}
    \bm{\hat{\alpha}}
    =
  \begin{pmatrix}
    0, \alpha_{1},\alpha_{2},\alpha_{3}\\
    \beta_{1},\beta_{2},\beta_{3}
  \end{pmatrix}
    &=
  \begin{pmatrix}
      0, A_{0}, A_{-+-+}, A_{-++-}\\
      A_{++--}, A_{++++}, 1-A_{0}
  \end{pmatrix},
    \label{sol1_rc0}
    \\
    {\rm or}\quad&
    %\quad
 =\begin{pmatrix}
      0, A_{0}-A_{2}, A_{-+-+}, A_{-+--}\\
      A_{++--}, A_{++-+}, 1-A_{0}-A_{2}
  \end{pmatrix},
    \label{sol2_rc0}
  \end{align}
 up to the symmetries of the two difference equations.  Thus
  $Rc_{0}(A)$ and $Rc ^{(0)}({\hat{\bm{\alpha}}})$ 
  are essentially the same  in these cases.
  \label{sol_rc0}
\end{prp}
This can be obtained by solving the system $$H_{0}(n;A)=H^{(0)}(n;\hat{{\bm{\alpha}}})\quad n=1,2,...$$ with unknown $ \bm{\hat{\alpha}}$. Actual process is as follows: The numerator of $H_{0}(n;A)-H^{(0)}(n;\hat{{\bm{\alpha}}})$ is a
polynomial in $n$ of degree 14. The coefficient of the top term
decomposes as
\begin{align*}  {\rm constant}\times(\alpha_{1}-\beta_{1}-\beta_{2}+1-A_{1})
(\alpha_{1}-\beta_{1}-\beta_{2}+1+A_{1}).\end{align*}
If we set 
$\displaystyle \alpha_{1}=\beta_{1}+\beta_{2}-1\pm A_{1} $, 
then the second top term decomposes as
\begin{align*}   {\rm constant}\times(2A_{1}-1)(2A_{1}+1)
  (\alpha_{2}+\alpha_{3}-\beta_{1}-\beta_{2}+2A_{0}).\end{align*}
Since we assume $A_{1}$ is a free parameter, we have
\begin{align} \alpha_{2}+\alpha_{3}-\beta_{1}-\beta_{2}+2A_{0}=0. \label{eq_lmm1_proof_prop10.8}
\end{align}
Eliminating $\beta_2$ and $\beta_3$ from $H^{(0)}(n; \hat{\bm{\alpha}})$ by making use of this equality and the balance condition \begin{align}  \beta_{1}+\beta_{2}+\beta_{3}  -\alpha_{1}-\alpha_{2}-\alpha_{3}  =1
  \label{balanced_condition_terminate}
\end{align}
in Proposition \ref{Rc0}, we get an expression
\begin{align*}
  H^{(0)}(n; \hat{\bm{\alpha}})
  =\dfrac{
    -\gamma^{(0)} +O\left( \left( n-A_{0}\right)^{2} \right)
}{
  \gamma^{(0)} +O\left( \left( n-A_{0} \right)^{2} \right)  },\end{align*}
where
\begin{align}
  \gamma^{(0)}:=  A_{0}^{2}(\alpha_{1}-A_{0})^{2}  (\alpha_{2}-\beta_{1}+A_{0})^{2}
  (\alpha_{3}-\beta_{1}+A_{0})^{2}.
  \label{eq_lmm2_proof_prop10.8}  \end{align}
Now, assume $\gamma^{(0)} \neq 0.$
Then, $H^{(0)}(A_0; \hat{\bm{\alpha}})=-1,$
whereas $H_0(A_0; A$) is parameter dependent:
\begin{align*}
  H_{0}\left(A_{0}; A\right)
=
\dfrac{-\left( A_{0}^{2}-A_{1}^{2}+A_{2}^{2}+A_{3}^{2}-1 \right)^{2}}{4A_{0}^{2}A_{2}^{2}A_{3}^{2}}.
\end{align*}
So it does not happen that 
$H_0(n; A)=H^{(0)}(n; \hat{\bm{\alpha}})$ in the case $\gamma^{(0)} \neq 0.$
Therefore, $\gamma^{(0)}$ must be equal to 0. Thanks to the symmetry of $Rc^{(0)}(\hat{\bm{\alpha}})$ (cf. Remark \ref{symAc0}), we have only to study two cases:
\begin{align*}
  \text{Case 1}:\
  \alpha_{1}=A_{0}
  \quad{\rm and}\quad 
  \text{Case 2}:\
  \beta _{1}=\alpha _{3}+A_{0}.
\end{align*}
\par\medskip\noindent
\underline{Case 1}: 
Eliminating $\beta _{2}$ and $\beta _{3}$ from $H^{(0)}(n; \hat{\bm{\alpha}} )$ by (\ref{eq_lmm1_proof_prop10.8}) and (\ref{balanced_condition_terminate}), we have
\begin{align*}
  &H^{(0)}\left(n; \hat{\bm{\alpha}}\right)
=\dfrac{
  -\left\{4(n-A_{0})^{2}-1\right\}
\epsilon^{(0)}_{1}(n)\epsilon^{(0)}_{1}(n+1)}
{4n(n-2A_{0})
  (n+\alpha_{2}-\beta_{1})(n-\alpha_{2}+\beta_{1}-2A_{0})
  (n+\alpha_{3}-\beta_{1})(n-\alpha_{3}+\beta_{1}-2A_{0})},\\
&\text{where}\quad
\epsilon^{(0)}_{1}(n)=2\left(
n^{2}
-2A_{0}n+2A_{0}\beta_{1}+
\alpha_2\alpha_3+\alpha_2\beta_1+\alpha_3\beta_1-\beta_1^2
-\alpha_2-\alpha_3+1
-n\right).\end{align*}
Equating this and  $H_{0}(n; A)$, we conclude that $\hat{\bm{\alpha}}$ satisfies
\begin{itemize}
  \item
    $\displaystyle
    \epsilon^{(0)}_{1}(0)=A_{0}^{2}+A_{1}^{2}-A_{2}^{2}-A_{3}^{2}+2A_{0}+1, \phantom{+\dfrac{1}{2}}
  $
  \item
    $\displaystyle
  \left\{
    \pm \left( -\alpha_{2}+\beta_{1}-A_{0}\right)-A_{0},\ 
    \pm \left( -\alpha_{3}+\beta_{1}-A_{0}\right)-A_{0}  \right\}
  =
  \left\{    -A_{0}\pm A_{2},\     -A_{0}\pm A_{3}  \right\}.  $\end{itemize}
These lead to (\ref{sol1_rc0}), up to the symmetries of $Rc_0(A)$ and $Rc^{(0)}(\hat{\bm{\alpha}})$.
%'Ώ̐« (\ref{symmetryRc0}),(\ref{symmetry_trivial_terminate_4f3_2200220})
\par\medskip\noindent
\underline{Case 2}: In the same way we have
\begin{align*}
  &H^{(0)}\left(n; \hat{\bm{\alpha}}\right)
=
\dfrac{
  -\left\{4(n-A_{0})^{2}-1\right\}
\epsilon^{(0)}_{2}(n)\epsilon^{(0)}_{2}(n+1)}
{4n(n-2A_{0})
  (n-\alpha_{1})
  (n+\alpha_{1}-2A_{0})
  (n-\alpha_{2}+\alpha_{3}-A_{0})
  (n+\alpha_{2}-\alpha_{3}-A_{0})
},
\\
&\text{where}\quad
\epsilon^{(0)}_{2}(n)
=
2\bigr(
n^{2}
-2A_{0}n+2A_{0}^{2}-A_{0}\alpha_{1}+2A_{0}\alpha_{2}+2A_{0}\alpha_{3}
\\
&\phantom{\text{where}\quad
\epsilon^{(0)}_{1}(n)}
\qquad
-\alpha_{1}\alpha_{2}-\alpha_{1}\alpha_{3}+2\alpha_2\alpha_3
+\alpha_{1}-\alpha_{2}-\alpha_{3}+1-A_{0}
-n
\bigl),
\end{align*}
and the conditions 
\begin{itemize}
  \item
    $\displaystyle
    \epsilon^{(0)}_{2}(0)=A_{0}^{2}+A_{1}^{2}-A_{2}^{2}-A_{3}^{2}+2A_{0}+1, \phantom{+\dfrac{1}{2}}
  $
  \item
    $\displaystyle
  \left\{
    \pm \left( \alpha_{1}-A_{0}\right)-A_{0},\ 
    \pm \left( \alpha_{2}-\alpha_{3}\right)-A_{0}
  \right\}
  =
  \left\{
    -A_{0}\pm A_{2},\
    -A_{0}\pm A_{3}
  \right\},
  $
\end{itemize}
which lead to (\ref{sol2_rc0}), up to the symmetries of $Rc_0(A)$ and $Rc^{(0)}(\hat{\bm{\alpha}})$. \par\noindent
These conclude the proof the proposition.

\par\medskip
This proposition together with Corollary \ref{solspace} implies
\begin{prp}\label{solRc0} 
\begin{align*}
  (\text{Solutions of $Rc_{0}(A)$})
  &=
  \dfrac{
    \left( A_{----} \right)_{n}
    \left( A_{--++} \right)_{n}
  }{
    n!
    \left( \frac{3}{2}-A_{0} \right)_{n}
  }
  (\text
    {Solutions of 
    $Rc ^{(0)}({\hat{\bm{\alpha}}})$ 
    with $\bm{\hat{\alpha}}$ as (\ref{sol1_rc0})
    }
  )
  \\
  &=
  \dfrac{
    \left( 1-A_{0}-A_{2} \right)_{n}
    \left( A_{--+-} \right)_{n}
    \left( A_{--++} \right)_{n}
  }{
    n!
    \left( 1-A_{0} \right)_{n}
    \left( \frac{3}{2}-A_{0} \right)_{n}
  }
  (\text
    {Solutions of 
    $Rc ^{(0)}({\hat{\bm{\alpha}}})$ 
    with $\bm{\hat{\alpha}}$ as (\ref{sol2_rc0}) } ).
\end{align*}
\end{prp}
Picking the solutions with initial condition $C_{-1}=0,C_{0}=1$ up from the right hand-sides, we get
\begin{thm}\label{solZ4F3} The solution $f^{(0,0)}(A;x)$ of $Z(A)$ holomorphic at $x=0$ 
  (with normalization $f^{(0,0)}(A;0)=C_{0}=1$) can be expressed as
   \begin{align}
    &f^{(0,0)}(A;x) 
    \notag
    \\
    &=
    \sum_{n=0}^{\infty}x^{n}
\dfrac{
    \left(A_{----}% \frac{-A_{0}-A_{1}-A_{2}-A_{3}+1}{2} 
\right)_{n}
    \left(A_{--++}% \frac{-A_{0}-A_{1}+A_{2}+A_{3}+1}{2} 
\right)_{n}
  }{
    n!
    \left(  \frac{3}{2}-A_{0} \right)_{n}
  }
  %\\  &\times
\pFq{4}{3}{
      -n,
      A_{0}-n,
      A_{-+-+},%\frac{-A_{0}+A_{1}-A_{2}+A_{3}+1}{2},
      A_{-++-}%\frac{-A_{0}+A_{1}+A_{2}-A_{3}+1}{2}
    }{
      A_{++--}%\frac{A_{0}+A_{1}-A_{2}-A_{3}+1}{2}
-n,
      A_{++++}%\frac{A_{0}+A_{1}+A_{2}+A_{3}+1}{2}
-n,
      1-A_{0}
    }{1}
    \label{expression_z_0}
    \\
    \begin{split}
    &=
    \sum_{n=0}^{\infty}x^{n}
  \dfrac{
    \left( 1-A_{0}-A_{2} \right)_{n}
    \left( A_{--+-} \right)_{n}
    \left( A_{--++} \right)_{n}
  }{
    n!
    \left( 1-A_{0} \right)_{n}
    \left( \frac{3}{2}-A_{0} \right)_{n}
  }
\pFq{4}{3}{
      -n,
      A_{0}-A_{2}-n,
      A_{-+-+},%\frac{-A_{0}+A_{1}-A_{2}+A_{3}+1}{2},
      A_{-+--}%\frac{-A_{0}+A_{1}+A_{2}-A_{3}+1}{2}
    }{
      A_{++--}%\frac{A_{0}+A_{1}-A_{2}-A_{3}+1}{2}
-n,
      A_{++-+}%\frac{A_{0}+A_{1}+A_{2}+A_{3}+1}{2}
-n,
      1-A_{0}-A_{2}
    }{1}.
    \label{expression_z_0_another}
  \end{split}
    \end{align}
   %$f^{(0,0)}(A;x)$
   It is invariant under the symmetry of $Rc_0(A)$ in Remark $\ref{symAc0}$.
\label{sol_space_rc0}
\end{thm}
\begin{rmk}[Another way of deriving (\ref{expression_z_0_another}) from (\ref{expression_z_0})] By the way, the following identity is known:
  \begin{align}
    \pFq{4}{3}{-n,a,b,c} {d,e,f}{1}
    &=
    \dfrac{
      \left( e-a \right)_{n}\left( f-a \right)_{n}
    }{
      \left( e \right)_{n}\left( f \right)_{n}
    }
    \pFq{4}{3}{ -n,a,d-b,d-c} {d,a+1-n-e,a+1-n-f}{1}, 
    \label{trans_4F3_1}
  \end{align}
  where  $ a+b+c-n+1=d+e+f $ \ $($see Theorem $3.3.3$ in $\cite{AAR})$.
  Applying this transformation formula to the expression $(\ref{expression_z_0})$ by putting
  \begin{align*}
    \begin{pmatrix}
      a, b, c\\
      d, e, f
    \end{pmatrix}
    =
    \begin{pmatrix}
      A_{-+-+}, A_{0}-n, A_{-++-}\\
      A_{++--}-n, A_{++++}-n, 1-A_{0}
    \end{pmatrix}
  \end{align*}
we get the expression $(\ref{expression_z_0_another})$.
\end{rmk}

\begin{rmk}\label{discovermc}
  It is easily seen that 
  \begin{align}\label{2F12F14F3}
    \pFq{2}{1}{\alpha,\beta}{\gamma}{x}  \pFq{2}{1}{a,b}{c}{x}
    =  \sum_{n=0}^{\infty}     x^{n}
    \dfrac{(a)_{n}(b)_{n}}{(c)_{n}n!}
    \pFq{4}{3}  {-n,1-c-n,\alpha,\beta}{1-a-n,1-b-n,\gamma}{1},
      \end{align}
which implies  
\begin{align*}
&\pFq{2}{1}{
      A_{-+-+},%\frac{-A_{0}+A_{1}-A_{2}+A_{3}+1}{2},
      A_{-++-}%\frac{-A_{0}+A_{1}+A_{2}-A_{3}+1}{2}
    }{
      1-A_{0}
    }
   x% {\dfrac{t-1}{-2}}
    %\\   &
    \pFq{2}{1}{
      A_{----},%\frac{-A_{0}-A_{1}-A_{2}-A_{3}+1}{2},
      A_{--++}%\frac{-A_{0}-A_{1}+A_{2}+A_{3}+1}{2}
  }{
      1-A_{0}
}x%{\dfrac{t-1}{-2}}   
\\%\end{align*}\begin{align*}
&=\sum_{n=0}^{\infty}x^{n}
\dfrac{
    \left(A_{----}% \frac{-A_{0}-A_{1}-A_{2}-A_{3}+1}{2} 
\right)_{n}
    \left(A_{--++}% \frac{-A_{0}-A_{1}+A_{2}+A_{3}+1}{2} 
\right)_{n}
  }{
    n!
    \left( 1-A_{0}\right)_{n}
  }
  %\\  &\times
\pFq{4}{3}{
      -n,
      A_{0}-n,
      A_{-+-+},%\frac{-A_{0}+A_{1}-A_{2}+A_{3}+1}{2},
      A_{-++-}%\frac{-A_{0}+A_{1}+A_{2}-A_{3}+1}{2} 
    }{
      A_{++--}%\frac{A_{0}+A_{1}-A_{2}-A_{3}+1}{2}
-n,
      A_{++++}%\frac{A_{0}+A_{1}+A_{2}+A_{3}+1}{2}
-n,
      1-A_{0}
    }{1}.
\end{align*}%\label{Gauss2}
Compare this with  \eqref{expression_z_0}. It leads the authors to the discovery of the relationship between $\tilde{Z}(A)$ and $L(A)$:\ $\tilde{Z}(A)=mc_{1/2}(L(A))$ in Theorem $\ref{midconvZto2G}$ $($see also \S $\ref{RL}$$)$.
\end{rmk}
Let us define the following operation to state the above discovery impressively.
\begin{dfn}The operation $[X^n]$ is defined to pick up the coefficient of $X^n$ from a series $\sum_{n=0}^{\infty}C_{n}X^{n}$, that is,
\begin{gather*}
  [X^{n}]\sum_{n=0}^{\infty}C_{n}X^{n}:=C_{n}.
\end{gather*}
\end{dfn}
Then, the expression (\ref{expression_z_0}) can be rephrased as
 \begin{align} \begin{split} \label{expression_z_1}
  f^{(0,0)}(A;x)
  &=
  \sum _{n=0} ^{\infty}
  x^n%\left(\dfrac{1-t}{2}\right)^{n}
  \dfrac{(1-A_{0})_{n}}{\left( \frac{3}{2}-A_{0} \right)_{n}}
  [X^{n}]
  \pFq{2}{1}{A_{-+-+}, A_{-++-}}{1-A_{0}}{X}
  \pFq{2}{1}{A_{----}, A_{--++}}{1-A_{0}}{X}.
   \end{split}
\end{align} 
\begin{rmk}
  By combining 
  (\ref{expression_z_0_another}) and (\ref{2F12F14F3}),
  we can also obtain another expression of
  $f^{(0,0)}(A;x)$:
  \begin{align*}
  f^{(0,0)}(A;x)
  &=
  \sum _{n=0} ^{\infty}
  x^n%\left(\dfrac{1-t}{2}\right)^{n}
  \dfrac{(1-A_{0}-A_{2})_{n} (1-A_{0}+A_{2})_{n}}
  {\left( 1-A_{0} \right)_{n}\left( \frac{3}{2}-A_{0} \right)_{n}}
  [X^{n}]
  \pFq{2}{1}{A_{-+-+}, A_{-+--}}{1-A_{0}-A_{2}}{X}
  \pFq{2}{1}{A_{--+-}, A_{--++}}{1-A_{0}+A_{2}}{X}.
  \end{align*}
\end{rmk}

\subsection{Other local solutions expressed in terms of $f^{(0,0)}$}
Recall 
\begin{align*}
\tilde{Z}(A) = {\rm Ad}\left(x^{-A_0+\frac12}\right)Z(A)
=x^{-A_0+\tfrac12} \circ Z(A) \circ x^{A_0-\tfrac12}
\end{align*}
and that $\tilde{Z}(A)$ has symmetries
\begin{align}
&A_j\to-A_j\quad (j=0,1,2,3)
\label{symmetry_tildeZ_1}
%\quad &{\rm and} \quad 
;\quad A_2\longleftrightarrow A_3 ;\quad
\\
%\quad &{\rm and} \quad 
&(x,A_0,A_1)\longleftrightarrow (1-x,A_1,A_0).
\label{symmetry_tildeZ_2}
\end{align}
In particular, the symmetry $A_{0}\rightarrow -A_{0}$ for $\tilde{Z}(A)$
implies that if $f(A;x)$ is a solution of $Z(A)$, then 
$x^{-A_{0}+1/2}f(A;x)$ and 
$x^{A_{0}+1/2}f(-A_{0},A_{1},A_{2},A_{3};x)$ satisfy $\tilde{Z}(A)$,
that is,
$x^{2A_{0}}f(-A_{0},A_{1},A_{2},A_{3};x)$ also solves  $Z(A)$.
Therefore, we get the following:
\begin{prp}
  Let 
  $f^{(0,2A_{0})}(A;x)$ 
  be the normalized local solution
  of $Z(A)$ at $x=0$ with exponent $2A_{0}$.
  Then, 
  $f^{(0,2A_{0})}(A;x)$ 
  is expressible as
  \begin{align*}
  f^{(0,2A_{0})}(A;x)
  =
  x^{2A_{0}}f^{(0,0)}(-A_{0},A_{1},A_{2},A_{3};x).
  \end{align*}
\end{prp}
Similarly,
the symmetry (\ref{symmetry_tildeZ_2}) leads to the following:
\begin{prp}
  Let   $f^{(1,\pm A_{1})}(A;x)$   be the normalized local solutions
  of $Z(A)$ at $x=1$ with exponents $\frac{1}{2}\pm A_{1}$,
  respectively.   Then, 
  $f^{(1,\pm A_{1})}(A;x)$   are expressed as
  \begin{align*}
    f^{(1,\pm A_{1})}(A;x)
  &=
  x^{A_{0}-\frac12}(1-x)^{\frac{1}{2}\pm A_{1}}
  f^{(0,0)}(\mp A_{1},A_{0},A_{2},A_{3};1-x),
  \end{align*}
respectively.
\end{prp}
Thus, we have obtained the normalized local solutions at $x=0$ with exponents $0, 2A_0$, and at $x=1$ with exponents $1/2\pm A_1$.

\subsection{Local solutions at infinity I: using invariants of the difference equations}\label{infty}
In this section, we start to find the normalized local solutions of $Z(A)$ at $x=\infty$. Recall that the local exponents at $x=\infty$ are
$$1-A_0 \pm A_2,\ \ 1-A_0 \pm A_3.$$
We find the normalized local solution with exponent $1-A_0+A_2$.
By substituting the expression
$$f^{(\infty, +A_{2})}(A;x)
:=
\left( \dfrac{1}{x} \right)^{1-A_0+A_2}\sum_{n=0}^\infty C_n \left( \dfrac{1}{x} \right)^n, \quad \text{where\quad }C_{0}=1.$$
 into $Z(A)$, we see that coefficients $C_n$ satisfy the following 3-term recurrence relation
\begin{align*}
Rc_{\infty}(A):
\begin{array}{ll}
&
C_n
=
\dfrac{
\left\{2(n+A_2)-1\right\}^2\left\{ 2n^{2}+4A_{2}n -A_0^2+A_1^2+A_{2}^{2}-A_3^2+1-2(n+A_{2})\right\}}{
4n(n+2A_2)(n+A_2+A_3)(n+A_2-A_3)
}
C_{n-1}
\\[3mm]
& 
\phantom{C_n}
-
\dfrac{(2n+2A_2-1)(2n+2A_2-3)(n+A_0+A_2-1)(n-A_0+A_2-1)}{
4n(n+2A_2)(n+A_2+A_3)(n+A_2-A_3)
}
C_{n-2}.
\end{array}
\end{align*}
Thus the invariant $H_\infty(n;A)$ of $Rc_\infty(A)$ is given by
\begin{align*}
H_\infty(n;A)
=
H_{0}(n; -A_{2}, A_{1}, A_{0}, A_{3}).
\end{align*}
\begin{rmk}[Symmetry] $Rc_\infty(A)$ is invariant under $ A_{j}\rightarrow -A_{j}
\ (j=0, 1, 3)$, not under $A_2\leftrightarrow A_3$.
 \label{symmetryRcinf} \end{rmk}
Hence, Proposition \ref{sol_rc0} yields the following proposition;
\begin{prp}The two invariants  $H^{(0)}(n; \hat{\bm{\alpha}})$ and  $H_{\infty}(n; A)$ 
   of the difference equations $Rc ^{(0)}({\hat{\bm{\alpha}}})$ and $Rc_{\infty}(A)$, respectively,  agree if and only if   
    \begin{align}    \bm{\hat{\alpha}}    =
  \begin{pmatrix}
    0, \alpha_{1},\alpha_{2},\alpha_{3}\\
    \beta_{1},\beta_{2},\beta_{3}
  \end{pmatrix}
    &=
  \begin{pmatrix}
      0, -A_{2}, A_{-+++}, A_{+++-}\\
      A_{-+--}, A_{++-+}, 1+A_{2}
  \end{pmatrix},
    \label{sol1_rcinfty}
    \\
   {\rm or} \quad   &   
  =\begin{pmatrix}
      0, -A_{0}-A_{2}, A_{-+++}, A_{-++-}\\
      A_{-+--}, A_{-+-+}, 1-A_{0}+A_{2}
  \end{pmatrix},
    \label{sol2_rcinfty}
    \\
 {\rm or} \quad   &   
  =\begin{pmatrix}
      0, -A_{2}-A_{3}, A_{+++-}, A_{-++-}\\
      A_{-+--}, A_{++--}, 1+A_{2}-A_{3}
  \end{pmatrix},
    \label{sol3_rcinfty}
  \end{align}
  \label{sol_rcinfty}
up to the symmetries of the difference equations.
\end{prp}
In a similar way to \S \ref{SolveAc0}, by computing the ratio of the coefficients of the two difference equations 
$Rc_{\infty}(A)$ and $Rc ^{(0)}({\hat{\bm{\alpha}}})$,
we get
\begin{prp} 
\begin{align*}
  &(\text{Solutions of $Rc_{\infty}(A)$})
  \\
  &=
  \dfrac{
    \left( \frac{1}{2}+A_{2} \right)_{n}
    \left( A_{--+-} \right)_{n}
    \left( A_{+-++} \right)_{n}
  }{
    n!
    \left( 1+A_{2}-A_{3} \right)_{n}
    \left( 1+A_{2}+A_{3} \right)_{n}
  }
  (\text
    {Solutions of 
    $Rc ^{(0)}({\hat{\bm{\alpha}}})$ 
    with $\bm{\hat{\alpha}}$ as (\ref{sol1_rcinfty})
    }
  )
  \\
  &=
  \dfrac{
    \left( \frac{1}{2}+A_{2} \right)_{n}
    \left( 1-A_{0}+A_{2} \right)_{n}
    \left( A_{+-++} \right)_{n}
    \left( A_{+-+-} \right)_{n}
  }{
    n!
    \left( 1+A_{2} \right)_{n}
    \left( 1+A_{2}-A_{3} \right)_{n}
    \left( 1+A_{2}+A_{3} \right)_{n}
  }
  (\text
    {Solutions of 
    $Rc ^{(0)}({\hat{\bm{\alpha}}})$ 
    with $\bm{\hat{\alpha}}$ as (\ref{sol2_rcinfty})
    }
  )
  \\
  &=
  \dfrac{
    \left( \frac{1}{2}+A_{2} \right)_{n}
    \left( A_{--++} \right)_{n}
    \left( A_{+-++} \right)_{n}
  }{
    n!
    \left( 1+A_{2} \right)_{n}
    \left( 1+A_{2}+A_{3} \right)_{n}
  }
  (\text
    {Solutions of 
    $Rc ^{(0)}({\hat{\bm{\alpha}}})$ 
    with $\bm{\hat{\alpha}}$ as (\ref{sol3_rcinfty})
    }
  ).
\end{align*}
\label{sol_space_rcinfty}
\end{prp}
Picking 
the solutions with initial condition $C_{-1}=0, C_{0}=1$
up from the right hand-sides, 
we get
\begin{thm}\label{series_sol_z_infty_+A2} 
  $f^{(\infty,+A_{2})}(A;x)$
  %,
  %which is the normalized solution
  %of $Z(A)$ at $x=\infty$ with exponent $-A_{0}+1+A_{2}$,
  is given as
\begin{align*}
  f^{(\infty,+A_{2})}(A;x)
  =  
  \left( \dfrac1x \right)^{1-A_{0}+A_{2}}
    \sum_{n=0}^{\infty}
      &\left( \dfrac1x \right)^{n}
  \dfrac{
    \left( \frac{1}{2}+A_{2} \right)_{n}
    \left( A_{--+-} \right)_{n}
    \left( A_{+-++} \right)_{n}
  }{
    n!
    \left( 1+A_{2}-A_{3} \right)_{n}
    \left( 1+A_{2}+A_{3} \right)_{n}
  }
  \\
  &\times\pFq{4}{3}{
      -n,
      -A_{2}-n,
      A_{-+++},
      A_{+++-}
    }{
      A_{-+--}-n,
      A_{++-+}-n,
      1+A_{2}
    }{1}
    \\
 =     
  \left( \dfrac1x \right)^{1-A_{0}+A_{2}}
    \sum_{n=0}^{\infty}
      &\left( \dfrac1x \right)^{n}
  \dfrac{
    \left( \frac{1}{2}+A_{2} \right)_{n}
    \left( 1-A_{0}+A_{2} \right)_{n}
    \left( A_{+-++} \right)_{n}
    \left( A_{+-+-} \right)_{n}
  }{
    n!
    \left( 1+A_{2} \right)_{n}
    \left( 1+A_{2}-A_{3} \right)_{n}
    \left( 1+A_{2}+A_{3} \right)_{n}
  }
  \\
  &\times\pFq{4}{3}{
      -n,
      -A_{0}-A_{2}-n,
      A_{-+++},
      A_{-++-}
    }{
      A_{-+--}-n,
      A_{-+-+}-n,
      1-A_{0}+A_{2}
    }{1}
    \\
 =     
  \left( \dfrac1x \right)^{1-A_{0}+A_{2}}
    \sum_{n=0}^{\infty}
      &\left( \dfrac1x \right)^{n}
  \dfrac{
    \left( \frac{1}{2}+A_{2} \right)_{n}
    \left( A_{--++} \right)_{n}
    \left( A_{+-++} \right)_{n}
  }{
    n!
    \left( 1+A_{2} \right)_{n}
    \left( 1+A_{2}+A_{3} \right)_{n}
  }
  \\
  &\times\pFq{4}{3}{
      -n,
      -A_{2}-A_{3}-n,
      A_{+++-},
      A_{-++-}
    }{
      A_{-+--}-n,
      A_{++--}-n,
      1+A_{2}-A_{3}
    }{1}.
    \end{align*}
 The function $x^{1-A_{0}+A_{2}}f^{(\infty,+A_{2})}(A;x)$ is invariant under the change of parameters $A_j\to -A_j\ (j=0,1,3)$.%s in Remark \ref{symmetryRcinf}.
    \end{thm}
Note that one of these three expressions yield the other two by the help of the transformation formula (\ref{trans_4F3_1}) as stated before. Note also that these expressions are rephrased in terms of products of $_2F_1$
using (\ref{2F12F14F3}):
\begin{thm}
  $f^{(\infty,+A_{2})}(A;x)$
  is expressible as
    \begin{align*}
      f^{(\infty, +A_{2})}(A;x)
      &=\left( \dfrac1x \right)^{1-A_{0}+A_{2}}
    \sum_{n=0}^{\infty}
      \left( \dfrac1x \right)^{n}
  \dfrac{
    \left( \frac{1}{2}+A_{2} \right)_{n}
    \left( 1+A_{2} \right)_{n}
  }{
    \left( 1+A_{2}-A_{3} \right)_{n}
    \left( 1+A_{2}+A_{3} \right)_{n}
  }
  \\
  &
  \phantom{
      \left( \dfrac1x \right)^{-A_{0}+1+A_{2}}
    \sum_{n=0}^{\infty}
  }
  \times
  [X^{n}]
  \pFq{2}{1}{
      A_{-+++},
      A_{+++-}
    }{
      1+A_{2}
    }{X}
  \pFq{2}{1}{
      A_{+-++},
      A_{--+-}
    }{
      1+A_{2}
    }{X}
  \\
  &=
      \left( \dfrac1x \right)^{1-A_{0}+A_{2}}
    \sum_{n=0}^{\infty}
      \left( \dfrac1x \right)^{n}
  \dfrac{
    \left( \frac{1}{2}+A_{2} \right)_{n}
    \left( 1-A_{0}+A_{2}\right)_{n}
    \left( 1+A_{0}+A_{2}\right)_{n}
  }{
    (1+A_{2})_{n}(1+A_{2}-A_{3})_{n}(1+A_{2}+A_{3})_{n}
  }
  \\
  &
  \phantom{
      \left( \dfrac1x \right)^{-A_{0}+1+A_{2}}
    \sum_{n=0}^{\infty}
  }
  \times
  [X^{n}]
  \pFq{2}{1}{
      A_{-+++},
      A_{-++-}
    }{
      1-A_{0}+A_{2}
    }{X}
  \pFq{2}{1}{
      A_{+-++},
      A_{+-+-}
    }{
      1+A_{0}+A_{2}
    }{X}
  \\
  &=
      \left( \dfrac1x \right)^{1-A_{0}+A_{2}}
    \sum_{n=0}^{\infty}
      \left( \dfrac1x \right)^{n}
  \dfrac{
    \left( \frac{1}{2}+A_{2} \right)_{n}
  }{
    (1+A_{2})_{n}
  }
  \\
  &
  \phantom{
      \left( \dfrac1x \right)^{-A_{0}+1+A_{2}}
    \sum_{n=0}^{\infty}
  }
  \times
  [X^{n}]
  \pFq{2}{1}{
      A_{+++-},
      A_{-++-}
    }{
      1+A_{2}-A_{3}
    }{X}
  \pFq{2}{1}{
      A_{+-++},
      A_{--++}
    }{
      1+A_{2}+A_{3}
    }{X}.
    \end{align*}
    \label{series_sol_z_infty_+A2_another}
\end{thm}
The other three local solutions at $x=\infty$ are easily obtained from Theorem \ref{series_sol_z_infty_+A2} by recalling symmetries of $Z(A)$
$$A_j\to-A_j\quad (j=1,2,3)\quad {\rm and}\quad A_2\leftrightarrow A_3.$$
\begin{thm}
  Let
  $f^{(\infty, -A_{2})}(A;x)$
  ,
  $f^{(\infty, \pm A_{3})}(A;x)$
  be normalized local solutions of $Z(A)$ at $x=\infty$
  with exponents $1-A_{0}-A_{2}$, $1-A_{0}\pm A_{3}$, respectively.
  Then, these can be expressed as
  \begin{align*}
  f^{(\infty, -A_{2})}(A;x)  &=
  f^{(\infty, +A_{2})}(A_{0},A_{1},-A_{2},A_{3};x),
  \\
  f^{(\infty, \pm A_{3})}(A;x)
  &=
  f^{(\infty, +A_{2})}(A_{0},A_{1},\pm A_{3},A_{2};x),
  \end{align*}
  where 
  $f^{(\infty, +A_{2})}(A;x)$ are given in Theorem 
  $\ref{series_sol_z_infty_+A2}$ and 
 $\ref{series_sol_z_infty_+A2_another}$.
\end{thm}
\subsection{Correspondence of solutions via the Riemann-Liouville transformation}\label{RL}
For a linear differential operator $P$, the middle convolution $mc_\mu P$ of $P$ with parameter $\mu$ is defined in \S\ref{midconv}  as the linear differential operator $\partial^{-\mu}\circ P\circ\partial^\mu$.
\par\medskip
On the other hand, for a function $u(x)$, the notion of Riemann-Liouville transformation of $u$ with parameter $\mu$ is defined as the function in $x$:
$$\left(I_\gamma^{\mu}u\right)(x):= \frac1{\Gamma(\mu)}\int_\gamma u(s)(x-s)^{\mu-1}ds,$$
where $\gamma$ is a cycle.%loaded cycle in the sense of \cite{Yos}.
\par\medskip
It is known ([Hara2]) that if $u$ is a solution of $P$, then the function given by the integral above is a solution of $mc_\mu P$.
\par\medskip
If $u(x)$ is given locally around $x=0$ as
$$u=x^\alpha\sum_{j=0}^\infty c_nx^n\quad (\alpha\notin \Z),$$
and if $\mu\notin\Z$, then we can choose $\gamma$ as a path from $0$ to $x$:
$$\left(I_{[0,x]}^{\mu}u\right)(x)=\frac1{\Gamma(\mu)}\int_0^xu(s)(x-s)^{\mu-1}ds.$$
Since we assume $\alpha,\mu\notin\Z$, we can apply the beta function formula to get
$$\left(I_{[0,x]}^{\mu}u\right)(x)=\frac{\Gamma(1+\alpha)}{\Gamma(1+\alpha+\mu)}x^{\alpha+\mu}\sum_{n=0}^\infty\frac{(1+\alpha)_n}{(1+\alpha+\mu)_n}c_nx^n.$$
\subsection{Partial correspondence of local solutions at $x=0,1$}
Recall Theorem \ref{midconvZto2G}:
$$\tilde Z(A):=x^{-A_0+1/2}Z(A)x^{A_0-1/2}=mc_{\frac12}L(A).$$
This suggests us to apply the above transformation formula for $$\alpha=-A_0,\quad \mu=1/2$$
and for $u$ the product of the two Gauss hypergeometric series multiplied by $x^{-A_0}$:
$$u=x^{-A_0}\sum c_nx^n,\quad c_n=[X^n]\pFq{2}{1}{A_{-+-+}, A_{-++-}}{1-A_{0}}{x}
  \pFq{2}{1}{A_{----}, A_{--++}}{1-A_{0}}{x}.$$
We get a solution of $\tilde{Z}(A)$:
$$\left(I_{[0,x]}^{1/2}u\right)(x)=\displaystyle \frac{\Gamma(1-A_0)}{\Gamma(3/2-A_0)}x^{-A_0+1/2}\sum\frac{(1-A_0)_n}{(3/2-A_0)_n}c_nx^n,$$
which is, by the definition of $f^{(0,0)}(A,x)$, equal to
$$\displaystyle\frac{\Gamma(1-A_0)}{\Gamma(3/2-A_0)}x^{-A_0+1/2}f^{(0,0)}(A,x).$$
This rediscovers the expression (\ref{expression_z_1}), and we have
\begin{prp}Via the Riemann-Liouville transformation $I_{[0,x]}^{1/2}$ above, the local solution at $x=0$ to $L(A)$, the product of two Gauss equation, of exponent $\pm A_0$ is sent to the local solution at $x=0$ to $\tilde{Z}(A)$ of exponent $\tfrac12\pm A_0$.\end{prp}
\par\medskip
Since the equation is stable under the change
$$(x,A_0, A_1, A_2,A_3) \mapsto (1-x, A_1, A_0, A_2, A_3)$$
as in Remark \ref{symmZA}, the happening at $x=1$ reduces to that at $x=0$.
\subsection{Local solutions at infinity II: using middle convolution}
Since the local exponents  $$\alpha=-1\pm A_2,\ 1\pm A_3$$ of $L(A)$ at $\infty$ are non integral, the corresponding solutions with exponents $$\alpha=-1/2\pm A_2,\ 3/2\pm A_3$$ of $\tilde{Z}(A)$ are obtained via the Riemann-Liouville transformation with parameter $1/2$ from those of $L(A)$. Just apply to the local solution at infinity
\begin{align*}
u(x) = x^{\alpha} \sum_{n=-\infty}^0 c_n x^n= \sum_{n=-\infty}^0 c_n x^{n+\alpha}\qquad \alpha\in\{-1\pm A_2,\ 1\pm A_3\}
\end{align*}
to get
$$ \frac{\Gamma(\alpha+1)}{\Gamma(\alpha+\tfrac32)} 
x^{\alpha+\tfrac12}\sum_{n=-\infty}^0 c_n \frac{(-\alpha-\tfrac12)_n}
{(-\alpha)_n} x^{n}.$$

\section{Local solutions of $Z(A)$ at $x=0$ with exponent $A_0\pm1/2$}
In this section,
local solutions of $Z(A)$ at $x=0$ with exponent $A_0\pm \frac12$ are constructed.
\subsection{Recurrence relation $Rc_1(A)$}
The coefficients of a local solution
$$x^{A_0-\frac{1}{2}}\sum_{n=0}^\infty C_nx^n$$
satisfy
\begin{align*}
Rc_{1}(A):
\begin{array}{ll}
&
C_n
=
\dfrac
{(n-1)^2\left(2n^{2}-4n-A_{0}^2+A_{1}^2-A_{2}^2-A_{3}^2+\frac{5}{2}\right)}
{n(n-1)\left(n+A_{0}-\frac{1}{2}\right)\left(n-A_{0}-\frac{1}{2}\right)}
C_{n-1}
\\
&
\phantom{C_n}
-
\dfrac
{\left(n+A_{2}-\frac{3}{2}\right)\left(n-A_{2}-\frac{3}{2}\right)\left(n+A_{3}-\frac{3}{2}\right)\left(n-A_{3}-\frac{3}{2}\right)}
{n(n-1)\left(n+A_{0}-\frac{1}{2}\right)\left(n-A_{0}-\frac{1}{2}\right)}
C_{n-2}.
\end{array}
\qquad (n=2, 3, \ldots)
\end{align*}
For arbitrary given $C_0$ and $C_1$, 
remaining coefficients $C_n(n\ge2)$ are uniquely determined.
The invariant of $Rc_{1}(A)$ is given by
\begin{align*}
  H_{1}(n;A)
  =
  \dfrac{
    -n(n-1)
    \left( 2n^{2}-4n-A_{0}^{2}+A_{1}^{2}-A_{2}^{2}-A_{3}^{2}+\frac{5}{2} \right)
    \left( 2n^{2}   -A_{0}^{2}+A_{1}^{2}-A_{2}^{2}-A_{3}^{2}+\frac{1}{2} \right)
  }{
    \left( n+A_{0}-\frac{1}{2} \right)
    \left( n-A_{0}-\frac{1}{2} \right)
    \left( n+A_{2}-\frac{1}{2} \right)
    \left( n-A_{2}-\frac{1}{2} \right)
    \left( n+A_{3}-\frac{1}{2} \right)
    \left( n-A_{3}-\frac{1}{2} \right)
  }.
\end{align*}
\begin{rmk}[Symmetry] $Rc_{1}(A)$ is invariant under
\begin{align}
  A_j\to-A_j\quad (j=0,1,2,3)\quad {\rm and}\quad A_2\leftrightarrow A_3.
  \label{symmetry_Rc1}
\end{align}
\end{rmk}

\subsection{Special values of non-terminating $_4F_3$ at $x=1$}
The difference equation $Rc_0(A)$ for holomorphic solutions at $x=0$ was solved by special values of the {\it terminating} series $_4F_3(*;1)$, which satisfy the difference equation $Rc^{(0)}(\hat{\bm{\alpha}})$. The key was to find the parameters $\hat{\bm{\alpha}}$ so that the invariants of the two difference equations
$$H_0(n;A)\quad{\rm and}\quad H^{(0)}(n;\hat{\bm{\alpha}}) $$
agree.

For the difference equation $Rc_1(A)$ for general parameters $A=(A_0,A_{1},A_{2},A_3)$, we can not find $\hat{\bm{\alpha}}$ so that the invariants $H_1(n;A)$ and $H^{(0)}(n;\hat{\bm{\alpha}})$ agree. In other words, the {\it terminating} series $_4F_3(*;1)$ can not solve our equation $Rc_1(A)$.

We introduce some  special values of non-terminating $_4F_3$ at 1: We first introduce
$$
  {}_{4}\tilde{f}_{3}
\biggl(\genfrac..{0pt}{}
{a_{0}, a_{1}, a_{2}, a_{3}}
{b_{0}, b_{1}, b_{2}, b_{3}};x\biggr)
:=
  \sum _{n=0} ^{\infty}
  \dfrac
  {\Gamma(a_{0}+n)\Gamma(a_{1}+n)\Gamma(a_{2}+n)\Gamma(a_{3}+n)}
  {\Gamma(b_{0}+n)\Gamma(b_{1}+n)\Gamma(b_{2}+n)\Gamma(b_{3}+n)}
  x^{n},
$$  
  and\footnote{This series is known to be convergent at $x=1$ if
$$\Re( b_0+b _1 + b _2 +b _3 - a _0 - a _1 - a _2 - a _3 )>1.$$}
\begin{align*}
  \pfq{4}{3}
{a_{0}, a_{1}, a_{2}, a_{3}}
{b_{1}, b_{2}, b_{3}}{x}
  &:= {}_{4}\tilde{f}_{3}
\biggl(\genfrac..{0pt}{}
{a_{0}, a_{1}, a_{2}, a_{3}}
{1, b_{1}, b_{2}, b_{3}};x\biggr)
  = \dfrac
  {\Gamma(a_{0})\Gamma(a_{1})\Gamma(a_{2})\Gamma(a_{3})}
  {\Gamma(b_{1})\Gamma(b_{2})\Gamma(b_{3})}
  {}\pFq{4}{3}
{a_{0}, a_{1}, a_{2}, a_{3}}
{b_{1}, b_{2}, b_{3}}{x}.
\end{align*}
For
$${\bm{\alpha}}
=
\begin{pmatrix}
    \alpha_{0},\alpha_{1},\alpha_{2},\alpha_{3}\\
    \beta_{1},\beta_{2},\beta_{3}
\end{pmatrix},$$
we define
  \begin{align*}
      y_{0}(\bm{\alpha})
&:=
{}_{4}f_{3}(\bm{\alpha}; 1)
:=
\pfq{4}{3}{\alpha_{0}, \alpha_{1}, \alpha_{2}, \alpha_{3}}
{\beta_{1}, \beta_{2}, \beta_{3}}{1}, \\
  y_{i}(\bm{\alpha})
  &:=
  {}_{4}\tilde{f}_{3}
\biggl(\genfrac..{0pt}{}
  {\alpha_{0}+1-\beta_{i}, \alpha_{1}+1-\beta_{i}, \alpha_{2}+1-\beta_{i}, \alpha_{3}+1-\beta_{i}}
{2-\beta_{i}, \beta_{1}+1-\beta_{i}, \beta_{2}+1-\beta_{i}, \beta_{3}+1-\beta_{i}};1\biggr),
 &(i=1, 2, 3)
\\
  y_{i+4}(\bm{\alpha})
  &:=
  -
  {}_{4}\tilde{f}_{3}
\biggl(\genfrac..{0pt}{}
  {\alpha_{i}, \alpha_{i}+1-\beta_{1}, \alpha_{i}+1-\beta_{2}, \alpha_{i}+1-\beta_{3}}
{\alpha_{i}+1-\alpha_{0}, \alpha_{i}+1-\alpha_{1}, \alpha_{i}+1-\alpha_{2}, \alpha_{i}+1-\alpha_{3}};1\biggr).
&(i=0, 1, 2, 3)
  \end{align*}
  From now on we always assume that the parameters are balanced:
  \begin{align}
    \beta _{1}+\beta _{2}+\beta _{3}-
    \alpha _{0}-\alpha _{1}-\alpha _{2}-\alpha _{3}=1;
    \label{balanced_condition_nonterminate}
  \end{align}
  so all the infinite series above are convergent.
%Of course, any $y_{i}$ $(i=1,2,\cdots,7)$ can be written in terms of $_4f_3$.
\subsection{Difference equation $Rc^{(1)}(\bm{\alpha})$\ :\ an extension of $Rc^{(0)}(\hat{\bm{\alpha}})$}
Set 
\begin{align*}
  \bm{e}_{1}
  =
  \begin{pmatrix}
    1, 0, 0, 0\\
    1, 0, 0
  \end{pmatrix}
  ,\ 
  \bm{e}_{2}
  =
  \begin{pmatrix}
    0, 1, 0, 0\\
    0, 1, 0
  \end{pmatrix}
  ,\ 
  \bm{e}_{12}
  =
  \begin{pmatrix}
    1, 1, 0, 0\\
    1, 1, 0
  \end{pmatrix}.
\end{align*}
as in \S\ref{3termrelfor4F3}. We have
\begin{prp}[\cite{EO}]%, compare with Proposition \ref{AARW}]
\begin{align*}
  y_{i}(\bm{\alpha})
  &=
  U_{1}^{(1)}(\bm{\alpha})
  y_{i}(\bm{\alpha}+\bm{e}_{1})
  +
  V_{1}^{(1)}(\bm{\alpha})
  y_{i}(\bm{\alpha}+\bm{e}_{12})
  +
  \dfrac{1}
  {\alpha_{0}(\beta_{3}-\alpha_{0}-1)},
  %\label{contig_rel_4f3_seed1_nonterminate}
  \\
  y_{i}(\bm{\alpha})
  &=
  U_{2}^{(1)}(\bm{\alpha})
  y_{i}(\bm{\alpha}+\bm{e}_{2})
  +
  V_{2}^{(1)}(\bm{\alpha})
  y_{i}(\bm{\alpha}+\bm{e}_{12})
  +
  \dfrac{1}
  {\alpha_{1}(\beta_{3}-\alpha_{1}-1)}
 %\label{contig_rel_4f3_seed2_nonterminate}
\end{align*}
hold for any $i$ $(i=0, 1, \dots , 7)$.
Here,
\begin{align*}
  &
  U_{1}^{(1)}(\bm{\alpha})
  =
  \dfrac{-(\beta_{1}-\alpha_{1})(\beta_{1}+\beta_{2}-\alpha_{2}-\alpha_{3})}
  {\alpha_{0}(\beta_{3}-\alpha_{0}-1)},
  &
  &V_{1}^{(1)}(\bm{\alpha})
  =
  \dfrac{-(\beta_{2}-\alpha_{2})(\beta_{2}-\alpha_{3})}
  {\alpha_{0}(\beta_{3}-\alpha_{0}-1)},
  &
  \\
  &
  U_{2}^{(1)}(\bm{\alpha})
  =
  U_{1}^{(1)}(\bm{\alpha})
  |_{\alpha_{0}\leftrightarrow \alpha_{1}, \beta_{1}\leftrightarrow \beta_{2}},
  &
  &
  V_{2}^{(1)}(\bm{\alpha})
  =
  V_{1}^{(1)}(\bm{\alpha})
  |_{\alpha_{0}\leftrightarrow \alpha_{1}, \beta_{1}\leftrightarrow \beta_{2}}.
  &
\end{align*} 
\end{prp}
As Proposition \ref{AARW} led to Lemma \ref{q0} and Proposition \ref{Rc0}, we have
\begin{lmm}%[Compare with Lemma \ref{q0}]
\label{prp_diff_eq_4f3_2200220_non-homo}
Assume condition \eqref{balanced_condition_nonterminate}.
Set $(n; \bm{\alpha}):=-n\bm{e}_{12}+\bm{\alpha}$.
$y_{i}(n; \bm{\alpha}) \ ( i=0,1,\dots , 7)$ satisfies the non-homogeneous linear difference equation
\begin{gather*}
  D_{n}
  =
  q_{1}^{(1)}(n;{\bm{\alpha}})D_{n-1}
  +
  q_{2}^{(1)}(n;{\bm{\alpha}})D_{n-2}
  +
  q_{0}^{(1)}(n;{\bm{\alpha}}),
  %\label{diff_eq_4f3_2200220_non-homo}
\end{gather*}  
where 
\begin{align*} 
q_1^{(1)}({\bm{\alpha}})&=
  U^{(1)}_{1}({\bm{\alpha}})
  U^{(1)}_{2}({\bm{\alpha}}+{\bm{e}_{1}})
  +
  V^{(1)}_{1}({\bm{\alpha}})
  +
  \dfrac
  {U^{(1)}_{1}({\bm{\alpha}})V^{(1)}_{2}({\bm{\alpha}}+{\bm{e}_{1}})}
  {U^{(1)}_{1}({\bm{\alpha}}+{\bm{e}_{12}})}, 
     \\
q_2^{(1)}({\bm{\alpha}})&=
  \dfrac{ 
  -
    U^{(1)}_{1}({\bm{\alpha}})
  V^{(1)}_{1}({\bm{\alpha}}+{\bm{e}_{12}})
  V^{(1)}_{2}({\bm{\alpha}}+{\bm{e}_{1}})
}
{U^{(1)}_{1}({\bm{\alpha}}+\bm{e}_{12})},
\\
q_0^{(1)}({\bm{\alpha}})&=
\dfrac{(\beta_{1}-\alpha_{2}+1)(\beta_{1}-\alpha_{3}+1)(\beta_{1}+\beta_{2}-\alpha_{2}-\alpha_{3})}
{\alpha_{0}\alpha_{1}(\beta_{3}-\alpha_{0}-1)(\beta_{3}-\alpha_{1}-1)(\beta_{1}+\beta_{2}-\alpha_{2}-\alpha_{3}+2)}
+
\dfrac{\alpha_{0}\alpha_{1}-\beta_{1}(\beta_{1}+\beta_{2}-\alpha_{2}-\alpha_{3})}
{\alpha_{0}\alpha_{1}(\beta_{3}-\alpha_{0}-1)(\beta_{3}-\alpha_{1}-1)}.
\end{align*}
\end{lmm}
%This generalization of Proposition \ref{diffeqforterm4F3} yields the following:
\begin{prp}%[Compare with Proposition \ref{Rc0}]
\label{prp_diff_eq_4f3_2200220_homo}
Assume condition (\ref{balanced_condition_nonterminate}).
\begin{align*}
  D_{ij}^{(1)}(n;\bm{\alpha})
  &:=
  y_{i}(n;\bm{\alpha})
  -
  y_{j}(n;\bm{\alpha})\qquad i, j \in \{0, 1, 2, \dots ,7\}
\end{align*}
satisfies the homogeneous linear difference equation
\begin{gather*}
  Rc ^{(1)}(\bm{\alpha}):\
  D_{n}
  =
  q_{1}^{(1)}(n;{\bm{\alpha}})D_{n-1}
  +
  q_{2}^{(1)}(n;{\bm{\alpha}})D_{n-2}.
  \label{diff_eq_4f3_2200220_homo}
\end{gather*}  
The invariant of $Rc^{(1)}(\bm{\alpha})$ is
\begin{align*}
  H^{(1)}(n; \bm{\alpha})
  :=
  \dfrac{
    q_{1}^{(1)}(n;\bm{\alpha})
    q_{1}^{(1)}(n+1;\bm{\alpha})
  }{
    q_{2}^{(1)}(n+1;\bm{\alpha})
  }.
\end{align*}
\end{prp}
%We should remark $Rc^{(1)}(\bm{\alpha})$ can be regarded as an extension of $Rc^{(0)}(\hat{\bm{\alpha}})$.
Note that this invariant is a generalization of the former one: %In particular,
$$H^{(1)}(n; \bm{\alpha})|_{\alpha_{0}=0}=H^{(0)}(n; \hat{\bm{\alpha}}).$$
\begin{rmk}[Symmetry]
The homogeneous linear difference equation $Rc ^{(1)}(\bm{\alpha})$ has symmetries
\begin{align}
  \begin{split}
  \alpha_0\leftrightarrow \alpha_1\ ,
  \quad 
  \alpha_2\leftrightarrow \alpha_3\ ,
  \quad 
  \beta_1\leftrightarrow \beta_2\ ,
  \label{symmetry_trivial_4f3_2200220}
  \end{split}
  \\
  \begin{split}
  &\bm{\alpha}\rightarrow 
  \left( \alpha_{0}+1-\beta_{3},\alpha_{1}+1-\beta_{3},\alpha_{2}+1-\beta_{3},\alpha_{3}+1-\beta_{3},\beta_{1}+1-\beta_{3},\beta_{2}+1-\beta_{3},2-\beta_{3} \right)\ ,
  \label{symmetry_non-trivial1_4f3_2200220}
  \end{split}
  \\
  \begin{split}
  &\bm{\alpha}\rightarrow
  \left( \alpha_{0},\alpha_{0}+1-\beta_{3},\alpha_{0}+1-\beta_{1},\alpha_{0}+1-\beta_{2},\alpha_{0}+1-\alpha_{2},\alpha_{0}+1-\alpha_{3},\alpha_{0}+1-\alpha_{1} \right).
  \label{symmetry_non-trivial2_4f3_2200220}
  \end{split}
\end{align}
\label{symmetry_4f3_2200220}
\end{rmk}

\subsection{Local solutions at zero II: solving $Rc_{1}(A)$}
%In this subsection, we construct local solutions of $Z(A)$ at $x=0$ with exponent $A_0\pm \frac12$.
%In the previous section, we introduced $H^{(1)}(n; \bm{\alpha})$, which was an extension of $H^{(0)}(n; \hat{\bm{\alpha}})$.
We eventually find $\bm{\alpha}$ to solve $Rc_1(A)$.
\begin{prp}%[Compare with Proposition \ref{hatal}]
 The invariant $ H_{1}(n;A)$ of the equation $Rc_1(A)$ and the invariant  $H^{(1)}(n; {\bm{\alpha}})$ of the equation $Rc^{(1)}(\bm{\alpha})$ agree if and only if
  $  %\displaystyle
  \bm{\alpha} =(\alpha_{0},\alpha_{1},\alpha_{2},\alpha_{3};\beta_{1},\beta_{2},\beta_{3})$ is equal to
   \begin{align}
    \bm{\alpha}
    &=
    \left( 
      \dfrac{1}{2}, A_{0}+\dfrac{1}{2}, A_{+---}, A_{+-++};
      A_{+--+}+\dfrac{1}{2}, A_{+-+-}+\dfrac{1}{2}, A_{0}+1
    \right),
    \label{sol1_rc1}
    \\
   {\rm or}\quad &
    %\quad
   = \left( 
      \dfrac{1}{2}, A_{2}+\dfrac{1}{2}, A_{--+-}, A_{+-++};
      A_{--++}+\dfrac{1}{2}, A_{+-+-}+\dfrac{1}{2}, A_{2}+1
    \right),
    \label{sol2_rc1}
    \\
    {\rm or} \quad& 
   =    \left( 
      A_{0}+\dfrac{1}{2}, A_{2}+\dfrac{1}{2}, A_{+-++}, A_{+-+-};
      A_{+-++}+\dfrac{1}{2}, A_{+-+-}+\dfrac{1}{2}, A_{0}+A_{2}+1
    \right),
    \label{sol3_rc1}
    \\
   {\rm or}\quad &  
    =    \left( 
      A_{2}+\dfrac{1}{2}, A_{3}+\dfrac{1}{2}, A_{+-++}, A_{--++};
      A_{+-++}+\dfrac{1}{2}, A_{--++}+\dfrac{1}{2}, A_{2}+A_{3}+1
    \right),
    \label{sol4_rc1}
  \end{align}
  up to the symmetries
  \eqref{symmetry_Rc1},
  \eqref{symmetry_trivial_4f3_2200220},
  \eqref{symmetry_non-trivial1_4f3_2200220},
  \eqref{symmetry_non-trivial2_4f3_2200220}
 of the difference equations.
  \label{sol_rc1}
\end{prp}
Proof is parallel to that of Proposition \ref{hatal}, and goes as follows: The numerator of $H_1(n; A)- H^{(1)}(n; \bm{\alpha})$ is a
polynomial in $n$ of degree 14. The coefficient of the top term
decomposes as
\begin{align*}  {\rm constant}\times(\alpha_{0}+\alpha_{1}-\beta_{1}-\beta_{2}+1-A_{1})
  (\alpha_{0}+\alpha_{1}-\beta_{1}-\beta_{2}+1+A_{1}).\end{align*}
If we set 
$\displaystyle \alpha_{0}=-\alpha_{1}+\beta_{1}+\beta_{2}-1\pm A_{1} $, 
then the second top term decomposes as
\begin{align*}   {\rm constant}\times(2A_{1}-1)(2A_{1}+1)
  (\alpha_{2}+\alpha_{3}-\beta_{1}-\beta_{2}+1).\end{align*}
Since we assume $A_{1}$ is a free parameter, we have
\begin{align}  \alpha_{2}+\alpha_{3}-\beta_{1}-\beta_{2}+1=0.
  \label{eq_lmm1_proof_prop11.6}
\end{align}
Eliminating $\beta_2$ and $\beta_3$ from $H^{(1)}(n; \bm{\alpha})$ by making use of the condition (\ref{balanced_condition_nonterminate}) and the relation (\ref{eq_lmm1_proof_prop11.6}), we get an expression
\begin{align*}
  H^{(1)}(n; \bm{\alpha})
  =\dfrac{
  -\gamma^{(1)} +O\left( \left( n-\frac12 \right)^{2} \right)
}{
  \gamma^{(1)} +O\left( \left( n-\frac12 \right)^{2} \right)
  }, \end{align*}
where
\begin{align*}
  \gamma^{(1)} :=(2\alpha_0-1)^2(2\alpha_1-1)^2
  (2\alpha_2-2\beta_1+1)^2(2\alpha_3-2\beta_1+1)^2.
  %\label{eq_lmm2_proof_prop11.6}
\end{align*}
Now, assume $\gamma^{(1)}\neq 0$. Then, $H^{(1)}\left(\frac12; \bm{\alpha}\right)=-1,$ whereas $H_{1}\left(\frac12; A\right)$ is parameter dependent:
\begin{align*}
H_{1}\left(\frac12; A\right)
=
\dfrac{-\left( A_{0}^{2}-A_{1}^{2}+A_{2}^{2}+A_{3}^{2}-1 \right)^{2}}{4A_{0}^{2}A_{2}^{2}A_{3}^{2}}.
\end{align*} So it does not happen that $H_{1}\left(n; A\right)=H^{(1)}(n;\bm{\alpha})$ in the case $\gamma^{(1)}\neq 0$. %is parameter dependent:
%$$H_{1}\left(\frac12; A\right)=-\frac{( A_{0}^{2}-A_{1}^{2}+A_{2}^{2}+A_{3}^{2}-1)^{2}}{4A_{0}^{2}A_{2}^{2}A_{3}^{2}}=-1,$$
Thus $\gamma^{(1)}$ must be 0.
Thanks to the symmetry (\ref{symmetry_trivial_4f3_2200220}) of $Rc^{(1)}(\bm{\alpha})$, we have only to consider two cases:
$$  \text{Case 1}:\
  \alpha_{0}=\dfrac{1}{2}
  \quad{\rm and}\quad
  \text{Case 2}:\
  \beta _{1}=\alpha _{2}+\dfrac{1}{2}.$$
\medskip
\noindent
\underline{Case 1}: 
Eliminating $\beta_2$ and $\beta_3$ from $H^{(1)}(n; \bm{\alpha})$ by the condition (\ref{balanced_condition_nonterminate}) and the relation (\ref{eq_lmm1_proof_prop11.6}), we get an expression
\begin{align*}
  &H^{(1)}\left(n; \bm{\alpha}\right)=
  \dfrac{-n(n-1)\epsilon^{(1)}_{1}(n)\epsilon^{(1)}_{1}(n+1)}
{(n-\alpha_{1})(n+\alpha_{1}-1)
  (n+\alpha_{2}-\beta_{1})(n-\alpha_{2}+\beta_{1}-1)
  (n+\alpha_{3}-\beta_{1})(n-\alpha_{3}+\beta_{1}-1)},\\
&\text{where}\quad
\epsilon^{(1)}_{1}(n)
=2\left(
n^{2}-2n-\alpha_1\alpha_2-\alpha_1\alpha_3+\alpha_2\alpha_3+\alpha_2\beta_1+\alpha_3\beta_1-\beta_1^2+\alpha_1-\alpha_2-\alpha_3+\beta_1+1
\right).\end{align*}
Equating this and $H_{1}(n; A)$, we conclude that $\bm{\alpha}$ satisfies
\begin{itemize}
  \item
    $\displaystyle
    \epsilon^{(1)}_{1}(0)=-A_{0}^{2}+A_{1}^{2}-A_{2}^{2}-A_{3}^{2}+\frac52,
  $
  \item
    $\displaystyle
  \left\{
    \pm \left( \alpha_{1}-\dfrac{1}{2}\right)-\dfrac{1}{2},\ 
    \pm \left( -\alpha_{2}+\beta_{1}-\dfrac{1}{2}\right)-\dfrac{1}{2},\ 
    \pm \left( -\alpha_{3}+\beta_{1}-\dfrac{1}{2}\right)-\dfrac{1}{2}
  \right\}
  \\[2pt]
  =
  \left\{
   \pm A_{0}-\frac12,\ \pm A_{2}-\frac12,\ \pm A_{3}-\frac12
   \right\}.  $
\end{itemize}
These lead to (\ref{sol1_rc1}) and (\ref{sol2_rc1}), up to the symmetries (\ref{symmetry_Rc1}) and (\ref{symmetry_trivial_4f3_2200220}).
\par\medskip\noindent
\underline{Case 2}:  In the same way we have
\begin{align*}
&H^{(1)}\left(n; \bm{\alpha}\right)
=
\dfrac{-n(n-1)\epsilon^{(1)}_{2}(n)\epsilon^{(1)}_{2}(n+1)}
{
  (n-\alpha_{0})(n+\alpha_{0}-1)
  (n-\alpha_{1})(n+\alpha_{1}-1)
  \left( n-\alpha_{2}+\alpha_{3}-\frac{1}{2} \right)
  \left( n+\alpha_{2}-\alpha_{3}-\frac{1}{2} \right)},\\
&\text{where}\quad
\epsilon^{(1)}_2(n)
=2\left(
  n^{2}-2n+\alpha_0\alpha_1-\alpha_0\alpha_2-\alpha_0\alpha_3-\alpha_1\alpha_2-\alpha_1\alpha_3+2\alpha_2\alpha_3+\frac{\alpha_{0}}{2}+\frac{\alpha_{1}}{2}+1
\right),\end{align*}
and the conditions
\begin{itemize}
  \item
    $\displaystyle
  \epsilon^{(1)}_2(0)=-A_{0}^{2}+A_{1}^{2}-A_{2}^{2}-A_{3}^{2}+\frac52,
  $
  \item
    $\displaystyle
  \left\{
    \pm \left( \alpha_{0}-\dfrac{1}{2}\right)-\dfrac{1}{2},\ 
    \pm \left( \alpha_{1}-\dfrac{1}{2}\right)-\dfrac{1}{2},\ 
    \pm \left( \alpha_{2}-\alpha_{3}\right)-\dfrac{1}{2}
  \right\}
  =
  \left\{
   \pm A_{0}-\frac12,\ \pm A_{2}-\frac12,\ \pm A_{3}-\frac12
  \right\},
  $
\end{itemize}
which lead to (\ref{sol3_rc1}) and (\ref{sol4_rc1}), up to the symmetries
(\ref{symmetry_Rc1}) and (\ref{symmetry_trivial_4f3_2200220}). \par\noindent
These complete the proof of the proposition.
\par\medskip\noindent
This proposition together with Corollary \ref{solspace} implies 
\begin{prp} %[Compare with Proposition \ref{Rc0}]
\begin{align*}
  &(\text{Solutions of $Rc_{1}(A)$})
  \\
  &=
  \dfrac{
    \left( \frac{1}{2} \right)_{n}
  }{
    n!
  }
  (\text
    {Solutions of $Rc^{(1)}(\bm{\alpha})$
    with $\bm{\alpha}$ as \eqref{sol1_rc1}
    }
  )
  \\
  &=
  \dfrac{
    \left( \frac{1}{2} \right)_{n}
    \left( -A_{2}+\frac{1}{2} \right)_{n}
    \left( A_{2}+\frac{1}{2} \right)_{n}
  }{
    n!
    \left( -A_{0}+\frac{1}{2} \right)_{n}
    \left( A_{0}+\frac{1}{2} \right)_{n}
  }
  (\text
    {Solutions of $Rc^{(1)}(\bm{\alpha})$
    with $\bm{\alpha}$ as \eqref{sol2_rc1}
    }
  )
  \\
  &=
  \dfrac{
    \left( -A_{2}+\frac{1}{2} \right)_{n}
    \left( A_{2}+\frac{1}{2} \right)_{n}
  }{
    n!
    \left( \frac{1}{2} \right)_{n}
  }
  (\text
    {Solutions of $Rc^{(1)}(\bm{\alpha})$
    with $\bm{\alpha}$ as \eqref{sol3_rc1}
    }
  )
  \\
  &=
  \dfrac{
    \left( -A_{2}+\frac{1}{2} \right)_{n}
    \left( A_{2}+\frac{1}{2} \right)_{n}
    \left( -A_{3}+\frac{1}{2} \right)_{n}
    \left( A_{3}+\frac{1}{2} \right)_{n}
  }{
    n!
    \left( \frac{1}{2} \right)_{n}
    \left( -A_{0}+\frac{1}{2} \right)_{n}
    \left( A_{0}+\frac{1}{2} \right)_{n}
  }
  (\text
    {Solutions of $Rc^{(1)}(\bm{\alpha})$
    with $\bm{\alpha}$ as \eqref{sol4_rc1}
    }
  ).
\end{align*}
\label{sol_space_rc1}
\end{prp}
So far we got many solutions of the difference equation $Rc_1(A)$, whose solution space is two dimensional. Among these, there are many linearly independent pairs, but it is not so obvious to pick two independent ones.
Set 
  \begin{align*}
    W_{1,ij}^{(1)}(n; A)
    &:=
  \dfrac{
    \left( \frac{1}{2} \right)_{n}
  }{
    n!
  }
    D_{ij}^{(1)}(n;\bm{\alpha})
    \Bigl|_{\bm{\alpha}=(\ref{sol1_rc1})}
    ,\
  \\
  W_{2,ij}^{(1)}(n; A)
  &:=
  \dfrac{
    \left( \frac{1}{2} \right)_{n}
    \left( -A_{2}+\frac{1}{2} \right)_{n}
    \left( A_{2}+\frac{1}{2} \right)_{n}
  }{
    n!
    \left( -A_{0}+\frac{1}{2} \right)_{n}
    \left( A_{0}+\frac{1}{2} \right)_{n}
  }
    D_{ij}^{(1)}(n;\bm{\alpha})
    \Bigl|_{\bm{\alpha}=(\ref{sol2_rc1})}
    ,\
  \\
  W_{3,ij}^{(1)}(n; A)
  &:=
  \dfrac{
    \left( -A_{2}+\frac{1}{2} \right)_{n}
    \left( A_{2}+\frac{1}{2} \right)_{n}
  }{
    n!
    \left( \frac{1}{2} \right)_{n}
  }
    D_{ij}^{(1)}(n;\bm{\alpha})
    \Bigl|_{\bm{\alpha}=(\ref{sol3_rc1})}
    ,\
  \\
  W_{4,ij}^{(1)}(n; A)
  &:=
  \dfrac{
    \left( -A_{2}+\frac{1}{2} \right)_{n}
    \left( A_{2}+\frac{1}{2} \right)_{n}
    \left( -A_{3}+\frac{1}{2} \right)_{n}
    \left( A_{3}+\frac{1}{2} \right)_{n}
  }{
    n!
    \left( \frac{1}{2} \right)_{n}
    \left( -A_{0}+\frac{1}{2} \right)_{n}
    \left( A_{0}+\frac{1}{2} \right)_{n}
  }
    D_{ij}^{(1)}(n;\bm{\alpha})
    \Bigl|_{\bm{\alpha}=(\ref{sol4_rc1})},
 \end{align*}
and 
$$  \Omega:=\{(k,i,j)\ | 
  \ k=1, 2, 3, 4;\ i, j=0, 1, \dots, 7;\ i\neq j\}.$$
  We have
\begin{thm}%[Compare with Theorem \ref{solZ4F3}]
  \label{sol_space_Z_0_basis}
For  $(k,i,j)\in \Omega$,
  \begin{align*}
    f_{k,ij}^{(0,A_{0}-1/2)}(A;x)
    :=
    x^{A_{0}-\frac12}
    \sum_{n=0}^{\infty}
    W_{k,ij}^{(1)}(n; A)x^{n}
  \end{align*}
  is a solution of $Z(A)$ at $x=0$ with exponent $A_0-\frac12$.
\end{thm}
Examples of linearly independent pairs:
$$\{f_{1,01}^{(0,A_{0}-1/2)},f_{1,02}^{(0,A_{0}-1/2)}\},\ \{f_{1,03}^{(0,A_{0}-1/2)},f_{1,54}^{(0,A_{0}-1/2)}\},\ \{f_{2,03}^{(0,A_{0}-1/2)},f_{2,54}^{(0,A_{0}-1/2)}\}.$$

\section{Local solutions of $Q(A)$}Let $\displaystyle P_{Q}(A;x) \left(  = P_{Q}(A_{0}, A_{1}, A_{2}, A_{3}; x) \right) $ denote the space of solutions of $Q(A)$ appeared in \S 6.2.
Remark \ref{symmetryQ} implies the following equivalence.
\begin{prp} 
\begin{align}
P_{Q}(A_{0},A_{1},A_{2},A_{3};x)
&=
P_{Q}(\pm A_{0},\pm A_{1},A_{2}, \pm A_{3};x)
\label{solution_symmetryQ_trivial}
\\
&=
P_{Q}(\pm A_{1},\pm A_{0},A_{2}, \pm A_{3};1-x)
\label{solution_symmetryQ_1}\\
&=
\left( -\dfrac{1}{x} \right)^{1+2A_{2}}
P_{Q}\left(\pm A_{3},\pm A_{1},A_{2}, \pm A_{0};\dfrac{1}{x}\right)
\label{solution_symmetryQ_infty}
\\
&=
\left( 1-x \right)^{-1-2A_{2}}
P_{Q}\left(\pm A_{0},\pm A_{3},A_{2}, \pm A_{1};\dfrac{x}{x-1}\right).
\label{solution_symmetryQ_pfaff}
\end{align}
\end{prp}
To help understand the following propositions, we tabulate several expressions of the Riemann scheme of $Q(A)$:
\begin{align*}
 \left\{
  \begin{matrix}
   x=0&x=1&x=\infty\\
   0&0&1+2A_2\\
   -A_0-A_2&-A_1-A_2&1+A_2-A_3\\
   A_0-A_2&A_1-A_2&1+A_2+A_3
  \end{matrix}
 \right\}
%\end{align*}\begin{align*}
    &=
  \left(-\dfrac{1}{x}\right)^{1+2A_{2}}
 \left\{
  \begin{matrix}
    \frac{1}{x}=\infty &\frac{1}{x}=1&\frac{1}{x}=0\\[3pt]
    1+2A_{2}&0& 0\\
    1+A_{2}-A_0&-A_1-A_2& -A_3-A_{2}\\
    1+A_{2}+A_{0}&A_1-A_2& A_{3}-A_2
  \end{matrix}
 \right\}
 \\[5pt]
  &=
  (1-x)^{-(1+2A_{2})}
  \begin{Bmatrix}
    \frac{x}{x-1}=0 & \frac{x}{x-1}=\infty & \frac{x}{x-1}=1\\[3pt]
    0            & 1+2A_{2}         &  0        \\
    -A_{0}-A_{2} & 1+A_{2}-A_{1}    & -A_{3}-A_{2}  \\
    A_{0}-A_{2}  & 1+A_{2}+A_{1}    & A_{3}-A_{2}
  \end{Bmatrix}.
\end{align*}
Thanks to the identities
(\ref{solution_symmetryQ_1}) and 
(\ref{solution_symmetryQ_infty}), 
three linearly independent local solutions at $x=0$ give those at other singular points.
We find a holomorphic solution $f_{Q}^{(0,0)}(A;x)$
and a solution $f_{Q}^{(0,\pm)}(A;x)$ of local exponent $\pm A_{0}-A_{2}$ at $x=0$ as follows.

\subsection{Holomorphic solution  $f_{Q}^{(0,0)}(A;x)$ to $Q(A)$ at $x=0$}
Set
\begin{gather*}
f_{Q}^{(0,0)}(A;x)
=
\sum_{n=0}^{\infty} C_{n}\, x^{n},
\quad
C_{0}=1.
%\label{f_Q_0_0}
\end{gather*}
The coefficients $C_{n}$ 
satisfy the recurrence relation
\begin{align*}
  RcQ(A):
\begin{array}{ll}
  &C_{n}
  =
  \dfrac
  {\left(n+A_2-\frac{1}{2}\right) 
   \left( 2n^{2}+2(2A_{2}-1)n-A_{0}^{2}+A_{1}^{2}+A_{2}^{2}-A_{3}^{2}-2A_{2}+1\right)
  }{n(n+A_{0}+A_{2})(n-A_{0}+A_{2})}
  C_{n-1}
  \\
  &\phantom{C_{n}}-
  \dfrac
  {(n+2A_{2}-1)(n+A_{2}+A_{3}-1)(n+A_{2}-A_{3}-1)
  }{n(n+A_{0}+A_{2})(n-A_{0}+A_{2})}
  C_{n-2},
\end{array}
\end{align*}
whose invariant will be called $H_{Q}(n; A)$.
\begin{rmk}[Symmetry]
  $RcQ(A)$ is invariant under
  \begin{align}
    A_{j} \rightarrow -A_{j} 
    \quad(j=0, 1, 3).
    \label{symmetryRcQ}
  \end{align}
\end{rmk}
\begin{prp}
  Let $A$ be generic.
  Then, the two invariants $H_{Q}(n; A)$ and $H^{(0)}(n; \hat{\bm{\alpha}})$ 
  of the difference equations
  $RcQ(A)$ and $Rc^{(0)}(\hat{\bm{\alpha}})$
  agree
  as functions in $n$ if and only if
  \begin{align*}
   \hat{\bm{\alpha}}
    =
  \begin{pmatrix}
    0, \alpha_{1},\alpha_{2},\alpha_{3}\\
    \beta_{1},\beta_{2},\beta_{3}
  \end{pmatrix}
    &=
    \begin{pmatrix}
      0, A_{0}-A_{2}, A_{+-+-}, A_{+-++}\\A_{+---},A_{+--+},1+A_{0}+A_{2}
  \end{pmatrix},
  %\label{sol1_rcq}
  \\
 {\rm or}\quad  &\ %\phantom{=}
  =\begin{pmatrix}
    0, -A_{2}+A_{3}, A_{--++}, A_{+-++}\\
    A_{---+}, A_{+--+}, 1+A_{2}+A_{3}
  \end{pmatrix}
  %\label{sol2_rcq}
  ,\\
{\rm or}\quad &\ %\phantom{=}
  =\begin{pmatrix}
    0, -A_{2}, A_{--++}, A_{+-+-}\\
    A_{----},A_{+--+},1+A_{2}
  \end{pmatrix},
  %\label{sol3_rcq}
  \end{align*}
  up to the symmetries of two difference equations,
  that is, (\ref{symmetryRcQ}) and Remark \ref{symAc0}.
\label{solution_rcQ}
\end{prp}
Doing the same as we got Theorem \ref{solZ4F3} from Proposition \ref{sol_rc0},  we can obtain the following theorem:
\begin{thm} The holomorphic solution  $f_{Q}^{(0,0)}(A;x)$ to  $Q(A)$ at $x=0$ has the following expressions:
\begin{align}
  &f_{Q}^{(0,0)}(A;x)
  \notag
  \\
  &=
  \sum_{n=0}^{\infty}
    x^{n}
  \dfrac{
    \left( 1+2A_{2} \right)_{n}
    \left( A_{-+++} \right)_{n}
    \left( A_{-++-} \right)_{n}
  }{
    \left( 1+A_{2} \right)_{n}
    (1-A_{0}+A_{2})_{n}
    n!
  }
  \pFq{4}{3}{
      -n,
      A_{0}-A_{2}-n,
      A_{+-+-},
      A_{+-++}
    }{
      A_{+---}-n,
      A_{+--+}-n,
      1+A_{0}+A_{2}
    }{1}
    \label{expression_q_0_0_1}
    \\
  &=
  \sum_{n=0}^{\infty}
    x^{n}
  \dfrac{
    \left( 1+2A_{2} \right)_{n}
    \left( 1+A_{2}+A_{3} \right)_{n}
    \left( A_{+++-} \right)_{n}
    \left( A_{-++-} \right)_{n}
  }{
    \left( 1+A_{2} \right)_{n}
    (1-A_{0}+A_{2})_{n}
    (1+A_{0}+A_{2})_{n}
    n!
  }
  \pFq{4}{3}{
      -n,
      -A_{2}+A_{3}-n,
      A_{--++},
      A_{+-++}
    }{
      A_{---+}-n,
      A_{+--+}-n,
      1+A_{2}+A_{3}
    }{1}
    \label{expression_q_0_0_2}
    \\
  &=
    \sum_{n=0}^{\infty}
    x^{n}
  \dfrac{
    \left( 1+2A_{2} \right)_{n}
    \left( A_{++++} \right)_{n}
    \left( A_{-++-} \right)_{n}
  }{
    (1-A_{0}+A_{2})_{n}
    (1+A_{0}+A_{2})_{n}
    n!
  }
  \pFq{4}{3}{
      -n,
      -A_{2}-n,
      A_{--++},
      A_{+-+-}
    }{
      A_{----}-n,
      A_{+--+}-n,
      1+A_{2}
    }{1}.
    \label{expression_q_0_0_3}
  \end{align}
  \label{expression_q_0_0}
\end{thm}
\begin{rmk}
With the help of the transformation formula \eqref{trans_4F3_1}, one of the three
\eqref{expression_q_0_0_1},
\eqref{expression_q_0_0_2},
\eqref{expression_q_0_0_3}
implies the other two.
\end{rmk}

\begin{rmk}
  By using the formula \eqref{2F12F14F3},
  the expressions 
  \eqref{expression_q_0_0_1},
  \eqref{expression_q_0_0_2},
  \eqref{expression_q_0_0_3}
can be written also as
\begin{align}
  &f_{Q}^{(0,0)}(A;x)
  \notag
  \\
  &=
    \sum_{n=0}^{\infty}
      x^{n}
  \dfrac{
    \left( 1+2A_{2} \right)_{n}
  }{
    \left( 1+A_{2}\right)_{n}
  }
  [Z^{n}]
  \pFq{2}{1}{A_{-+++},A_{-++-}}
  {1-A_{0}+A_{2}}{Z}
  \pFq{2}{1}{A_{+-+-},A_{+-++}}
  {1+A_{0}+A_{2}}{Z}
    \label{expression_q_0_0_1_another}
  \\  
  &=
    \sum_{n=0}^{\infty}
      x^{n}
  \dfrac{
    (1+2A_{2})_{n}
    (1+A_{2}+A_{3})_{n}
    (1+A_{2}-A_{3})_{n}
  }{
    (1+A_{2})_{n}
    (1-A_{0}+A_{2})_{n}
    (1+A_{0}+A_{2})_{n}
  }
  [Z^{n}]
  \pFq{2}{1}{A_{+++-},A_{-++-}}
  {1+A_{2}-A_{3}}{Z}
  \pFq{2}{1}{A_{--++},A_{+-++}}
  {1+A_{2}+A_{3}}{Z}
    \label{expression_q_0_0_2_another}
  \\  
  &=
  \sum_{n=0}^{\infty}
      x^{n}
  \dfrac{
    \left( 1+A_{2}\right)_{n} \left( 1+2A_{2} \right)_{n} }{
      (1-A_{0}+A_{2})_{n}(1+A_{0}+A_{2})_{n}
  }
  [Z^{n}]
  \pFq{2}{1}{A_{++++},A_{-++-}}{1+A_{2}}{Z}
  \pFq{2}{1}{A_{--++},A_{+-+-}}{1+A_{2}}{Z}.
    \label{expression_q_0_0_3_another}
    \end{align}
    \label{expression_q_0_0_another}
 \end{rmk}
\begin{rmk}\label{trivialsymmetryfQ}
  Applying the trivial symmetry \eqref{solution_symmetryQ_trivial}) to $f_{Q}^{(0,0)}(A;x)$, we have
\begin{align*}
  f_{Q}^{(0,0)}(A;x)
  =
  f_{Q}^{(0,0)}(\pm A_{0}, \pm A_{1},A_{2}, {\pm} A_{3};x).
\end{align*}
\end{rmk}

\subsection{Local  solution  to $Q(A)$ at $x=0$ with exponent $A_{0}-A_{2}$}
Set
\begin{gather*}
f_{Q}^{(0,+)}(A;x)
=
x^{A_{0}-A_{2}}
\sum_{n=0}^{\infty}
C_{n}x^{n},\qquad C_0=1.
%\label{f_Q_0_+}
\end{gather*}
The coefficients $C_n$ satisfy the recurrence relation:
\begin{align*}
\begin{array}{ll}
&
C_n
=
\dfrac
{\left(n+A_{0}-\frac{1}{2}\right)\left(2n^{2}+2(2A_{0}-1)n+A_{0}^2+A_{1}^2-A_{2}^2-A_{3}^2-2A_{0}+1\right)}
{n\left(n+2A_{0}\right)\left(n+A_{0}-A_{2}\right)}
C_{n-1}
\\
& 
\phantom{C_n}
-
\dfrac
{\left(n+A_{0}+A_{2}-1\right)\left(n+A_{0}+A_{3}-1\right)\left(n+A_{0}-A_{3}-1\right)}
{n\left(n+2A_{0}\right)\left(n+A_{0}-A_{2}\right)}
C_{n-2}.
\end{array}
\end{align*}
The same argument as in the  previous section leads to
\begin{thm}
The local solution   $f_{Q}^{(0,+)}(A;x)$ to $Q(A)$  at $x=0$ with exponent $A_{0}-A_{2}$ has the following expressions:
  \begin{align}
    &f_{Q}^{(0,+)}(A;x)
    \notag
    \\
    &=
      x^{A_0-A_{2}}
    \sum_{n=0}^{\infty}
      x^{n}
  \dfrac{
    \left( 1+A_{0}+A_{2} \right)_{n}
  }{
    \left( 1+A_{0}\right)_{n}
  }
  [Z^{n}]
  \pFq{2}{1}{A_{+--+},A_{+---}}{1+A_{0}-A_{2}}{Z}
  \pFq{2}{1}{A_{++++},A_{+++-}}{1+A_{0}+A_{2}}{Z}
  \label{expression_q_0_+_1_another}
  \\  
    &=
      x^{A_0-A_{2}}
    \sum_{n=0}^{\infty}
      x^{n}
  \dfrac{
    (1+A_{0}-A_{3})_{n}
    (1+A_{0}+A_{3})_{n}
  }{
    (1+A_{0})_{n}
    (1+A_{0}-A_{2})_{n}
  }
  [Z^{n}]
  \pFq{2}{1}{A_{+--+},A_{+-++}}{1+A_{0}+A_{3}}{Z}
  \pFq{2}{1}{A_{++--},A_{+++-}}{1+A_{0}-A_{3}}{Z}
  \label{expression_q_0_+_2_another}
  \\  
  &=
      x^{A_0-A_{2}}
    \sum_{n=0}^{\infty}
      x^{n}
  \dfrac{
    (1+A_{0})_{n}
  }{
    (1+A_{0}-A_{2})_{n}
  }
  [Z^{n}]
  \pFq{2}{1}{A_{+-+-},A_{+--+}}{1+A_{0}}{Z}
  \pFq{2}{1}{A_{++--},A_{++++}}{1+A_{0}}{Z}.
  \label{expression_q_0_+_3_another}
    \end{align}
  \label{expression_q_0_+_another}
\end{thm}
\begin{rmk}
 As in the previous section,  the formula  \eqref{2F12F14F3} changes the expressions 
  \eqref{expression_q_0_+_1_another} ,
  \eqref{expression_q_0_+_2_another}
  and
  \eqref{expression_q_0_+_3_another}
  into those in terms of  $_4F_3(1)$.
\end{rmk}
\begin{rmk}
  Applying the trivial symmetry \eqref{solution_symmetryQ_trivial} to $f_{Q}^{(0,+)}(A;x)$, we have
\begin{align*}
  f_{Q}^{(0,+)}(A;x)
  =
  f_{Q}^{(0,+)}(A_{0}, \pm A_{1},A_{2}, {\pm} A_{3};x).
\end{align*}
\label{expression_q_0_+_symmetry}
\end{rmk}
%\newpage
\subsection{Pfaff transforms of the solutions of $Q(A)$}
From formula \eqref{solution_symmetryQ_pfaff} we get the following Pfaff transformations:
\begin{prp}
  \begin{align*}
    f_{Q}^{(0,0)}(A;x)
    &=
  (1-x)^{-(1+2A_{2})}
    f_{Q}^{(0,0)}
    \left( {\pm}A_{0},\pm A_{3},A_{2},\pm A_{1};\frac{x}{x-1}\right),
  \\
    f_{Q}^{(0,+)}(A;x)
    &=
  C\times
  (1-x)^{-(1+2A_{2})}
    f_{Q}^{(0,+)}
    \left( A_{0},\pm A_{3},A_{2},\pm A_{1};\frac{x}{x-1}\right),
  \end{align*}
where  $C$ is a constant depending on the choice of the branch of $x^{A_0-A_2}$.% around $x=0$.
\label{pfaff_trans_q_0}
\end{prp}

\subsection{Local solutions of  $Q(A)$}
Applying the symmetries
(\ref{solution_symmetryQ_trivial}),
(\ref{solution_symmetryQ_1}) and (\ref{solution_symmetryQ_infty}) 
to local solutions $f_{Q}^{(0,0)}(A; x)$ and $f_{Q}^{(0,+)}(A; x)$ of $Q(A)$,
we obtain series expressions of other local solutions.
\begin{prp} Local solutions of $Q(A)$ are tabulated in Table $1$.
\begin{table}[htb]
\begin{center}
  \begin{tabular}{c|l||l}
    names of solutions& (point,  exponent) & \ series expression\\[2mm]
    $f_{Q}^{(0,0)}(A;x)$ & $(0, 0)$ 
    & \ Theorem \ref{expression_q_0_0}, Remark \ref{expression_q_0_0_another}, \ref{trivialsymmetryfQ}, Proposition \ref{pfaff_trans_q_0}
    \\[1mm]
    $f_{Q}^{(0,+)}(A;x)$ & $(0, A_{0}-A_{2})$
    & \ Theorem \ref{expression_q_0_+_another},  Remark \ref{expression_q_0_+_symmetry},  Proposition \ref{pfaff_trans_q_0}
    \\[1mm]
    $f_{Q}^{(0,-)}(A;x)$ & $(0, -A_{0}-A_{2})$ 
    &\ 
    $\displaystyle
    f_{Q}^{(0,+)}(-A_{0}, \pm A_{1}, A_{2}, {\pm}A_{3};x)$
    \\[1mm]
    $f_{Q}^{(1,0)}(A;x)$ & $(1, 0)$
    &\ 
    $\displaystyle
    f_{Q}^{(0,0)}\left(\pm A_{1}, \pm A_{0}, A_{2}, \pm A_{3}; 1-x\right)$
    \\[1mm]
    $f_{Q}^{(1,+)}(A;x)$ & $(1, A_{1}-A_{2})$
    &\ 
    $\displaystyle
    f_{Q}^{(0,+)}\left( A_{1}, \pm A_{0}, A_{2}, \pm A_{3}; 1-x\right)$
    \\[1mm]
    $f_{Q}^{(1,-)}(A;x)$ & $(1, -A_{1}-A_{2})$
    &\ 
    $\displaystyle
    f_{Q}^{(0,+)}\left( -A_{1}, \pm A_{0}, A_{2}, \pm A_{3}; 1-x\right)$
    \\[1mm]
    $f_{Q}^{(\infty,0)}(A;x)$ & $(\infty, 2A_{2}+1)$
    &\ 
    $\displaystyle
    \left(-\dfrac{1}{x}\right)^{1+2A_{2}}
    f_{Q} ^{(0,0)}
    \left(\pm A_{3}, \pm A_{1},A_{2},\pm A_{0};\dfrac{1}{x}\right)$
    \\[1mm]
    $f_{Q}^{(\infty,+)}(A;x)$ & $(\infty, A_{2}+A_{3}+1)$
    &\ 
    $\displaystyle
    \left(-\dfrac{1}{x}\right)^{1+2A_{2}}
    f_{Q}^{(0,+)}\left( A_{3}, \pm A_{1},A_{2},\pm A_{0};\dfrac{1}{x}\right)$
    \\[1mm]
    $f_{Q}^{(\infty,-)}(A;x)$ & $(\infty, A_{2}-A_{3}+1)$
    &\ 
    $\displaystyle
    \left(-\dfrac{1}{x}\right)^{1+2A_{2}}
    f_{Q}^{(0,+)}\left( -A_{3}, \pm A_{1},A_{2},\pm A_{0};\dfrac{1}{x}\right)$
    \\
  \end{tabular}
  \caption{Names of local solutions of $Q(A)$ and their series expressions}
\end{center}
\end{table}
\label{expression_q_0_-}
\end{prp}
%%%%%%%%%%%%%%%%%%%%%%%%%%%%%%%%
\section{Local solutions of the Dotsenko-Fateev equation}
The Dotsenko-Fateev equation $S$ is obtained (Propositions \ref{QtoR} ) from the equation $Q(A)$ by the change of unknown 
\begin{align}
  z=x^{-A_0-A_2}(x-1)^{-A_1-A_2}w,
  \label{relation_q_df}
\end{align}
where $z$: solution of $Q(A)$, $w$: solution of $S$, 
and the parameter change 
given in Propositions \ref{RtoDF}. So we get expressions of local solutions of Dotsenko-Fateev  (D-F  for short) equation from those obtained in the previous section.

Table \ref{correspond_solution_Q_DF} tabulates the local solutions of $Q(A)$ appeared in the previous section, and names $f_{DF}^{(*,*)}(a,b,c,g;x)$ of the corresponding solutions of the 
Dotsenko-Fateev equation.
\begin{table}[htb]
\begin{center}
  \begin{tabular}{c|l||c|l}
    solutions of $Q(A)$& (point,  exponent) &  solutions of D-F equation & (point, exponent) \\[2mm]
    $f_{Q}^{(0,0)}(A;x)$ & $(0, 0)$ &  $f_{DF}^{(0,1)}(a, b, c, g; x)$ & $(0, a+c+1)$\\[1mm]
    $f_{Q}^{(0,+)}(A;x)$ & $(0, A_{0}-A_{2})$ & $f_{DF}^{(0,2)}(a, b, c, g; x)$ & $(0, 2a+2c+g+2)$\\[1mm]
    $f_{Q}^{(0,-)}(A;x)$ & $(0, -A_{0}-A_{2})$ & $f_{DF}^{(0,0)}(a, b, c, g; x)$ & $(0, 0)$\\[1mm]
    $f_{Q}^{(1,0)}(A;x)$ & $(1, 0)$ &  $f_{DF}^{(1,1)}(a, b, c, g; x)$ & $(1, b+c+1)$\\[1mm]
    $f_{Q}^{(1,+)}(A;x)$ & $(1, A_{1}-A_{2})$ & $f_{DF}^{(1,2)}(a, b, c, g; x)$ & $(1, 2b+2c+g+2)$\\[1mm]
    $f_{Q}^{(1,-)}(A;x)$ & $(1, -A_{1}-A_{2})$ & $f_{DF}^{(1,0)}(a, b, c, g; x)$ & $(1, 0)$\\[1mm]
    $f_{Q}^{(\infty,0)}(A;x)$ & $(\infty, 2A_{2}+1)$ &  $f_{DF}^{(\infty,1)}(a, b, c, g; x)$ & $(\infty, -a-b-2c-g-1)$\\[1mm]
    $f_{Q}^{(\infty,+)}(A;x)$ & $(\infty, A_{2}+A_{3}+1)$ & $f_{DF}^{(\infty,0)}(a, b, c, g; x)$ & $(\infty, -2c)$\\[1mm]
    $f_{Q}^{(\infty,-)}(A;x)$ & $(\infty, A_{2}-A_{3}+1)$ & $f_{DF}^{(\infty,2)}(a, b, c, g; x)$ & $(\infty, -2a-2b-2c-g-2)$\\
  \end{tabular}
  \caption{Names of local solutions of D-F equations corresponding to those of $Q(A)$}
  \label{correspond_solution_Q_DF}
\end{center}
\end{table}

\subsection{Local solutions of the Dotsenko-Fateev equation at $x=0$}
\begin{prp}
$(1)$ 
%\begin{prp}
\begin{align*}
 &f_{DF}^{(0,0)}(a,b,c,g;x)/ (1-x)^{b+c+1}
 \\
 &=
% (1-x)^{b+c+1}
    \sum_{n=0}^{\infty}
      x^{n}
  \dfrac{
    \left( -a-c-g \right)_{n}
  }{
    \left( -a-c-\frac{g}{2}\right)_{n}
  }
  [Z^{n}]
  \pFq{2}{1}{-c,-a-b-c-\frac{g}{2}-1}{-a-c}{Z}
  \pFq{2}{1}{b+1,-a-\frac{g}{2}}{-a-c-g}{Z}
  \\
 &=
% (1-x)^{b+c+1}
    \sum_{n=0}^{\infty}
      x^{n}
  \dfrac{
    (-2a-b-c-g-1)_{n}
    (b-c+1)_{n}
  }{
      \left(- a-c \right)_{n}
      \left(- a-c-\frac{g}{2} \right)_{n}
  }
  [Z^{n}]
  \pFq{2}{1}{-c,-c-\frac{g}{2}}{b-c+1}{Z}
  \pFq{2}{1}{-a,-a-\frac{g}{2}}{-2a-b-c-g-1}{Z}
  \\
 &=
% (1-x)^{b+c+1}
 \sum_{n=0}^{\infty}
      x^{n}
  \dfrac{
    \left( -a-c-\frac{g}{2}\right)_{n}
  }{
    \left( -a-c \right)_{n}
  }
  [Z^{n}]
  \pFq{2}{1}{-c,-a-b-c-g-1}{-a-c-\frac{g}{2}}{Z}
  \pFq{2}{1}{b+1,-a}{-a-c-\frac{g}{2}}{Z}.
\end{align*}
\label{series_sol_df_0_0}
%\end{prp}
$(2)$ 
\begin{align*}
 &f_{DF}^{(0,1)}(a,b,c,g;x)/x^{a+c+1}(1-x)^{b+c+1}
 \\
 &=
 %x^{a+c}(1-x)^{b+c+1}
    \sum_{n=0}^{\infty}
      x^{n}
  \dfrac{
    \left( 1-g \right)_{n}
  }{
    \left( 1-\frac{g}{2}\right)_{n}
  }
  [Z^{n}]
  \pFq{2}{1}{b+1,-a-\frac{g}{2}}{-a-c-g}{Z}
  \pFq{2}{1}{a+1,-b-\frac{g}{2}}{a+c+2}{Z}
  \\
 &=
 %x^{a+c}(1-x)^{b+c+1}
    \sum_{n=0}^{\infty}
      x^{n}
  \dfrac{
    (1-g)_{n}
    (a+b+2)_{n}
    (-a-b-g)_{n}
  }{
      \left(1- \frac{g}{2} \right)_{n}
      \left(- a-c-g \right)_{n}
    \left( a+c+2 \right)_{n}
  }
  [Z^{n}]
  \pFq{2}{1}{c+1,-a-\frac{g}{2}}{-a-b-g}{Z}
  \pFq{2}{1}{a+1,-c-\frac{g}{2}}{a+b+2}{Z}
  \\
 &=
% x^{a+c}(1-x)^{b+c+1}
 \sum_{n=0}^{\infty}
      x^{n}
  \dfrac{
    \left( 1-\frac{g}{2}\right)_{n}
    \left( 1-g \right)_{n}
  }{
    \left( -a-c-g \right)_{n}\left( a+c+2 \right)_{n}
  }
  [Z^{n}]
  \pFq{2}{1}{a+1,-a-b-c-g-1}{1-\frac{g}{2}}{Z}
  \pFq{2}{1}{b+1,c+1}{1-\frac{g}{2}}{Z}.
\end{align*}
\label{series_sol_df_0_a+c}
%\end{prp}
%
$(3)$ 
%\begin{prp}
  \begin{align*}
    f_{DF}^{(0,2)}(a,b,c,g;x)
    =
    x^{2a+2c+g+2}f_{DF}^{(0,0)}
    \left(-c-\frac{g}{2}-1,a+b+c+\frac{g}{2}+1,-a-\frac{g}{2}-1,g; x\right).
  \end{align*}
  \label{expression_df_0_2a+2c+g}
\end{prp}
Proof:\ By the relation (\ref{relation_q_df}) and Proposition \ref{expression_q_0_-}, we get (1) and (3).
From the expressions
(\ref{expression_q_0_0_1_another}), 
(\ref{expression_q_0_0_2_another}), 
(\ref{expression_q_0_0_3_another})
of the holomorphic solution $f_Q^{(0,0)}(A;x)$ of $Q(A)$, we have (2).
%By the relation (\ref{relation_q_df}) and Proposition \ref{expression_q_0_-}, we get (3).
\qed

\begin{rmk}
By using the trivial symmetry \eqref{solution_symmetryQ_trivial} of the equation $Q(A)$, we can obtain other expressions of   $f_{DF}^{(0,i)}(a,b,c,g;x)$ $(i=0,1,2)$, which we omit. 
\end{rmk}
Relation (\ref{relation_q_df}) between $Q(A)$ and D-F, 
and the Pfaff transformation (\ref{solution_symmetryQ_pfaff}) 
lead to other expressions of $f_{DF}^{(0,0)}$ and $f_{DF}^{(0,1)}$:

\begin{prp}
\phantom{}\\
$(1)$\
$\displaystyle
f_{DF}^{(0,0)}(a,b,c,g;x)
=(1-x)^{2c}
f_{DF}^{(0,0)}\left(a,-a-b-c-g-2,c,g;\dfrac{x}{x-1}\right).$
\\
$(2)$\
$\displaystyle
f_{DF}^{(0,1)}(a,b,c,g;x)
=C\times  (1-x)^{2c}
f_{DF}^{(0,1)}\left(a,-a-b-c-g-2,c,g;\dfrac{x}{x-1}\right),
$
\\
where $C$ is a constant depending on the choice of a branch around $x=0$;
or more precisely, 
\begin{align*}
 &f_{DF}^{(0,1)}(a,b,c,g;x)/ x^{a+c+1}(1-x)^{b+c+g}
 \\
 &=
% x^{a+c}(1-x)^{b+c+g}
    \sum_{n=0}^{\infty}
    \left(\dfrac{x}{x-1}\right)^{n}
  \dfrac{
    \left( 1-g \right)_{n}
  }{
    \left( 1-\frac{g}{2}\right)_{n}
  }
  [Z^{n}]
  \pFq{2}{1}{c+1,-b-\frac{g}{2}}{a+c+2}{Z}
  \pFq{2}{1}{b+1,-c-\frac{g}{2}}{-a-c-g}{Z}
  \\
 &=
% x^{a+c}(1-x)^{b+c+g}
    \sum_{n=0}^{\infty}
    \left(\dfrac{x}{x-1}\right)^{n}
  \dfrac{
    \left( 1-g \right)_{n}
    \left( b+c+2 \right)_{n}
    \left( -b-c-g \right)_{n}
  }{
    \left( 1-\frac{g}{2} \right)_{n}
    \left( a+c+2 \right)_{n}
    \left( -a-c-g\right)_{n}
  }
  [Z^{n}]
  \pFq{2}{1}{c+1,-a-\frac{g}{2}}{b+c+2}{Z}
  \pFq{2}{1}{a+1,-c-\frac{g}{2}}{-b-c-g}{Z}
  \\
 &=
% x^{a+c}(1-x)^{b+c+g}
    \sum_{n=0}^{\infty}
    \left(\dfrac{x}{x-1}\right)^{n}
  \dfrac{
    \left( 1-\frac{g}{2} \right)_{n}
    \left( 1-g \right)_{n}
  }{
    \left( a+c+2 \right)_{n}
    \left( -a-c-g \right)_{n}
  }
  [Z^{n}]
  \pFq{2}{1}{a+1,b+1}{1-\frac{g}{2}}{Z}
  \pFq{2}{1}{c+1,-a-b-c-g-1}{1-\frac{g}{2}}{Z}.
\end{align*}
\label{series_sol_df_0_0_pfaff}
\end{prp}
\subsection{Local solutions of the  Dotsenko-Fateev equation at $x=1$}
By the symmetry of (\ref{solution_symmetryQ_1}) of $Q(A)$
and the relation (\ref{relation_q_df}) between $Q(A)$ and D-F,
we have
\begin{prp}
\phantom{}\\
$(1)$\
$\displaystyle
f_{DF}^{(1,0)}\left(a,b,c,g;x\right)
=f_{DF}^{(0,0)}\left(b,a,c,g;1-x\right).
$\\[3pt]
$(2)$\
$\displaystyle
f_{DF}^{(1,1)}\left(a,b,c,g;x\right)
=f_{DF}^{(0,1)}\left(b,a,c,g;1-x\right).
$\\[3pt]
$(3)$\
$\displaystyle
f_{DF}^{(1,2)}\left(a,b,c,g;x\right)
=f_{DF}^{(0,2)}\left(b,a,c,g;1-x\right)  \\[2pt]
  \phantom{(3)\ f_{DF}^{(1,2)}\left(a,b,c,g;x\right)}
  =  (1-x)^{2b+2c+g+2}
  f_{DF}^{(0,0)}\left(-c-\frac{g}{2}-1,a+b+c+\frac{g}{2}+1,-b-\frac{g}{2}-1,g;1-x\right).
$
\end{prp}

\subsection{Solutions of the Dotsenko-Fateev equation at $x=\infty$}
The relation (\ref{solution_symmetryQ_infty})
and the relation (\ref{relation_q_df}) between $Q(A)$ and D-F lead to
\begin{prp}
\phantom{}\\
$(1)$\
$\displaystyle
f_{DF}^{(\infty,0)}\left(a,b,c,g;x\right)
=\left(-\dfrac{1}{x}\right)^{-2c}
f_{DF}^{(0,0)}\left(-a-b-c-g-2,b,c,g;\dfrac{1}{x}\right).
$\\[2pt]
$(2)$\
$\displaystyle
f_{DF}^{(\infty,1)}\left(a,b,c,g;x\right)
=\left(-\dfrac{1}{x}\right)^{-2a-2b-2c-g-2}
f_{DF}^{\left( 0,1 \right)}\left(a,c,b,g;\dfrac{1}{x}\right).
$\\[2pt]
$(3)$\
$\displaystyle
f_{DF}^{(\infty,2)}\left(a,b,c,g;x\right)
=\left(-\dfrac{1}{x}\right)^{-2a-2b-2c-g-2}
  f_{DF}^{\left( 0,0 \right)}\left(a,c,b,g;\dfrac{1}{x}\right).$
\end{prp}

\par\bigskip\noindent
    {Acknowledgement:} The authors thank D. Zagier for presenting them with the challenging system $Z_3(A),a_1=a_2=a_3=0$. They also thank N. Takayama for instructing them about computer systems as well as various computational skills, and K. Mimachi for telling them his results on the Dotsenko-Fateev equation.
\par
%%%%%%%%%%%% Authors addresses %%%%%%%%%%%%%
%%%%%\newpage%
\par\bigskip
\small
\noindent
\author{Akihito Ebisu}\\
Faculty of Information and Computer Science,\
Chiba Institute of Technology,\\
Chiba 275-0023, Japan\\
email: {akihito.ebisu@p.chibakoudai.jp}
\par\vskip2ex\noindent
\author{Yoshishige Haraoka}\\
Department of Mathematics,\
Kumamoto University,\\
Kumamoto 860-8555, Japan\\
email: {haraoka@kumamoto-u.ac.jp}
\par\vskip2ex\noindent
\author{Masanobu Kaneko}\\
   Department of Mathematics,\
   Kyushu University,\\
   Fukuoka 819-0395, Japan\\
email: {mkaneko@math.kyushu-u.ac.jp}
\par\vskip2ex\noindent
\author{Hiroyuki Ochiai}\\
   Department of Mathematics,\
   Kyushu University,\\
   Fukuoka 819-0395, Japan\\
email: {ochiai@imi.kyushu-u.ac.jp}
\par\vskip2ex\noindent
\author{Takeshi Sasaki}\\
   Kobe University,\\
   Kobe 657-8501, Japan\\
email: {sasaki@math.kobe-u.ac.jp}
\par\vskip2ex\noindent
\author{Masaaki Yoshida}\\
   Kyushu University,\\
   Fukuoka 819-0395, Japan\\
email: {myoshida@math.kyushu-u.ac.jp}
\end{document}